\documentclass[opre,nonblindrev]{informs4_hide}

\OneAndAHalfSpacedXI

\usepackage{amsmath,amsfonts,amssymb}
\usepackage{mathtools}
\allowdisplaybreaks
\usepackage{bm}
\usepackage[mathscr]{euscript}
\usepackage{enumitem}
\usepackage{xcolor}
\usepackage{graphicx}
\usepackage{multirow}
\usepackage{booktabs}
\usepackage{microtype}
\usepackage{fix-cm}
\usepackage[ruled,vlined]{algorithm2e}
\usepackage{physics}
\usepackage{xfrac}
\usepackage{nicematrix}
\usepackage{float}
\usepackage{url}
\usepackage{pgfplots}
\pgfplotsset{compat=1.16}

\DeclareMathOperator{\pr}{\mathbb P}
\DeclareMathOperator{\E}{\mathbb E}
\DeclareMathOperator{\Var}{\mathrm{Var}}

\DeclareMathOperator{\ind}{\mathbb I}

\newcommand{\eps}{\varepsilon}
\newcommand{\R}{\mathbb R}

\newcommand{\calG}{\mathcal G}

\newcommand{\wt}{\widetilde}

\usepackage{natbib}
\usepackage{bibunits}
\bibpunct[, ]{(}{)}{,}{a}{}{,}%
\def\bibfont{\normalsize}

\TheoremsNumberedThrough
\ECRepeatTheorems
\EquationsNumberedThrough
\MANUSCRIPTNO{}

\begin{document}

\RUNAUTHOR{Huang et al.}
\RUNTITLE{Two-Scale Diffusion Approximations for Many-Server Queues with Overdispersed Arrivals}
\TITLE{Two-Scale Diffusion Approximations for Many-Server Queues with Overdispersed Arrivals}

\ARTICLEAUTHORS{
\AUTHOR{Junfei Huang}
\AFF{Department of Decisions, Operations and Technology, CUHK Business School, Chinese University of Hong Kong, Shatin, Hong Kong SAR,
\EMAIL{junfeih@cuhk.edu.hk}}
\AUTHOR{Jinghan Ji, Xiaowei Zhang}
\AFF{Department of Industrial Engineering and Decision Analytics, The Hong Kong University of Science and Technology, Clear Water Bay, Hong Kong SAR,
\EMAIL{alvin.ji@connect.ust.hk}, \EMAIL{xiaoweiz@ust.hk}}
}

\ABSTRACT{We study steady-state diffusion approximations for a many-server queue driven by overdispersed arrivals from a doubly stochastic Poisson process. Its arrival intensity follows a generalized Cox--Ingersoll--Ross (CIR) process, where $\alpha \in [0, 1)$ governs how the arrival-count variance scales as $\lambda^{1+\alpha}$ with the baseline arrival rate $\lambda$ under Taylor's law. For $\alpha > 0$, arrival-intensity fluctuations grow faster with $\lambda$ than counting fluctuations from arrivals and departures, so normalizing by $\lambda^{(1+\alpha)/2}$ pairs a nonvanishing intensity volatility with a vanishing queue-length volatility. This two-scale formulation differs critically from conventional diffusion approximations, which either apply a single $\lambda^{1/2}$ scaling with nonvanishing volatilities or drop the vanishing queue-length noise in the heavy-traffic limit. Retaining both scales, we compare two diffusions: one replaces the scaled intensity with its limiting Ornstein--Uhlenbeck (OU) process, and the other preserves the exact CIR intensity. Using Stein's method, we establish steady-state error bounds for smooth test functions of arbitrary polynomial growth that depend on $\alpha$ in opposite directions. Both diffusions achieve the classical $O(\lambda^{-1/2})$ rate at $\alpha = 0$. For $\alpha \in (0, 1)$, however, the OU rate degrades to the sharp order $O(\lambda^{-(1-\alpha)/2})$, whereas preserving the CIR intensity improves the rate to $O(\lambda^{-(2+\alpha)/4})$. Our proofs establish new derivative bounds for the associated Poisson equations that remain valid when the diffusion has vanishing volatility and nonsmooth drift.
}


\maketitle

\section{Introduction}
\label{sec:intro}

Managing demand uncertainty is central to capacity planning in service systems, where service resources represent both the main operating expense and performance bottleneck. Empirical studies across modern service sectors, including contact centers, healthcare facilities, cloud platforms, and digital communication networks, show that real-world arrivals of customers, patients, and requests follow complex stochastic dynamics \citep{Oreshkin2016, Hu2025, Daw2025Contact, Daw2025Batch}. Across these settings, arrivals consistently display two critical features that complicate capacity planning at large scale. First, arrival counts over fixed time intervals are \emph{overdispersed} and follow \emph{Taylor's law} \citep{Taylor2019, Hong2023}: their variance exceeds the mean and grows as a power function of that mean. Second, arrivals are \emph{temporally dependent}: periods of high or low demand persist across consecutive intervals before gradually returning toward the long-run average \citep{DawPender2018, DawPender2022, Hong2026}. When large variance and temporal persistence occur together, periods of elevated demand are both larger and longer-lasting than independent-arrival models predict, substantially increasing congestion and the capacity needed to maintain service quality.

To capture both features within a unified arrival process, \citet{Hong2026} recently proposed a doubly stochastic Poisson process (DSPP) driven by a generalized Cox--Ingersoll--Ross (CIR) arrival intensity. In this model, the stochastic intensity reverts exponentially toward a baseline arrival rate $\lambda$ to account for decaying temporal dependence, while a dispersion parameter $\alpha \in [0, 1)$ governs how the variance of arrival counts scales at rate $\lambda^{1+\alpha}$ under Taylor's law. This formulation covers both linear variance scaling ($\alpha = 0$) and power-law overdispersion ($\alpha > 0$) within one parsimonious model. Despite the practical prevalence of overdispersed and temporally dependent demand \citep{Hong2023}, however, queueing analysis for systems driven by such DSPP arrivals remains underexplored: existing work relies on heavy-traffic limits as $\lambda \to \infty$ \citep{SunLiu2021} or infinite-server approximations to guide capacity sizing \citep{Hong2026}, leaving pre-limit (finite-$\lambda$) performance evaluation for finite-server queues open.

Evaluating and optimizing steady-state performance in these finite-server systems requires approximations that are both analytically transparent and scale-free as $\lambda$ grows. Unlike infinite-server or equal-rate baselines, the exact joint process generally has no closed-form stationary distribution, obscuring how overdispersion, temporal persistence, and safety capacity jointly govern congestion. Moreover, direct numerical evaluation of the exact system scales poorly with $\lambda$: discrete-event simulation requires long run lengths alongside $O(\lambda)$ event rates per unit time and suffers from superlinear variance growth, while direct Markovian discretization must resolve $O(\lambda^{(1+\alpha)/2})$ discrete queue states coupled with a continuous arrival intensity and break down once arrival or service dynamics depart from tractable algebraic structure.

Diffusion approximations offer a natural analytical alternative, yet constructing an accurate diffusion model at finite $\lambda$ presents a multi-scale problem when $\alpha > 0$. Congestion in a DSPP-driven queue is governed by two sources of randomness that scale differently with $\lambda$: conditional on the current arrival and departure rates, \emph{counting fluctuations} from individual arrivals, service completions, and abandonments accumulate over $O(\lambda)$ events per unit time and have standard deviation of order $\lambda^{1/2}$, whereas \emph{intensity fluctuations} under Taylor's law have standard deviation of order $\lambda^{(1+\alpha)/2}$. Following \citet{Hong2026}, normalizing the centered arrival intensity and arrival process by $\lambda^{(1+\alpha)/2}$ keeps intensity fluctuations of order one while, for $\alpha > 0$, causing the counting fluctuations to vanish and the volatility of the scaled intensity to approach a constant as $\lambda$ increases.

Under this scaling, the scaled intensity converges to an Ornstein--Uhlenbeck (OU) process and the scaled arrival process converges to an integrated OU process, a Gaussian limit that offers high analytical tractability. Passing directly to this $\lambda \to \infty$ limit in the queueing dynamics, however, treats the queue length as deterministic once the arrival-intensity path is given. To retain both sources of randomness at finite $\lambda$, a natural starting point is a \emph{two-scale} OU-based diffusion approximation: a coupled two-dimensional stochastic differential equation (SDE) that pairs the OU arrival-intensity equation with a queue-length equation whose $O(\lambda^{-\alpha/2})$ volatility matches the stationary variance rate of counting fluctuations.

We apply Stein's method \citep{Gurvich2014, BravermanDai2017, HuangGurvich2018} to quantify the steady-state accuracy of this two-scale diffusion approximation. Comparing stationary generators through the associated Poisson equation avoids proving transient weak convergence and interchanging large-system and long-time limits, reducing the steady-state error to bounds on the partial derivatives of the Poisson-equation solution with respect to the queue length and the arrival intensity. Despite the constant volatility and linear drift of the OU intensity equation, two features of our two-dimensional diffusion place those derivative bounds outside existing Stein-method frameworks.

First, the volatility in the queue-length equation vanishes at rate $\lambda^{-\alpha/2}$ when $\alpha > 0$, whereas the volatility in the arrival-intensity equation remains of order one. Standard multidimensional derivative bounds scale inversely with the smallest volatility across both equations and therefore diverge as $\lambda$ grows, ignoring the fact that random fluctuations in the arrival intensity still propagate into the queue length through its drift.

Second, this vanishing volatility interacts with a nondifferentiable drift. Because the marginal departure rate per additional customer switches from $\mu$ (individual service rate) to $\theta$ (abandonment rate) once all servers become busy, the drift of the queue-length equation is nondifferentiable at the full-occupancy threshold unless $\mu = \theta$. Existing derivative bounds that track small volatilities require diffusions with twice-differentiable drifts and break down at this threshold, whereas existing Stein-method bounds for diffusions with piecewise-linear drifts assume that the volatility in every SDE equation remains bounded away from zero.

\subsection{Main Contributions}
\label{subsec:contributions}

First, we develop a pre-limit two-scale diffusion framework for DSPP-driven many-server queues across $\alpha \in [0, 1)$. Unlike traditional Stein-method analyses of Poisson or renewal queues (where a single $\lambda^{1/2}$ scaling keeps the volatility in every SDE equation of order one), our framework couples a variance-matched volatility of order $\lambda^{-\alpha/2}$ in the queue-length equation (which equals $\sqrt{2}$ at $\alpha = 0$ and vanishes for $\alpha > 0$) with an order-one volatility in the arrival-intensity equation. Once counting-fluctuation variance is matched, our stationary generator expansion reveals that steady-state accuracy is governed by the arrival-intensity dynamics, motivating a second two-scale diffusion that retains the exact CIR intensity alongside the baseline OU model.

Second, we show that the steady-state error bounds for the OU-based approximation motivated by \citet{Hong2026} and the CIR-preserving approximation depend on overdispersion in opposite directions. At $\alpha = 0$, both variance-matched diffusions recover the classical $O(\lambda^{-1/2})$ Erlang-A rate \citep{BravermanDaiFeng2017}. For $\alpha \in (0, 1)$, however, replacing the CIR intensity with an OU process loses the positive skewness of its Gamma stationary distribution, degrading the error bound to $O(\lambda^{-(1-\alpha)/2})$, a rate sharp over every stated smooth test class. Retaining the CIR intensity eliminates this skewness error and gives the bounds $O(\lambda^{-(2+\alpha)/4})$ when $\mu \neq \theta$ and $O(\lambda^{-(1+\alpha)/2})$ when $\mu = \theta$. Unlike at $\alpha = 0$, vanishing queue-length volatility separates the $\mu \neq \theta$ and $\mu = \theta$ rates when $\alpha > 0$. Moreover, when $\mu = \theta$, both variance-matched diffusions match all stationary joint moments of degree at most two for every admissible $\lambda$.

Third, to establish these guarantees under vanishing volatilities and a nondifferentiable drift, we develop new coordinate-specific and boundary-weighted derivative bounds for the Poisson equation of a two-dimensional diffusion. By combining short-time counting fluctuations with arrival-intensity fluctuations transmitted through the queue drift, we bound the intensity and mixed partial derivatives uniformly in the vanishing $O(\lambda^{-\alpha/2})$ queue-length volatility while controlling the time spent near full occupancy when $\mu \neq \theta$. Extending this analysis to the CIR-preserving diffusion introduces an additional state-dependent degeneracy near zero arrival intensity, where its square-root volatility approaches zero; we overcome this boundary degeneracy via a Lamperti pathwise sensitivity analysis, an adapted feedback Brownian perturbation that steers coupled variations to match initial queue-length derivatives, and uniform inverse-moment bounds for the Gamma stationary distribution of the arrival intensity.

\subsection{Related Works}
\label{subsec:literature}

\paragraph{Heavy-Traffic Limits and Pre-Limit Diffusion Approximations.}
Classical heavy-traffic theory establishes process-level and steady-state diffusion limits under Poisson or renewal arrivals \citep{HalfinWhitt1981, Garnett2002, Whitt2002, Whitt2004, GamarnikGoldberg2013, Dai2014}, while extensions to uncertain or stochastic arrival intensities similarly rely on asymptotic limits as $\lambda \to \infty$ \citep{Bassamboo2010, SunLiu2021}. To avoid sensitivity to asymptotic regime boundaries, a complementary literature constructs pre-limit diffusion approximations that retain system parameters directly \citep{WardGlynn2005, gurvich2014excursion} and uses Stein's method to prove steady-state error bounds across Markovian, phase-type, and renewal queues \citep{Gurvich2014, BravermanDaiFeng2017, BravermanDai2017, HuangGurvich2018, BravermanScully2026}. Whereas existing Stein-method guarantees rely on single-scale normalization with nonvanishing queue-length volatility, we establish pre-limit error bounds for two-scale DSPP-driven queues across dispersion regimes and show that replacing the arrival intensity with its heavy-traffic OU limit loses accuracy under strong overdispersion unless skewed intensity dynamics are preserved.

\paragraph{Poisson-Equation Derivative Bounds.}
Our proof machinery bridges two disjoint strands of the literature on Poisson equations for diffusions. Multivariate Stein bounds that allow small or degenerate volatilities assume smooth, contractive drifts \citep{Meckes2009, MackeyGorham2016, Gorham2019}, whereas Stein bounds for diffusions with piecewise-linear drifts at full occupancy require volatility bounded away from zero in every SDE equation \citep{Gurvich2014, BravermanDai2017, Jin2025}. We draw on drift-propagation ideas for degenerate SDEs \citep{DelarueMenozzi2010, Menozzi2011}, where randomness in a driving coordinate propagates into another through its drift. These results assume uniformly elliptic driving volatility and smooth coefficients; we extend this mechanism to simultaneously handle a nondifferentiable drift at full occupancy and a state-dependent CIR driving volatility that vanishes at the lower intensity boundary.

\subsection{Notation and Organization}

Throughout the paper, we use the following notation. For real $x, y$, write $x^+ = \max(x, 0)$ and $x \wedge y = \min(x, y)$, and let $\mathbb{I}_A$ (or $\mathbb{I}\{\cdot\}$) denote the indicator of an event or condition $A$. Subscripts on a function denote partial derivatives, and $\nabla h$ and $\nabla^2 h$ denote its gradient and Hessian. We use the Euclidean norm $|z| = (\sum_i z_i^2)^{1/2}$ for vectors and its induced operator norm $\|M\| = \sup_{|z|=1} |Mz|$ for matrices. We write $X(\infty)$ for a random variable with the stationary distribution of $X$, and $\pi(h) = \int h\,d\pi$ for the expectation of $h$ under a probability distribution $\pi$.

The remainder of the paper is organized as follows. Section~\ref{sec:model} formulates the queueing model and constructs the baseline OU-based diffusion. Section~\ref{sec:ou_results} establishes its steady-state error bound, proves rate sharpness, and presents the directional derivative bounds. Section~\ref{sec:cir_results} introduces the CIR-preserving diffusion, establishes its error bounds and the exact second-order moment agreement of both approximations when $\mu = \theta$, and outlines the boundary-weighted derivative bounds. Section~\ref{sec:discussion} discusses modeling implications, Section~\ref{sec:numerical-experiments} presents numerical experiments, and Section~\ref{sec:conclusion} concludes. All technical proofs appear in the e-companion.

\section{Problem Formulation}
\label{sec:model}

We consider a family of many-server queues indexed by the baseline arrival rate $\lambda \ge 1$. The $\lambda$th system has $n_\lambda$ identical servers; individual service times are independent exponential random variables with rate $\mu$, and customers' patience times are independent exponential random variables with abandonment rate $\theta$. Arrivals follow a doubly stochastic Poisson process (DSPP) $A_\lambda$ driven by a stochastic arrival intensity $I_\lambda$: conditional on the sample path $\{I_\lambda(s) : s \ge 0\}$, $A_\lambda$ is a nonhomogeneous Poisson process with cumulative rate $\int_0^t I_\lambda(s)\,ds$ over $[0, t]$. Following \citet{Hong2026}, $I_\lambda$ satisfies the generalized CIR stochastic differential equation (SDE)
\begin{equation}
dI_\lambda(t)
=
\kappa\bigl(\lambda - I_\lambda(t)\bigr)\,dt
+
\sigma\sqrt{\lambda^\alpha I_\lambda(t)}\,dB_2(t),
\label{eq:cir_sde}
\end{equation}
where $B_2$ is a standard Brownian motion, $\kappa$ is the mean-reversion rate, $\sigma$ governs the intensity volatility, and $\alpha \in [0, 1)$ is the dispersion parameter.

\begin{assumption}
\label{assump:model_primitives}
The service and abandonment rates satisfy $\mu > 0$ and $\theta > 0$. The arrival-intensity parameters satisfy $\kappa, \sigma > 0$ and $\alpha \in [0, 1)$. For the values of $\lambda$ under consideration, $2\kappa\lambda^{1-\alpha} \ge \sigma^2$, and $I_\lambda(0) > 0$ almost surely.
\end{assumption}

Under Assumption~\ref{assump:model_primitives}, the condition $2\kappa\lambda^{1-\alpha} \ge \sigma^2$ ensures that zero is an inaccessible lower boundary for $I_\lambda$, so $I_\lambda(t) > 0$ for all $t \ge 0$ almost surely. Moreover, $I_\lambda$ admits a unique Gamma stationary distribution with mean $\E[I_\lambda(\infty)] = \lambda$ and variance $\Var(I_\lambda(\infty)) = (\sigma^2/(2\kappa))\lambda^{1+\alpha}$.

Let $X_\lambda(t)$ denote the number of customers in the system at time $t$. When $X_\lambda(t) = x$, the total departure rate from service completions and abandonments is
\begin{equation}
d_\lambda(x)
=
\mu(x \wedge n_\lambda)
+
\theta(x - n_\lambda)^+,
\label{eq:death_rate}
\end{equation}
so the marginal departure rate per customer equals $\mu$ below full occupancy ($x < n_\lambda$) and switches to $\theta$ above full occupancy ($x > n_\lambda$). Let $N_a$ and $N_d$ be independent unit-rate Poisson processes, independent of $B_2$. Representing $X_\lambda$ via random time changes as
\[
X_\lambda(t)
=
X_\lambda(0)
+
N_a\!\left(\int_0^t I_\lambda(s)\,ds\right)
-
N_d\!\left(\int_0^t d_\lambda(X_\lambda(s))\,ds\right)
\]
yields the martingale decomposition and predictable quadratic variation
\begin{equation}
\begin{aligned}
X_\lambda(t)
&=
X_\lambda(0)
+
\int_0^t
\bigl[I_\lambda(s) - d_\lambda(X_\lambda(s))\bigr]\,ds
+
M_\lambda(t), \\
\langle M_\lambda \rangle_t
&=
\int_0^t
\bigl[I_\lambda(s) + d_\lambda(X_\lambda(s))\bigr]\,ds,
\end{aligned}
\label{eq:exact_martingale}
\end{equation}
where $M_\lambda$ is a zero-mean martingale capturing the counting fluctuations from individual arrivals, service completions, and abandonments.

\subsection{Two-Scale Normalization and Capacity Rule}
\label{subsec:queue-scales}

The model above differs from conventional queueing formulations in two structural ways that determine how $(X_\lambda, I_\lambda)$ must be scaled as $\lambda$ grows. First, under Poisson arrivals $X_\lambda(t)$ is Markovian on its own \citep{Gurvich2014, BravermanDaiFeng2017}, and under phase-type or renewal arrivals the auxiliary phase or residual interarrival variables fluctuate on the $O(\lambda^{-1})$ timescale between individual arrivals \citep{BravermanDai2017, BravermanScully2026}. Here, because $I_\lambda(t)$ reverts at an $O(1)$ rate $\kappa$ independent of $\lambda$, temporal dependence persists over $O(1)$ time horizons and makes the Markovian state $(X_\lambda(t), I_\lambda(t))$ genuinely two-dimensional. Second, decomposing the centered arrival count over $[0, t]$ as
\[
A_\lambda(t) - \lambda t
=
\left[
A_\lambda(t) - \int_0^t I_\lambda(s)\,ds
\right]
+
\int_0^t \bigl[I_\lambda(s) - \lambda\bigr]\,ds
\]
separates two uncorrelated stationary variance components of different orders: the first term measures conditional counting fluctuations and has stationary variance $\lambda t$, whereas the second measures integrated arrival-intensity fluctuations and has stationary variance
\[
\Var\!\left(\int_0^t \bigl[I_\lambda(s) - \lambda\bigr]\,ds\right)
=
\frac{\sigma^2}{\kappa^3}\lambda^{1+\alpha}
\bigl(\kappa t - 1 + e^{-\kappa t}\bigr),
\]
so $\Var(A_\lambda(t))$ grows at order $\lambda^{1+\alpha}$ under Taylor's law for every fixed $t > 0$. Whereas Poisson and renewal arrival counts have variance of order $\lambda$ and normalize on the single scale $\lambda^{1/2}$, here counting fluctuations of standard deviation $O(\lambda^{1/2})$ coexist with arrival-intensity fluctuations of standard deviation $O(\lambda^{(1+\alpha)/2})$, which dominate whenever $\alpha \in (0, 1)$.

To keep the $O(\lambda^{(1+\alpha)/2})$ arrival-intensity fluctuations of order one as $\lambda$ grows, we introduce the scaling parameters and centered, scaled state
\begin{equation}
\begin{aligned}
\delta_\lambda &= \lambda^{-(1+\alpha)/2}, &
\eps_\lambda &= \lambda^{-(1-\alpha)/2}, \\
\wt X_\lambda(t) &= \delta_\lambda\bigl(X_\lambda(t) - n_\lambda\bigr), \qquad &
\wt I_\lambda(t) &= \delta_\lambda\bigl(I_\lambda(t) - \lambda\bigr),
\end{aligned}
\label{eq:two_scale_defs}
\end{equation}
where $\delta_\lambda \eps_\lambda = \lambda^{-1}$ and $\lambda \delta_\lambda^2 = \lambda^{-\alpha}$. Centering $X_\lambda(t)$ at $n_\lambda$ places the full-occupancy threshold at $\wt X_\lambda(t) = 0$. Because $I_\lambda(t) > 0$ under Assumption~\ref{assump:model_primitives}, $(\wt X_\lambda(t), \wt I_\lambda(t))$ takes values in $\mathsf{D}_\lambda = \R \times (-\eps_\lambda^{-1}, \infty)$, where the inaccessible lower boundary $u = -\eps_\lambda^{-1}$ corresponds to $I_\lambda = 0$.

\begin{assumption}
\label{assump:staffing_c}
The number of servers $n_\lambda$ satisfies
\begin{equation}
\beta_\lambda = \delta_\lambda(n_\lambda \mu - \lambda),
\qquad
\sup_{\lambda \ge 1} |\beta_\lambda| < \infty,
\qquad
\beta_\lambda \to \beta \in \R.
\label{eq:staffing}
\end{equation}
Moreover, the volatility parameter $c_\lambda$ in the queue-length SDE introduced in Section~\ref{subsec:diffusion-construction} satisfies $0 < c_\lambda \le \bar{c}$ for a fixed finite constant $\bar{c} \ge \sqrt{2}$.
\end{assumption}

Equivalently, $n_\lambda = \lambda/\mu + \beta_\lambda \lambda^{(1+\alpha)/2}/\mu$, so $\beta_\lambda$ measures safety capacity above the offered load $\lambda/\mu$ in units of $\lambda^{(1+\alpha)/2}/\mu$ servers. This capacity rule recovers classical $\lambda^{1/2}$ safety capacity at $\alpha = 0$ \citep{HalfinWhitt1981, Garnett2002} and scales safety capacity with the $O(\lambda^{(1+\alpha)/2})$ standard deviation of the arrival intensity when $\alpha \in (0, 1)$ \citep{Hong2026}.

\subsection{OU-Based Diffusion Approximation}
\label{subsec:diffusion-construction}

Rewriting the exact system in the scaled coordinates $(\wt X_\lambda, \wt I_\lambda)$ motivates our diffusion formulation. Define the piecewise-linear scaled departure-drift function, the scaled CIR volatility factor, and the parameter bounds
\begin{equation}
\begin{gathered}
g(y) = \mu(y \wedge 0) + \theta y^+, \qquad
v_\lambda(u) = 1 + \eps_\lambda u, \\
a_0 = \min\{\mu, \theta\}, \qquad
a_1 = \max\{\mu, \theta\}, \qquad
B_\beta = \sup_{\lambda \ge 1} |\beta_\lambda| < \infty.
\end{gathered}
\label{eq:queue-drift-abbreviations}
\end{equation}
Using \eqref{eq:cir_sde}--\eqref{eq:queue-drift-abbreviations} and $\delta_\lambda \eps_\lambda = \lambda^{-1}$, the exact scaled pair $(\wt X_\lambda, \wt I_\lambda)$ satisfies
\begin{equation}
\begin{aligned}
\wt X_\lambda(t)
&=
\wt X_\lambda(0)
+
\int_0^t
\bigl[
-\beta_\lambda + \wt I_\lambda(s)
- g(\wt X_\lambda(s))
\bigr]\,ds
+
\delta_\lambda M_\lambda(t), \\
d\wt I_\lambda(t)
&=
-\kappa \wt I_\lambda(t)\,dt
+
\sigma\sqrt{v_\lambda(\wt I_\lambda(t))}\,dB_2(t),
\end{aligned}
\label{eq:exact_scaled_system}
\end{equation}
where the scaled counting-fluctuation martingale has predictable quadratic variation
\[
\langle \delta_\lambda M_\lambda \rangle_t
=
\int_0^t
\Bigl[
2\lambda^{-\alpha}
+ \delta_\lambda\bigl\{
\wt I_\lambda(s) + \beta_\lambda
+ g(\wt X_\lambda(s))
\bigr\}
\Bigr]\,ds.
\]

A third structural difference from single-scale queueing models appears when passing to the heavy-traffic limit $\lambda \to \infty$. Under $\lambda^{1/2}$ scaling ($\alpha = 0$), $\langle \delta_\lambda M_\lambda \rangle_t$ remains of order one as $\lambda \to \infty$, so heavy-traffic diffusion limits and pre-limit diffusion approximations share an $O(1)$ volatility in the queue-length equation \citep{Gurvich2014, BravermanDaiFeng2017, HuangGurvich2018}. When $\alpha \in (0, 1)$, however, $\eps_\lambda \to 0$ and $\wt I_\lambda$ converges weakly to an OU process $U$ \citep{Hong2026}, while $\delta_\lambda M_\lambda(t)$ vanishes in probability and collapses the limiting queue-length dynamics to the conditional ordinary differential equation $\dot{Y}(t) = -\beta + U(t) - g(Y(t))$. At any finite $\lambda$, by contrast, discrete arrivals and departures still drive random excursions across full occupancy ($y = 0$), where the derivative of $g(y)$ jumps from $\mu$ to $\theta$ whenever $\mu \neq \theta$.

To retain both sources of randomness at finite $\lambda$, we pair the limiting OU intensity $U$ with a Brownian counting-fluctuation term $c_\lambda B_1(t)$, where $B_1$ is a standard Brownian motion independent of $B_2$. This yields our baseline \emph{OU-based two-scale diffusion} $(Y_\lambda^{\mathrm{OU}}, U)$ on $\R^2$:
\begin{equation}
\begin{aligned}
dY_\lambda^{\mathrm{OU}}(t)
&=
\bigl[
-\beta_\lambda + U(t)
- g(Y_\lambda^{\mathrm{OU}}(t))
\bigr]\,dt
+ c_\lambda\,dB_1(t), \\
dU(t)
&=
-\kappa U(t)\,dt + \sigma\,dB_2(t).
\end{aligned}
\label{eq:ou_diffusion_general_c}
\end{equation}
Conservation of flow in stationarity gives $\E[d_\lambda(X_\lambda(\infty))] = \E[I_\lambda(\infty)] = \lambda$. Matching $\langle c_\lambda B_1 \rangle_t = c_\lambda^2 t$ to the resulting stationary counting-variance rate yields the \emph{variance-matched volatility}:
\begin{equation}
\E\langle \delta_\lambda M_\lambda \rangle_t
=
2\lambda \delta_\lambda^2 t
=
2\lambda^{-\alpha} t,
\qquad
c_\lambda^*
=
\sqrt{2}\,\lambda^{-\alpha/2},
\label{eq:variance_matched_c}
\end{equation}
which unifies the $\alpha = 0$ regime ($c_\lambda^* = \sqrt{2}$) and the overdispersed regime $\alpha \in (0, 1)$ ($c_\lambda^* \to 0$). Meanwhile, replacing $v_\lambda(u) = 1 + \eps_\lambda u$ with $v_\lambda(0) = 1$ in the intensity SDE preserves the exact first- and second-order stationary structure of $\wt I_\lambda$:
\[
\E[\wt I_\lambda(\infty)] = \E[U(\infty)] = 0,
\qquad
\operatorname{Cov}\bigl(\wt I_\lambda(0), \wt I_\lambda(t)\bigr)
=
\operatorname{Cov}\bigl(U(0), U(t)\bigr)
=
\frac{\sigma^2}{2\kappa}e^{-\kappa t}
\quad (t \ge 0).
\]
Unlike single-scale diffusion approximations that only replace discrete counting jumps with Brownian motion, \eqref{eq:ou_diffusion_general_c} simultaneously approximates counting fluctuations $\delta_\lambda M_\lambda$ via $c_\lambda B_1$ and replaces the Gamma stationary distribution of $I_\lambda(\infty)$ (after centering and scaling) with the symmetric Gaussian stationary distribution of $U$.

Under Assumptions~\ref{assump:model_primitives} and~\ref{assump:staffing_c}, both $(\wt X_\lambda, \wt I_\lambda)$ and $(Y_\lambda^{\mathrm{OU}}, U)$ admit unique stationary distributions, denoted by $\pi^X$ and $\pi^{\mathrm{OU}}$, whose polynomial moments of every fixed order are uniformly bounded for sufficiently large $\lambda$ (Lemmas~\ref{lem:diffusion_drift} and~\ref{lem:exact_scaled_moments}). To quantify steady-state approximation accuracy across joint queue-length and arrival-intensity observables, write $w(y, u) = 1 + |y| + |u|$, fix an integer $m \ge 2$ and a constant $H > 0$, and consider test functions $h \in C^2(\R^2)$ satisfying
\begin{equation}
|h(y, u)| \le H w(y, u)^m,
\qquad
|\nabla h(y, u)| \le H w(y, u)^{m-1},
\qquad
\|\nabla^2 h(y, u)\| \le H w(y, u)^{m-2}.
\label{eq:test_class}
\end{equation}
Every polynomial of total degree at most $m$ belongs to \eqref{eq:test_class} for a suitable $H$: taking $m = 2$ covers stationary means, variances, and queue--intensity covariances, whereas $m \ge 3$ admits higher-order moments sensitive to distributional skewness. The class also contains nonpolynomial functions: $h(y,u)=H\sin u$ satisfies \eqref{eq:test_class} for every $m\ge2$.

Our goal in the sequel is to characterize the pre-limit steady-state error $|\pi^X(h) - \pi^D(h)|$ across $\alpha \in [0, 1)$ for test functions satisfying \eqref{eq:test_class}, beginning with $D = \mathrm{OU}$ in Section~\ref{sec:ou_results}: how should $c_\lambda \in (0, \bar{c}]$ be calibrated, at what rate does the variance-matched error decay with $\lambda$, and can accuracy under strong overdispersion be improved in Section~\ref{sec:cir_results} by retaining the exact CIR arrival-intensity dynamics?

\section{Analysis of the OU-Based Diffusion Approximation}
\label{sec:ou_results}

\subsection{Steady-State Error Bounds and Variance Matching}
\label{subsec:ou-error-bounds}

The OU-based diffusion $(Y_\lambda^{\mathrm{OU}}, U)$ in \eqref{eq:ou_diffusion_general_c} approximates discrete counting fluctuations by the Brownian term $c_\lambda B_1$ and replaces the state-dependent CIR arrival-intensity volatility with the constant $\sigma$. To isolate the role of variance matching, we first establish a steady-state error bound for an arbitrary queue-length volatility $c_\lambda \in (0, \bar{c}]$ and then specialize to the variance-matched choice $c_\lambda^* = \sqrt{2}\,\lambda^{-\alpha/2}$ in \eqref{eq:variance_matched_c}.

\begin{theorem}
\label{thm:smooth_bound_general_c}
Suppose Assumptions~\ref{assump:model_primitives} and~\ref{assump:staffing_c} hold and $0 < c_\lambda \le \bar{c}$. Let $h$ satisfy \eqref{eq:test_class} for fixed $m \ge 2$ and $H > 0$. There exists a constant $C < \infty$, depending only on the fixed model parameters, $\bar{c}$, $m$, and $H$, such that, for all sufficiently large $\lambda$,
\begin{equation}
\left|\pi^X(h) - \pi^{\mathrm{OU}}(h)\right|
\le
C
\begin{cases}
\displaystyle
\left|\lambda^{-\alpha} - \frac{c_\lambda^2}{2}\right|
+ \lambda^{-(1-\alpha)/2},
& \mu = \theta, \\[3mm]
\displaystyle
c_\lambda^{-1/2}
\left|\lambda^{-\alpha} - \frac{c_\lambda^2}{2}\right|
+ \lambda^{-(1-\alpha)/2},
& \mu \ne \theta.
\end{cases}
\label{eq:ou_bound_second_order_only}
\end{equation}
Under variance matching, $c_\lambda = c_\lambda^* = \sqrt{2}\,\lambda^{-\alpha/2}$, with $\delta_\lambda$ and $\eps_\lambda$ as in \eqref{eq:two_scale_defs},
\begin{equation}
\left|\pi^X(h) - \pi^{\mathrm{OU}}(h)\right|
\le C
\begin{cases}
\lambda^{-(2+\alpha)/4} + \eps_\lambda, & \mu \ne \theta, \\
\delta_\lambda + \eps_\lambda, & \mu = \theta
\end{cases}
\ \le C'\eps_\lambda.
\label{eq:ou_bound_variance_matched}
\end{equation}
\end{theorem}

The intermediate decomposition in \eqref{eq:ou_bound_variance_matched} separates the error of approximating counting fluctuations from the error of replacing the CIR arrival intensity. Under variance matching, the counting-fluctuation contribution is $O(\delta_\lambda)$ when $\mu = \theta$ (where the queue-length drift is linear) and $O(\delta_\lambda (c_\lambda^*)^{-1/2}) = O(\lambda^{-(2+\alpha)/4})$ when $\mu \neq \theta$ (where nondifferentiability of $g(y)$ at full occupancy $y = 0$ introduces the factor $(c_\lambda^*)^{-1/2}$). By contrast, replacing the state-dependent CIR intensity volatility contributes $O(\eps_\lambda) = O(\lambda^{-(1-\alpha)/2})$. At $\alpha = 0$, all three terms reduce to the classical $O(\lambda^{-1/2})$ Erlang-A Stein rate \citep{Gurvich2014, BravermanDaiFeng2017}. Once $\alpha \in (0, 1)$, however, the two error sources move in opposite directions with $\alpha$: the counting-fluctuation error decays faster than $\lambda^{-1/2}$, whereas the arrival-intensity error $\eps_\lambda$ decays more slowly and dominates the overall $O(\lambda^{-(1-\alpha)/2})$ bound.

For a general volatility parameter $c \in (0, \bar{c}]$, write the counting-variance mismatch term in \eqref{eq:ou_bound_second_order_only} as
\begin{equation*}
\Psi_\lambda(c) = \left|\lambda^{-\alpha} - \frac{c^2}{2}\right|
\begin{cases}
c^{-1/2}, & \mu \ne \theta, \\
1, & \mu = \theta.
\end{cases}
\end{equation*}
Over $c \in (0, \bar{c}]$, $\Psi_\lambda(c)$ attains its unique minimum of zero at $c = c_\lambda^*$: choosing an order-one volatility when $\alpha > 0$ incurs an $O(1)$ variance mismatch, while letting $c_\lambda$ shrink faster than $\lambda^{-\alpha/2}$ when $\mu \neq \theta$ is penalized by $c_\lambda^{-1/2}$. Applying Theorem~\ref{thm:smooth_bound_general_c} to $h(y, u) = y$ and $h(y, u) = y^2$ yields error bounds for the stationary mean and variance of the scaled queue length.

\begin{corollary}
\label{cor:mean_variance}
Under the assumptions of Theorem~\ref{thm:smooth_bound_general_c}, there exists a constant $C < \infty$, independent of $\lambda$, such that, for all sufficiently large $\lambda$,
\begin{equation}
\left|
\E\wt X_\lambda(\infty)
- \E Y_\lambda^{\mathrm{OU}}(\infty)
\right|
+
\left|
\Var\bigl(\wt X_\lambda(\infty)\bigr)
- \Var\bigl(Y_\lambda^{\mathrm{OU}}(\infty)\bigr)
\right|
\le
C\bigl\{\Psi_\lambda(c_\lambda) + \eps_\lambda\bigr\}.
\label{eq:ou_queue_moment_bound}
\end{equation}
In particular, when $c_\lambda = c_\lambda^*$, the right-hand side reduces to $C\eps_\lambda = C\lambda^{-(1-\alpha)/2}$.
\end{corollary}

Because $X_\lambda(\infty) = n_\lambda + \delta_\lambda^{-1}\wt X_\lambda(\infty)$, multiplying the variance-matched mean-error bound in Corollary~\ref{cor:mean_variance} by $\delta_\lambda^{-1}=\lambda^{(1+\alpha)/2}$ gives an $O(\lambda^\alpha)$ upper bound in customer-count units. When $\mu=\theta$, Proposition~\ref{prop:cir-exact-moments} sharpens this conclusion: both the mean and variance agree exactly under variance matching.

\subsection{Sharpness of the OU Bound}
\label{subsec:ou-sharpness}

In stationarity, the exact scaled intensity $\wt I_\lambda$ and the OU process $U$ share the same mean, variance, and autocovariance function, so first- and second-order intensity statistics cannot detect the replacement of $\wt I_\lambda$ by $U$. The first unequal polynomial moment is the third central moment: $I_\lambda(\infty)$ follows a positively skewed Gamma distribution, whereas $U(\infty)$ follows a symmetric Gaussian distribution with zero third central moment.

\begin{proposition}
\label{prop:ou-sharpness}
Under Assumptions~\ref{assump:model_primitives} and~\ref{assump:staffing_c}, the stationary error for $h(y, u) = u^3$ satisfies, for every admissible $c_\lambda$,
\begin{equation}
\left|\pi^X(h) - \pi^{\mathrm{OU}}(h)\right|
= \frac{\sigma^4}{2\kappa^2}\eps_\lambda
= \frac{\sigma^4}{2\kappa^2}\lambda^{-(1-\alpha)/2}.
\label{eq:discussion-ou-cubic-lower}
\end{equation}
The $O(\lambda^{-(1-\alpha)/2})$ rate in \eqref{eq:ou_bound_variance_matched} is sharp on every test class \eqref{eq:test_class} with $m \ge 2$ and $H > 0$.
\end{proposition}

As shown in Section~\ref{subsec:ou-moment-proofs}, the stationary Gamma distribution gives \eqref{eq:discussion-ou-cubic-lower}, and rescaling $h(y,u)=u^3$ gives an admissible test for every $m\ge3$ and $H>0$. For $m=2$, the bounded test $h(y,u)=H\sin u$ has a nonzero stationary error of order $\eps_\lambda$. Both tests depend only on the arrival intensity, so changing $c_\lambda$ cannot remove their error. Whether a smooth observable of the queue length alone attains the OU rate remains open.

\subsection{Generator Comparison and Directional Derivative Bounds}
\label{subsec:ou-proof-outline}

To prove Theorem~\ref{thm:smooth_bound_general_c}, let $(P_t)_{t \ge 0}$ denote the transition semigroup of $(Y_\lambda^{\mathrm{OU}}, U)$, and let $f$ solve the Poisson equation associated with $h$:
\begin{equation*}
\calG_{\mathrm{OU}}f = h - \pi^{\mathrm{OU}}(h),
\qquad
f(y, u) = -\int_0^\infty \bigl(P_t h(y, u) - \pi^{\mathrm{OU}}(h)\bigr)\,dt.
\end{equation*}
For sufficiently large $\lambda$, the stationary generator identity $\pi^X(\calG_X f) = 0$ (Lemma~\ref{lem:petite_generator_identity}) yields $\pi^X(h) - \pi^{\mathrm{OU}}(h) = \pi^X[(\calG_{\mathrm{OU}} - \calG_X)f]$. Expanding the unit-customer jumps of $\calG_X$ in Taylor series of step size $\pm\delta_\lambda$ at an exact scaled state $(y, u)$ gives
\begin{equation}
(\calG_X - \calG_{\mathrm{OU}})f(y, u)
={} R_{\mathrm{queue}}f(y, u)
+ \left(\lambda^{-\alpha} - \frac{c_\lambda^2}{2}\right)f_{yy}(y, u)
+ \frac{\sigma^2}{2}\eps_\lambda u f_{uu}(y, u).
\label{eq:master-comparison}
\end{equation}
Each term in \eqref{eq:master-comparison} corresponds to a distinct modeling step. First, $R_{\mathrm{queue}}f(y, u)$ collects the integral Taylor remainders from upward and downward jumps of size $\delta_\lambda$ (which depend on queue-direction increments $f_{yy}(y \pm s\delta_\lambda, u) - f_{yy}(y, u)$ for $s \in [0, 1]$) together with the state-dependent counting-variance term $\tfrac{1}{2}\delta_\lambda[u + \beta_\lambda + g(y)]f_{yy}(y, u)$; although stationary flow balance centers the prefactor $u + \beta_\lambda + g(y)$, its product with $f_{yy}(y, u)$ need not have zero mean. Second, $(\lambda^{-\alpha} - c_\lambda^2/2)f_{yy}(y, u)$ records the mismatch in the mean counting-variance rate and vanishes identically when $c_\lambda = c_\lambda^*$. Third, $(\sigma^2/2)\eps_\lambda u f_{uu}(y, u)$ records the discrepancy between the exact CIR volatility $\sigma\sqrt{1 + \eps_\lambda u}$ and the constant OU volatility $\sigma$.

Bounding the expectation of \eqref{eq:master-comparison} under $\pi^X$ requires coordinate-specific estimates on $f_{yy}$, $f_{uu}$, and queue-direction increments of $f_{yy}$ that track their exact dependence on $c_\lambda \in (0, \bar{c}]$. Specifically, with $w = w(y, u)$ and $S = 1 + |x| + |y| + |u|$, Lemmas~\ref{lemma:stein-factor} and~\ref{lem:ou-fyy-lipschitz} establish that $f \in C^2(\R^2)$ satisfies, for all $x, y, u \in \R$,
\begin{equation}
\begin{aligned}
|f_{yy}(y, u)| &\le C H\bigl(1 + |\theta - \mu|c_\lambda^{-1/2}\bigr)w^{m-1}, \\
|f_{yu}(y, u)| + |f_{uu}(y, u)| & \le C H w^{m-1}, \\
|f_{yy}(y, u) - f_{yy}(x, u)| &\le C H |y - x|
\begin{cases}
c_\lambda^{-5/2} S^m, & \mu \ne \theta, \\
(1 + c_\lambda^{-1/2}) S^{m-2}, & \mu = \theta,
\end{cases}
\end{aligned}
\label{eq:ou_stein_bounds}
\end{equation}
where $C < \infty$ is independent of $\lambda$, $c_\lambda$, $h$, and $H$. Separating the $c_\lambda$-dependence across coordinates is essential: if one applied the queue-direction bound $O((c_\lambda^*)^{-1/2}) = O(\lambda^{\alpha/4})$ to $f_{uu}$ in \eqref{eq:master-comparison} when $\mu \neq \theta$, the resulting bound on the intensity discrepancy would be $O(\eps_\lambda (c_\lambda^*)^{-1/2}) = O(\lambda^{-(2-3\alpha)/4})$, which fails to vanish once $\alpha \ge 2/3$. Likewise, standard multidimensional Stein bounds that scale inversely with the minimum volatility across equations \citep{Meckes2009, Gorham2019} diverge as $c_\lambda^* \to 0$, whereas existing bounds for piecewise-linear drifts \citep{Gurvich2014, BravermanDai2017} rely solely on $c_\lambda B_1$ near $y = 0$ and would inflate $|f_{yy}|$ to $O(c_\lambda^{-1})$.

Lemmas~\ref{lemma:stein-factor} and~\ref{lem:ou-fyy-lipschitz} overcome both obstacles through two technical ideas. First, to bound the second derivatives in Lemma~\ref{lemma:stein-factor} when $\mu \neq \theta$, smoothing $g$ by a unit-mass mollifier $\rho_\eta$ gives $g_\eta''(y) = (\theta - \mu)\rho_\eta(y)$, which reduces $|(P_t h)_{ab}|$ to bounding the expected time $Y_\lambda^{\mathrm{OU}}(s)$ spends near $y = 0$. Over short horizons $s \in (0, 1]$, the queue-length Brownian motion $c_\lambda B_1(s)$ contributes fluctuations of standard deviation $c_\lambda\sqrt{s}$, while the intensity Brownian motion $B_2$ propagates through the drift $\int_0^s U(r)\,dr$ on the scale $s^{3/2}$ \citep{DelarueMenozzi2010}. Malliavin integration by parts with a primitive of the occupation test function gives a short-time occupation density bounded by a constant times $k_{c_\lambda}(s)$. The two Brownian perturbations used below leave the terminal arrival-intensity variation zero, so the same estimate holds with polynomial weights in the arrival intensity. The two noise scales balance at $s \asymp c_\lambda$:
\begin{equation*}
k_{c_\lambda}(s) = (c_\lambda^2 s + s^3)^{-1/2},
\qquad
\int_0^1 k_{c_\lambda}(s)\,ds = O(c_\lambda^{-1/2}),
\qquad
\int_0^1 s^j k_{c_\lambda}(s)\,ds = O(1)
\quad (j = 1, 2).
\end{equation*}
Integrating $k_{c_\lambda}(s)$ over $(0, 1]$ improves the naive single-noise bound $\int_0^1 (c_\lambda\sqrt{s})^{-1}\,ds = O(c_\lambda^{-1})$ to $O(c_\lambda^{-1/2})$ for $f_{yy}$, while each intensity derivative supplies an extra factor of $s$ from the drift time integral and keeps $f_{yu}$ and $f_{uu}$ uniform in $c_\lambda$.

Second, the same Malliavin integration-by-parts argument \citep{nualart2006malliavin} bounds queue-length derivatives of the smoothed semigroup and hence the increments in Lemma~\ref{lem:ou-fyy-lipschitz}. We select the smaller bound from two Brownian perturbations: perturbing $B_1$ directly has $L^2$ cost $O((c_\lambda\sqrt{t})^{-1})$ on short horizons $t \lesssim c_\lambda$, whereas perturbing $B_2$ along a cubic profile accumulates through the drift at $L^2$ cost $O(t^{-3/2})$ independently of $c_\lambda$ on horizons $t \gtrsim c_\lambda$. With $\mathcal{L}_\eta$ denoting the smoothed OU generator and $b_\eta(y, u) = -\beta_\lambda + u - g_\eta(y)$, differentiating the smoothed Poisson equation twice in $u$ and once in $y$ gives
\[
\begin{aligned}
(\mathcal{L}_\eta - 2\kappa)(f_\eta)_{uu}
&= (h_\eta)_{uu} - 2(f_\eta)_{yu},\\
\frac{c_\lambda^2}{2}(f_\eta)_{yyy}
&= (h_\eta)_y + g_\eta'(f_\eta)_y - b_\eta(f_\eta)_{yy}
+ \kappa u(f_\eta)_{yu} - \frac{\sigma^2}{2}(f_\eta)_{yuu}.
\end{aligned}
\]
Applying the $O(k_{c_\lambda}(t))$ derivative bound to the discounted semigroup representation of $(f_\eta)_{uu}$ bounds $|(f_\eta)_{yuu}|$ by $O(1 + c_\lambda^{-1/2})$, and dividing the second identity by $c_\lambda^2/2$ yields the weak derivative bound $|f_{yyy}| = O(c_\lambda^{-5/2})$ when $\mu \neq \theta$ (whereas for $\mu = \theta$, the linear drift gives $f_{yy} = -\int_0^\infty e^{-2\mu t} P_t h_{yy}\,dt$ directly). Taking the smaller of the size and increment bounds for $f_{yy}$ in Section~\ref{subsec:ou-main-proof} proves \eqref{eq:ou_bound_second_order_only} and \eqref{eq:ou_bound_variance_matched}. In Section~\ref{sec:cir_results}, we retain the exact CIR intensity dynamics to eliminate the $(\sigma^2/2)\eps_\lambda u f_{uu}$ bottleneck in \eqref{eq:master-comparison}.

\section{Analysis of the CIR-Preserving Diffusion Approximation}
\label{sec:cir_results}

Proposition~\ref{prop:ou-sharpness} shows that replacing the positively skewed stationary distribution of the CIR arrival intensity with a symmetric Gaussian distribution creates a sharp $O(\lambda^{-(1-\alpha)/2})$ error bottleneck. To eliminate this arrival-intensity error while retaining tractability for discrete queue jumps, we now preserve the exact CIR intensity dynamics and approximate only the counting fluctuations from arrivals, service completions, and abandonments.

\subsection{The CIR-Preserving Diffusion and Its Error Bounds}
\label{subsec:cir-construction}

Replacing the scaled counting-fluctuation martingale $\delta_\lambda M_\lambda(t)$ by $c_\lambda B_1(t)$ in \eqref{eq:exact_scaled_system} while retaining the scaled CIR intensity SDE yields the \emph{CIR-preserving two-scale diffusion} $(Y_\lambda^{\mathrm{CIR}}, Z_\lambda)$:
\begin{equation}
\begin{aligned}
dY_\lambda^{\mathrm{CIR}}(t)
&=
\bigl[
-\beta_\lambda + Z_\lambda(t)
- g(Y_\lambda^{\mathrm{CIR}}(t))
\bigr]\,dt
+ c_\lambda\,dB_1(t), \\
dZ_\lambda(t)
&=
-\kappa Z_\lambda(t)\,dt
+
\sigma\sqrt{1 + \eps_\lambda Z_\lambda(t)}\,dB_2(t),
\end{aligned}
\label{eq:cir_diffusion_general_c}
\end{equation}
where $B_1$ and $B_2$ are independent standard Brownian motions. Here $Z_\lambda$ obeys the same SDE as $\wt I_\lambda$, and $(Y_\lambda^{\mathrm{CIR}}, Z_\lambda)$ evolves on the state space $\mathsf{D}_\lambda = \R \times (-\eps_\lambda^{-1}, \infty)$, whose lower boundary $u = -\eps_\lambda^{-1}$ corresponds to zero unscaled intensity ($I_\lambda = 0$) and is inaccessible under Assumption~\ref{assump:model_primitives}. By Lemma~\ref{lem:diffusion_drift}, $(Y_\lambda^{\mathrm{CIR}}, Z_\lambda)$ admits a unique stationary distribution $\pi^{\mathrm{CIR}}$ whose polynomial moments of every fixed order are uniformly bounded for sufficiently large $\lambda$.

Because $Z_\lambda$ and $\wt I_\lambda$ share the exact same generator, the arrival-intensity volatility discrepancy $(\sigma^2/2)\eps_\lambda u f_{uu}$ in \eqref{eq:master-comparison} vanishes identically for every $c_\lambda > 0$. Consequently, the generator difference reduces to
\begin{equation}
(\calG_X - \calG_{\mathrm{CIR}})f(y, u)
= R_{\mathrm{queue}}f(y, u)
+ \left(\lambda^{-\alpha} - \frac{c_\lambda^2}{2}\right) f_{yy}(y, u),
\label{eq:cir-comparison}
\end{equation}
whose full expression is given in \eqref{eq:gen_diff_cir_general_c}. Choosing the variance-matched volatility $c_\lambda = c_\lambda^* = \sqrt{2}\,\lambda^{-\alpha/2}$ cancels the mean counting-variance mismatch $(\lambda^{-\alpha} - c_\lambda^2/2) f_{yy}$, leaving only $R_{\mathrm{queue}}f$: the two Taylor remainders for jumps of size $\pm\delta_\lambda$ and the zero-mean state-dependent counting-variance term.

\begin{theorem}
\label{thm:cir_smooth_bound}
Suppose Assumptions~\ref{assump:model_primitives} and~\ref{assump:staffing_c} hold and $0 < c_\lambda \le \bar c$. Let $h$ satisfy \eqref{eq:test_class} for fixed $m$ and $H$. There is a constant $C < \infty$, depending only on the fixed model parameters, $\bar c$, $m$, and $H$, such that, for all sufficiently large $\lambda$,
\begin{equation}
\left|\pi^X(h) - \pi^{\mathrm{CIR}}(h)\right|
\le C
\begin{cases}
\displaystyle
c_\lambda^{-1/2}
\left[
\lambda^{-(1+\alpha)/2} +
\left|\lambda^{-\alpha} - \frac{c_\lambda^2}{2}\right|
\right],
& \mu \ne \theta, \\[4mm]
\displaystyle
\lambda^{-(1+\alpha)/2} +
\left|\lambda^{-\alpha} - \frac{c_\lambda^2}{2}\right|,
& \mu = \theta.
\end{cases}
\label{eq:cir-general-c-bound}
\end{equation}
Under variance matching, $c_\lambda = c_\lambda^* = \sqrt{2}\,\lambda^{-\alpha/2}$,
\begin{equation}
\left|\pi^X(h) - \pi^{\mathrm{CIR}}(h)\right| \le C
\begin{cases}
\lambda^{-(2+\alpha)/4}, & \mu \ne \theta, \\
\lambda^{-(1+\alpha)/2}, & \mu = \theta.
\end{cases}
\label{eq:cir-smooth-bound}
\end{equation}
\end{theorem}

Comparing \eqref{eq:cir-smooth-bound} with Theorem~\ref{thm:smooth_bound_general_c} highlights two structural phenomena under overdispersion. First, the OU and CIR error bounds respond in opposite directions to $\alpha \in [0, 1)$. Whereas the OU bound $O(\lambda^{-(1-\alpha)/2})$ deteriorates as $\alpha$ increases because the third-moment skewness of $\wt I_\lambda(\infty)$ decays at the slower rate $\eps_\lambda = \lambda^{-(1-\alpha)/2}$, preserving the CIR intensity removes that error and ties accuracy to the scaled jump size $\delta_\lambda = \lambda^{-(1+\alpha)/2}$, which shrinks \emph{faster} relative to the $O(\lambda^{(1+\alpha)/2})$ normalization as $\alpha$ increases. As a result, while both models share the $O(\lambda^{-1/2})$ rate at $\alpha = 0$, the CIR-preserving bound improves upon the OU bound by a factor of $\lambda^{3\alpha/4}$ when $\mu \neq \theta$ and $\lambda^\alpha$ when $\mu = \theta$ for every $\alpha \in (0, 1)$.

Second, overdispersion separates the equal-rate ($\mu = \theta$) and unequal-rate ($\mu \neq \theta$) regimes in a way not seen in single-scale queues. In classical Stein-method analyses of Poisson-driven Erlang-A queues ($\alpha = 0$), the variance-matched volatility $c_\lambda^* = \sqrt{2}$ is of order one, so both $\mu = \theta$ and $\mu \neq \theta$ achieve the same $O(\lambda^{-1/2})$ error rate \citep{Gurvich2014, BravermanDaiFeng2017}. When $\alpha \in (0, 1)$, by contrast, a linear queue-length drift ($\mu = \theta$) attains the full jump-scale rate $O(\delta_\lambda) = O(\lambda^{-(1+\alpha)/2})$, whereas nondifferentiability of $g(y)$ at full occupancy ($\mu \neq \theta$) interacts with the vanishing volatility $c_\lambda^* = \sqrt{2}\,\lambda^{-\alpha/2}$ to introduce the factor $(c_\lambda^*)^{-1/2} \asymp \lambda^{\alpha/4}$, yielding $O(\delta_\lambda (c_\lambda^*)^{-1/2}) = O(\lambda^{-(2+\alpha)/4})$.

\subsection{Queue-Length Moments}
\label{subsec:cir-queue-moments}

Applying Theorem~\ref{thm:cir_smooth_bound} to $h(y, u) = y$ and $h(y, u) = y^2$ yields steady-state error bounds for the mean and variance of the scaled queue length.

\begin{corollary}
\label{cor:cir_moments}
Under the conditions of Theorem~\ref{thm:cir_smooth_bound} and variance matching $c_\lambda = c_\lambda^*$,
\begin{equation}
\left|\E\wt X_\lambda(\infty)
      - \E Y_\lambda^{\mathrm{CIR}}(\infty)\right|
+ \left|\Var\bigl(\wt X_\lambda(\infty)\bigr)
      - \Var\bigl(Y_\lambda^{\mathrm{CIR}}(\infty)\bigr)\right| \le C
\begin{cases}
\lambda^{-(2+\alpha)/4}, & \mu \ne \theta, \\
\lambda^{-(1+\alpha)/2}, & \mu = \theta.
\end{cases}
\label{eq:cir-moment-bound}
\end{equation}
\end{corollary}

In customer-count units, multiplying the mean-error bound in \eqref{eq:cir-moment-bound} by $\delta_\lambda^{-1}=\lambda^{(1+\alpha)/2}$ gives an $O(\lambda^{\alpha/4})$ upper bound when $\mu\ne\theta$, compared with $O(\lambda^\alpha)$ for OU. The CIR bound is $o(\lambda^{1/2})$ for every fixed $\alpha\in[0,1)$.

For $\mu=\theta$, the linear drift $g(y)=\mu y$ closes the stationary generator equations for polynomials of total degree at most two. The resulting identities sharpen the bounds in Corollaries~\ref{cor:mean_variance} and~\ref{cor:cir_moments} to exact agreement at every admissible $\lambda$.

\begin{proposition}
\label{prop:cir-exact-moments}
Suppose Assumptions~\ref{assump:model_primitives} and~\ref{assump:staffing_c} hold, $\mu = \theta = a$, and $c_\lambda = c_\lambda^*$. For every admissible $\lambda$ and every polynomial $p(y, u)$ of total degree at most two, $\pi^X(p) = \pi^{\mathrm{OU}}(p) = \pi^{\mathrm{CIR}}(p)$.
\end{proposition}

To see why variance matching produces exact second-order agreement, let $(Y_\infty,Z_\infty)$ denote the stationary state under any of $\pi^X$, $\pi^{\mathrm{OU}}$, and $\pi^{\mathrm{CIR}}$. The three scaled arrival intensities have mean zero and variance $\sigma^2/(2\kappa)$, although the OU stationary distribution differs from the exact CIR distribution. For $p(y,u)=y^2$, the second-order jump contribution in $\calG_Xy^2$ has stationary mean $2\lambda^{-\alpha}$ by flow balance, whereas either diffusion contributes $c_\lambda^2$. Variance matching equates these rates. The stationary generator identities for $y$, $yu$, and $y^2$ then give the same moments under all three distributions (Section~\ref{subsec:cir-moment-proofs}):
\begin{equation}
\E[Y_\infty] = -\frac{\beta_\lambda}{a},
\qquad
\E[Y_\infty Z_\infty] = \frac{\sigma^2}{2\kappa(a + \kappa)},
\qquad
\Var(Y_\infty) = \frac{\lambda^{-\alpha}}{a} + \frac{\sigma^2}{2\kappa a(a + \kappa)}.
\label{eq:cir_exact_moments_display}
\end{equation}

\subsection{Derivative Bounds near the Lower Intensity Boundary}
\label{subsec:cir-proof-outline}

Proving Theorem~\ref{thm:cir_smooth_bound} requires bounding $f_{yy}$ and its queue-direction increments over steps of size $\pm\delta_\lambda$ for the centered solution $f(y, u) = -\int_0^\infty (P_t^{\mathrm{CIR}} h(y, u) - \pi^{\mathrm{CIR}}(h))\,dt$ to $\calG_{\mathrm{CIR}}f = h - \pi^{\mathrm{CIR}}(h)$, where $(P_t^{\mathrm{CIR}})_{t \ge 0}$ is the transition semigroup of $(Y_\lambda^{\mathrm{CIR}}, Z_\lambda)$. Unlike the OU-based diffusion and degenerate SDEs with uniformly elliptic driving coordinates \citep{DelarueMenozzi2010, Menozzi2011}, the CIR intensity volatility $\sigma\sqrt{v_\lambda(u)}$, where $v_\lambda(u) = 1 + \eps_\lambda u$, approaches zero as $u$ nears the lower boundary $-\eps_\lambda^{-1}$. Differentiating this state-dependent volatility introduces inverse powers of $v_\lambda(u)$, so the Poisson-equation derivative bounds acquire a boundary weight $\omega(u) = 1 + v_\lambda(u)^{-q}$. Specifically, Lemma~\ref{lem:cir_queue_stein} shows that there exist fixed exponents $p \ge m$ and $q > 0$, independent of $\lambda$ and $c_\lambda$, such that for all sufficiently large $\lambda$, $x, y \in \R$, and $u > -\eps_\lambda^{-1}$ (with $w = w(y, u)$ and $S = 1 + |x| + |y| + |u|$),
\begin{equation}
\begin{aligned}
|f_{yy}(y, u)| &\le C H\bigl(1 + |\theta - \mu|c_\lambda^{-1/2}\bigr)w^p\omega(u), \\
|f_{yu}(y, u)| + |f_{uu}(y, u)| &\le C H w^p\omega(u), \\
|f_{yy}(y, u) - f_{yy}(x, u)| &\le C H |y - x| S^p \omega(u)
\begin{cases}
c_\lambda^{-5/2}, & \mu \ne \theta, \\
1 + c_\lambda^{-1/2}, & \mu = \theta.
\end{cases}
\end{aligned}
\label{eq:cir_stein_bounds}
\end{equation}
Notably, \eqref{eq:cir_stein_bounds} recovers the exact same powers of $c_\lambda$ as the constant-volatility OU bounds in \eqref{eq:ou_stein_bounds}, isolating the entire effect of CIR boundary degeneracy inside the multiplicative weight $\omega(u)$.

Establishing Lemma~\ref{lem:cir_queue_stein} and integrating $\omega(u)$ in \eqref{eq:cir-comparison} require three new analytical ingredients beyond the OU proof. First, to control pathwise sensitivities near $v_\lambda(u) = 0$, note that $Z_\lambda(s)$ responds nonlinearly to its initial state $u$, producing nonvanishing first and second flow derivatives $\partial_u Z_\lambda(s)$ and $\partial_u^2 Z_\lambda(s)$. Setting $V_s = v_\lambda(Z_\lambda(s))$ and applying the Lamperti transformation $R_s = 2\sqrt{V_s}/(\sigma\eps_\lambda)$ converts the CIR SDE into unit-volatility form:
\begin{equation*}
dR_s = \left(\frac{k_\lambda - 1/2}{R_s} - \frac{\kappa}{2}R_s\right)ds + dB_2(s),
\qquad
k_\lambda = \frac{2\kappa}{\sigma^2\eps_\lambda^2} = \frac{2\kappa}{\sigma^2}\lambda^{1-\alpha} \to \infty.
\end{equation*}
Differentiating $R_s$ along sample paths introduces inverse powers of $R_s$ and hence requires bounding $\sup_{s \le t} V_s^{-1}$. For sufficiently large $\lambda$, the growing shape parameter $k_\lambda$ gives $V_s^{-Q}$ a strong negative generator drift near zero, which, combined with supermartingale maximal inequalities over unit intervals, bounds pathwise inverse moments of every fixed order in terms of $1 + v_\lambda(u)^{-Q}$. In stationarity, moreover, both $v_\lambda(\wt I_\lambda(\infty))$ and $v_\lambda(Z_\lambda(\infty))$ follow a $\mathrm{Gamma}(k_\lambda, 1/k_\lambda)$ distribution with mean one whose inverse moments satisfy $\E[v_\lambda(\wt I_\lambda(\infty))^{-q}] = k_\lambda^q \Gamma(k_\lambda - q)/\Gamma(k_\lambda) \to 1$ as $\lambda \to \infty$, keeping the weighted polynomial moments $\pi^X(w^\ell \omega)$ uniformly bounded (Lemma~\ref{lem:cir_inverse_moments}).

Second, to obtain the short-time smoothing kernel $(c_\lambda^2 t + t^3)^{-1/2}$ via Malliavin integration by parts, the deterministic cubic perturbation of $B_2$ used in the OU proof breaks down. Under state-dependent volatility $\sigma\sqrt{V_s}$, perturbing $B_2$ by an adapted process $b_2(s)$ introduces the martingale term $\sigma\eps_\lambda(2\sqrt{V_s})^{-1}\mathcal{V}_s^u\,dB_2(s)$ into the intensity-variation equation $d\mathcal{V}_s^u$, preventing any deterministic open-loop profile from hitting the required terminal state variation $(\partial_y Y_t^\eta, 0)$ at time $t$. Instead, we construct a closed-loop adapted feedback perturbation $b_2(s)$ driven by the normalized tracking errors
\begin{equation}
\chi_s = 1 - \frac{\mathcal{V}_s^y}{\partial_y Y_s^\eta},
\qquad
\omega_s = \frac{(t - s)\mathcal{V}_s^u}{\partial_y Y_s^\eta},
\qquad
\E_{y,u}\int_0^t |b_2(s)|^2\,ds \le C t^{-3}\bigl[1 + v_\lambda(u)^{-Q}\bigr],
\label{eq:cir_feedback_errors}
\end{equation}
choosing the feedback gains so the deterministic drift matrix of $(\chi_s, \omega_s)$ has eigenvalues $-8$ and $-9$ on the logarithmic time scale $\log(t/(t-s))$. Weighting a quadratic Lyapunov function of $(\chi_s, \omega_s)$ by $1 + V_s^{-L}$ absorbs the singular $V_s^{-1}$ martingale terms and forces the weighted mean-square tracking error to decay as $((t - s)/t)^8$ as $s \uparrow t$, canceling the $(t - s)^{-4}$ feedback-gain singularity near $s = t$ to yield the $O(t^{-3})$ control-energy bound in \eqref{eq:cir_feedback_errors}. Pairing this feedback perturbation of $B_2$ on horizons $t \gtrsim c_\lambda$ with a direct perturbation of $B_1$ on horizons $t \lesssim c_\lambda$ preserves the short-time bound $(c_\lambda^2 t + t^3)^{-1/2}$ up to the spatial weight $w^p\omega(u)$.

Third, twice differentiating the smoothed CIR Poisson equation in $u$ shifts the drift of the generator by $+\sigma^2\eps_\lambda \partial_u$ (from differentiating $\tfrac{1}{2}\sigma^2(1 + \eps_\lambda u)\partial_{uu}$), and applying the weighted derivative bound to the resulting discounted semigroup representation yields $|f_{yyy}| = O(c_\lambda^{-5/2}w^p\omega(u))$ when $\mu \neq \theta$. Crucially, because customer arrivals and departures alter only the queue coordinate $y$ and leave the scaled intensity $u$ unchanged, all intermediate points $y \pm s\delta_\lambda$ in the jump Taylor remainders share the exact same boundary weight $\omega(u)$. Lemma~\ref{lem:cir_inverse_moments} therefore bounds the stationary expectation of $R_{\mathrm{queue}}f$ in Section~\ref{subsec:cir-proofs}, completing the proof of Theorem~\ref{thm:cir_smooth_bound}.

\section{Accuracy--Tractability Trade-Off and Modeling Implications}
\label{sec:discussion}

Table~\ref{tab:approximation-choice} and Figure~\ref{fig:rate-trifurcation} summarize the accuracy--tractability trade-off between the two variance-matched diffusions across $\alpha \in [0, 1)$. Defined on $\R^2$ with constant volatility matrix $\operatorname{diag}(c_\lambda^*, \sigma)$, the OU-based diffusion yields a centered Gaussian intensity with explicit transition densities while preserving the exact stationary mean, variance, and lag-$t$ autocovariance of $\wt I_\lambda$. It is therefore well suited to calculations governed by first- and second-order intensity structure and to the linearly dispersed regime $\alpha = 0$, where all three smooth-test bounds meet at the classical $O(\lambda^{-1/2})$ Erlang-A rate (Figure~\ref{fig:rate-trifurcation}a). Once $\alpha \in (0, 1)$, however, the decay exponents trifurcate. Replacing the positively skewed Gamma intensity distribution with a symmetric Gaussian distribution introduces a $\Theta(\lambda^{-(1-\alpha)/2})$ error in the third central intensity moment, creating a sharp $O(\lambda^{-(1-\alpha)/2})$ bottleneck across every smooth test class with $m \ge 2$ (Proposition~\ref{prop:ou-sharpness}) whose decay exponent $(1-\alpha)/2$ degrades toward zero as $\alpha \uparrow 1$. Retaining the CIR intensity on $\R \times (-\eps_\lambda^{-1}, \infty)$ with state-dependent volatility $\sigma\sqrt{v_\lambda(u)}$ eliminates this third-moment error and reverses the $\alpha$-dependence, raising the decay exponent to $(2+\alpha)/4$ when $\mu \neq \theta$ and $(1+\alpha)/2$ when $\mu = \theta$---improvements of $\lambda^{3\alpha/4}$ and $\lambda^\alpha$ over the OU bound.

\begin{table}[htbp]
\TABLE{Comparison of the variance-matched diffusions ($c_\lambda = c_\lambda^* = \sqrt{2}\,\lambda^{-\alpha/2}$) for fixed $\alpha \in [0, 1)$.\label{tab:approximation-choice}}
{\begin{tabular*}{\textwidth}{@{\extracolsep{\fill}}%
>{\raggedright\arraybackslash}p{.35\textwidth}%
>{\raggedright\arraybackslash}p{.28\textwidth}%
>{\raggedright\arraybackslash}p{.28\textwidth}@{}}
\toprule
Feature & OU-based diffusion & CIR-preserving diffusion \\
\midrule
State space
  & $\R^2$
  & $\R \times (-\eps_\lambda^{-1}, \infty)$ \\
Volatility matrix
  & $\operatorname{diag}(c_\lambda^*, \sigma)$
  & $\operatorname{diag}(c_\lambda^*, \sigma\sqrt{v_\lambda(u)})$ \\
Stationary distribution of scaled intensity
  & Centered Gaussian
  & Exact scaled Gamma \\
\addlinespace
Error in third central moment
  & $\Theta(\lambda^{-(1-\alpha)/2})$
  & $0$ \\
Smooth-test error bound, $\mu \ne \theta$
  & $O(\lambda^{-(1-\alpha)/2})$
  & $O(\lambda^{-(2+\alpha)/4})$ \\
Smooth-test error bound, $\mu = \theta$
  & $O(\lambda^{-(1-\alpha)/2})$
  & $O(\lambda^{-(1+\alpha)/2})$ ($0$ for degree $\le 2$) \\
\bottomrule
\end{tabular*}}
{\NOTE{Note.} Here $v_\lambda(u) = 1 + \eps_\lambda u$. Both diffusions match the exact stationary mean $0$, variance $\sigma^2/(2\kappa)$, and lag-$t$ autocovariance $\sigma^2 e^{-\kappa t}/(2\kappa)$ ($t \ge 0$) of $\wt I_\lambda$. The $O(\lambda^{-(1-\alpha)/2})$ OU rate is sharp for every $m \ge 2$ (Proposition~\ref{prop:ou-sharpness}), and the $\mu = \theta$ CIR error vanishes identically on all joint polynomials of degree at most two (Proposition~\ref{prop:cir-exact-moments}).}
\end{table}

\begin{figure}[htbp]
\FIGURE
{%
\definecolor{cirEqColor}{RGB}{25,65,135}%
\definecolor{cirNeqColor}{RGB}{0,128,115}%
\definecolor{ouColor}{RGB}{196,48,43}%
\definecolor{refGray}{RGB}{115,115,115}%
\begin{minipage}{\textwidth}
\centering
\begin{tikzpicture}
\begin{axis}[
  hide axis,
  scale only axis,
  xmin=0, xmax=1, ymin=0, ymax=1,
  width=0.95\textwidth, height=0.08\textwidth,
  legend columns=3,
  legend style={
    draw=gray!45,
    fill=white,
    font=\footnotesize,
    column sep=1.2em,
    /tikz/every even column/.append style={column sep=1.2em}
  }
]
\addlegendimage{cirEqColor, line width=1.25pt}
\addlegendentry{CIR-preserving ($\mu = \theta$)}
\addlegendimage{cirNeqColor, dashed, line width=1.25pt}
\addlegendentry{CIR-preserving ($\mu \neq \theta$)}
\addlegendimage{ouColor, dashdotted, line width=1.25pt}
\addlegendentry{OU-based (sharp, $m \ge 2$)}
\end{axis}
\end{tikzpicture}\endgraf
\vspace{2pt}
\begin{minipage}[t]{0.485\textwidth}
\centering
\begin{tikzpicture}
\begin{axis}[
  font=\scriptsize,
  width=0.96\textwidth,
  height=0.72\textwidth,
  xmin=0, xmax=1.03,
  ymin=0, ymax=1.06,
  xlabel={Dispersion parameter $\alpha$},
  ylabel={Scaled error decay exponent $\eta(\alpha)$},
  title={(a) Smooth-test error $O(\lambda^{-\eta(\alpha)})$},
  title style={font=\footnotesize, yshift=-2pt},
  label style={font=\footnotesize},
  tick label style={font=\scriptsize},
  xtick={0, 0.25, 0.5, 0.75, 1},
  ytick={0, 0.25, 0.5, 0.75, 1},
  yticklabels={$0$, $1/4$, $1/2$, $3/4$, $1$},
  grid=major,
  grid style={gray!18, thin},
  axis line style={gray!70},
  clip=false
]
  \addplot[cirEqColor, line width=1.25pt, domain=0:0.985, samples=50] {(1+x)/2};
  \addplot[cirNeqColor, dashed, line width=1.25pt, domain=0:0.985, samples=50] {(2+x)/4};
  \addplot[ouColor, dashdotted, line width=1.25pt, domain=0:0.985, samples=50] {(1-x)/2};
  \addplot[only marks, mark=*, mark size=1.8pt, draw=cirEqColor, fill=white, line width=0.9pt] coordinates {(1, 1)};
  \addplot[only marks, mark=*, mark size=1.8pt, draw=cirNeqColor, fill=white, line width=0.9pt] coordinates {(1, 0.75)};
  \addplot[only marks, mark=*, mark size=1.8pt, draw=ouColor, fill=white, line width=0.9pt] coordinates {(1, 0)};
  \addplot[only marks, mark=*, mark size=2.1pt, draw=black, fill=black] coordinates {(0, 0.5)};
  \node[font=\scriptsize, anchor=west, fill=white, inner sep=1.2pt] at (axis cs:0.03, 0.43) {Erlang-A};
  \node[font=\scriptsize, text=cirEqColor, anchor=south east] at (axis cs:0.56, 0.79) {$(1+\alpha)/2$};
  \node[font=\scriptsize, text=cirNeqColor, anchor=north west] at (axis cs:0.38, 0.58) {$(2+\alpha)/4$};
  \node[font=\scriptsize, text=ouColor, anchor=north east] at (axis cs:0.52, 0.23) {$(1-\alpha)/2$};
  \draw[<->, >=stealth, gray!80!black, thin] (axis cs:0.82, 0.105) -- (axis cs:0.82, 0.69)
    node[midway, right=1.5pt, font=\scriptsize, text=gray!85!black, fill=white, inner sep=1pt] {$3\alpha/4$};
\end{axis}
\end{tikzpicture}
\end{minipage}%
\hfill%
\begin{minipage}[t]{0.485\textwidth}
\centering
\begin{tikzpicture}
\begin{axis}[
  font=\scriptsize,
  width=0.96\textwidth,
  height=0.72\textwidth,
  xmin=0, xmax=1.03,
  ymin=-0.14, ymax=1.06,
  xlabel={Dispersion parameter $\alpha$},
  ylabel={Unscaled error exponent $\zeta(\alpha)$},
  title={(b) Unscaled queue error $O(\lambda^{\zeta(\alpha)})$},
  title style={font=\footnotesize, yshift=-2pt},
  label style={font=\footnotesize},
  tick label style={font=\scriptsize},
  xtick={0, 0.25, 0.5, 0.75, 1},
  ytick={0, 0.25, 0.5, 0.75, 1},
  yticklabels={$0$, $1/4$, $1/2$, $3/4$, $1$},
  grid=major,
  grid style={gray!18, thin},
  axis line style={gray!70},
  clip=false
]
  \addplot[refGray, densely dotted, line width=0.95pt, domain=0:0.985, samples=50] {(1+x)/2};
  \addplot[refGray!70, thin, densely dashed, domain=0:1] {0.5};
  \addplot[ouColor, dashdotted, line width=1.25pt, domain=0:0.985, samples=50] {x};
  \addplot[cirNeqColor, dashed, line width=1.25pt, domain=0:0.985, samples=50] {x/4};
  \addplot[cirEqColor, line width=1.25pt, domain=0:0.985, samples=50] {0};
  \addplot[only marks, mark=*, mark size=1.8pt, draw=refGray, fill=white, line width=0.8pt] coordinates {(1, 1)};
  \addplot[only marks, mark=*, mark size=1.8pt, draw=ouColor, fill=white, line width=0.9pt] coordinates {(1, 1)};
  \addplot[only marks, mark=*, mark size=1.8pt, draw=cirNeqColor, fill=white, line width=0.9pt] coordinates {(1, 0.25)};
  \addplot[only marks, mark=*, mark size=1.8pt, draw=cirEqColor, fill=white, line width=0.9pt] coordinates {(1, 0)};
  \addplot[only marks, mark=*, mark size=2.1pt, draw=black, fill=black] coordinates {(0, 0)};
  \node[font=\scriptsize, text=refGray, anchor=south west] at (axis cs:0.04, 0.78) {Safety capacity $(1+\alpha)/2$};
  \node[font=\scriptsize, text=refGray!85!black, anchor=south west, fill=white, inner sep=1pt] at (axis cs:0.76, 0.505) {$\lambda^{1/2}$ scale};
  \node[font=\scriptsize, text=ouColor, anchor=north west, xshift=1pt, yshift=3pt] at (axis cs:0.66, 0.66) {$\alpha$};
  \node[font=\scriptsize, text=cirNeqColor, anchor=south east] at (axis cs:0.78, 0.20) {$\alpha/4$};
  \node[font=\scriptsize, text=cirEqColor, anchor=south] at (axis cs:0.48, -0.12) {$0$ (exact $0$ for degree $\le 2$)};
\end{axis}
\end{tikzpicture}
\end{minipage}
\end{minipage}%
}
{Effect of the dispersion parameter $\alpha \in [0, 1)$ on steady-state approximation rates under variance matching ($c_\lambda = c_\lambda^* = \sqrt{2}\,\lambda^{-\alpha/2}$).\label{fig:rate-trifurcation}}
{Panel~(a) plots the decay exponent $\eta(\alpha)$ in $|\pi^X(h) - \pi^D(h)| = O(\lambda^{-\eta(\alpha)})$ from Theorems~\ref{thm:smooth_bound_general_c} and~\ref{thm:cir_smooth_bound}. Panel~(b) plots the corresponding unscaled customer-count exponent $\zeta(\alpha) = (1+\alpha)/2 - \eta(\alpha)$ obtained by multiplying scaled queue-length bounds by $\delta_\lambda^{-1} = \lambda^{(1+\alpha)/2}$, alongside the $O(\lambda^{(1+\alpha)/2})$ safety-capacity scale.}
\end{figure}

At finite system scale $\lambda$, accuracy also depends on the target observable and on whether congestion is reported in scaled or physical units. Although our smooth-test bounds hold uniformly over each class $(m, H)$, low-order moments can enjoy additional cancellation---most notably when $\mu = \theta$, where the variance-matched CIR diffusion reproduces the exact stationary mean, variance, and queue--intensity covariance at every admissible $\lambda$ (Proposition~\ref{prop:cir-exact-moments}). When translating scaled bounds into customer counts via $X_\lambda(\infty) = n_\lambda + \delta_\lambda^{-1}\wt X_\lambda(\infty)$, multiplying mean errors by $\delta_\lambda^{-1} = \lambda^{(1+\alpha)/2}$ shifts every exponent as shown in Figure~\ref{fig:rate-trifurcation}b (and scales variance errors by $\delta_\lambda^{-2} = \lambda^{1+\alpha}$). Under the OU approximation, the resulting $O(\lambda^\alpha)$ customer-count bound exceeds the classical $\lambda^{1/2}$ diffusive scale whenever $\alpha > 1/2$ and approaches the $O(\lambda^{(1+\alpha)/2})$ safety-capacity scale as $\alpha \uparrow 1$. Under the CIR-preserving approximation, by contrast, the unscaled mean error remains $O(\lambda^{\alpha/4}) = o(\lambda^{1/2})$ across all $\alpha \in [0, 1)$ when $\mu \neq \theta$, and vanishes identically when $\mu = \theta$.

Finally, our results point to a general principle for designing pre-limit diffusion approximations in multi-scale queueing systems. Retaining the vanishing volatility $c_\lambda^* = \sqrt{2}\,\lambda^{-\alpha/2}$ in the queue-length SDE permits direct steady-state comparison through the Poisson equation without interchanging large-system and steady-state limits \citep{Gurvich2014, BravermanDai2017}. At the same time, the generator comparison \eqref{eq:master-comparison} and our directional Stein bounds reveal that different coefficient approximations interact asymmetrically with vanishing noise: counting-variance mismatches multiply $f_{yy}$, which blows up as $(c_\lambda^*)^{-1/2}$ near full occupancy when $\mu \neq \theta$, whereas arrival-intensity volatility mismatches multiply $f_{uu}$, which stays $O(1)$ in $c_\lambda$ because intensity fluctuations act on the queue through its drift. Tracking which generator discrepancies hit volatility-sensitive derivatives thus provides a systematic criterion for deciding which finite-system terms a diffusion approximation must retain.

\section{Numerical Experiments}
\label{sec:numerical-experiments}

We complement our theoretical bounds with deterministic steady-state experiments examining how arrival-intensity skewness and vanishing counting-fluctuation volatility affect pure queue-length observables. To eliminate Monte Carlo simulation error, we compute the exact stationary distribution $\pi^X$ via a quasi-birth-and-death linear system that expands the arrival intensity in the generalized Laguerre polynomial eigenbasis of the CIR generator (Section~\ref{sec:ec-laguerre-qbd}); even under this spectral reduction, however, the discrete queue lattice spans $O(\lambda^{(1+\alpha)/2})$ active levels around $n_\lambda$. By contrast, we compute $\pi^{\mathrm{OU}}$ and $\pi^{\mathrm{CIR}}$ by solving their two-dimensional stationary Fokker--Planck equations on a fixed, $\lambda$-independent grid that pairs Gauss--Hermite or Gauss--Laguerre spectral collocation in the intensity coordinate with an exponentially fitted Chang--Cooper finite-volume scheme in the queue coordinate \citep{ChangCooper1970}. Throughout, we fix $\mu = 1$, $\kappa = 1$, and $\sigma = 1$, and specify $(\lambda, \alpha, \theta, \beta_\lambda)$ in each experiment.

\begin{figure}[t]
\FIGURE
{\includegraphics[width=\textwidth]{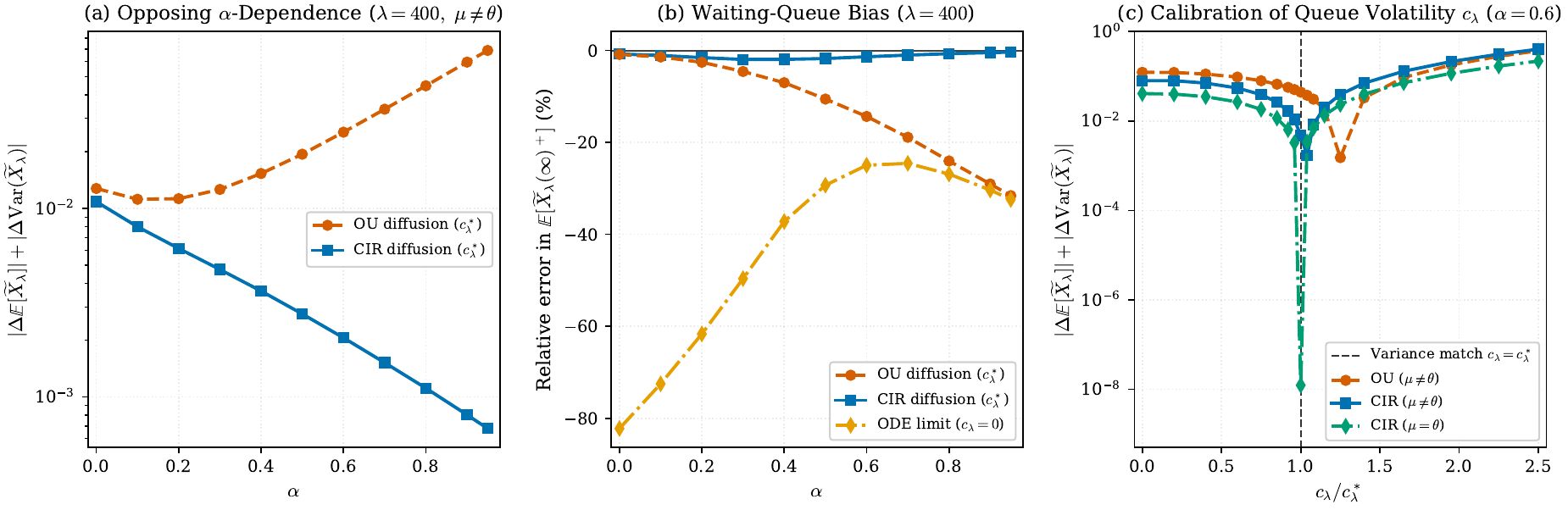}}
{Impact of overdispersion and queue-volatility calibration on steady-state queue-length accuracy.\label{fig:alpha-volatility}}
{Panels~(a) and~(b) fix $\lambda = 400$, $\theta = 0.5$, and $\beta_\lambda = 0.5$ across $\alpha \in [0, 0.95]$; Panel~(b) plots $100 \times \{\E[(Y_\lambda^D(\infty))^+] - \E[\wt X_\lambda(\infty)^+]\} / \E[\wt X_\lambda(\infty)^+]$, the signed relative error (in percent) in the expected scaled queue length. Panel~(c) fixes $\alpha = 0.6$, $\lambda = 200$, and $\beta_\lambda = 0.25$, comparing $\theta = 0.4$ ($\mu \neq \theta$) and $\theta = 1$ ($\mu = \theta$) across $c_\lambda / c_\lambda^* \in [0, 2.5]$.}
\end{figure}

Figure~\ref{fig:alpha-volatility} illustrates how overdispersion and volatility calibration govern accuracy at finite scale. In Panel~(a), the variance-matched OU and CIR approximations exhibit nearly identical queue-moment errors at $\alpha = 0$, matching their common $O(\lambda^{-1/2})$ benchmark in Theorems~\ref{thm:smooth_bound_general_c} and~\ref{thm:cir_smooth_bound}. Once $\alpha > 0$, however, the CIR error decreases by more than an order of magnitude whereas the OU error grows steadily, opening a gap of over two orders of magnitude by $\alpha = 0.95$. Panel~(b) explains this divergence via the signed relative error in the expected scaled queue length $\E[\wt X_\lambda(\infty)^+]$. Dropping counting fluctuations in the heavy-traffic ordinary differential equation (ODE) limit ($c_\lambda = 0$) misses random excursions across full occupancy and underestimates the expected queue length by over $80\%$ at small $\alpha$. Restoring $c_\lambda = c_\lambda^*$ under the OU intensity corrects this small-$\alpha$ bias, yet still underpredicts the expected queue length by more than $30\%$ at $\alpha = 0.95$: because $\theta < \mu$, upward arrival-intensity surges clear more slowly above full occupancy than downward troughs recover below it, so replacing the positively skewed Gamma intensity distribution with a symmetric Gaussian distribution systematically understates congestion. Retaining both $c_\lambda^*$ and the CIR intensity keeps the relative error within $2\%$ across the entire dispersion range.

Panel~(c) of Figure~\ref{fig:alpha-volatility} isolates the calibration role of $c_\lambda^* = \sqrt{2}\,\lambda^{-\alpha/2}$ by varying $c_\lambda / c_\lambda^*$ at $\alpha = 0.6$ and $\lambda = 200$. Under the CIR-preserving diffusion, the combined mean-and-variance error falls to machine zero at $c_\lambda / c_\lambda^* = 1$ when $\mu = \theta$ (verifying Proposition~\ref{prop:cir-exact-moments}) and attains a sharp minimum near $c_\lambda = c_\lambda^*$ when $\mu \neq \theta$, consistent with the variance-mismatch penalty $\Psi_\lambda(c_\lambda)$ in Theorem~\ref{thm:cir_smooth_bound}. Under the OU-based diffusion with $\mu \neq \theta$, by contrast, setting $c_\lambda = c_\lambda^*$ leaves an irreducible error floor driven by arrival-intensity skewness.\footnote{The shallow dip of the OU curve near $c_\lambda / c_\lambda^* \approx 1.25$ reflects a partial sign cancellation when $\theta < \mu$: inflating $c_\lambda$ above $c_\lambda^*$ artificially raises queue-length fluctuations, partially offsetting the downward congestion bias from missing positive Gamma skewness.}

\begin{figure}[t]
\FIGURE
{\includegraphics[width=0.75\textwidth]{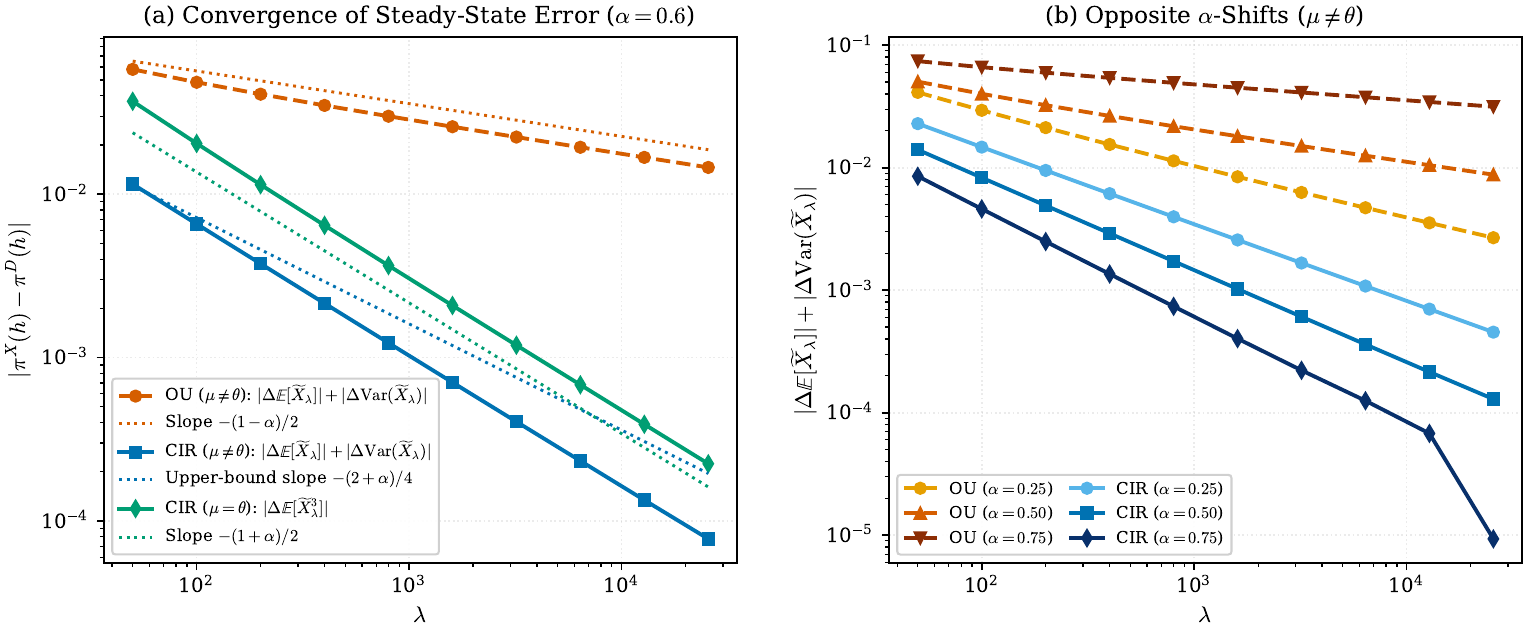}}
{Asymptotic decay of steady-state queue-length errors over $\lambda \in [50, 25600]$ ($\beta_\lambda = 0$).\label{fig:lambda-convergence}}
{Panel~(a) fixes $\alpha = 0.6$ and plots $|\E[\wt X_\lambda(\infty)] - \E[Y_\lambda^D(\infty)]| + |\Var(\wt X_\lambda(\infty)) - \Var(Y_\lambda^D(\infty))|$ at $\theta = 0.4$ ($\mu \neq \theta$) alongside the third-moment error $|\E[\wt X_\lambda(\infty)^3] - \E[Y_\lambda^{\mathrm{CIR}}(\infty)^3]|$ at $\theta = 1$ ($\mu = \theta$, where degree-two errors vanish by Proposition~\ref{prop:cir-exact-moments}), with theoretical reference slopes. Panel~(b) plots the same combined mean-and-variance error at $\theta = 0.4$ across $\alpha \in \{0.25, 0.50, 0.75\}$.}
\end{figure}

Figure~\ref{fig:lambda-convergence} tracks asymptotic convergence over $\lambda \in [50, 25600]$ ($\beta_\lambda = 0$) and sheds light on queue-only rate sharpness. Although Proposition~\ref{prop:ou-sharpness} proves sharpness of the $O(\lambda^{-(1-\alpha)/2})$ OU rate using arrival-intensity test functions, Panel~(a) shows that at $\alpha = 0.6$ and $\mu \neq \theta$, the pure queue-moment error under the OU diffusion decays at the exact slope $-(1-\alpha)/2 = -0.20$: nondifferentiability of the departure drift at full occupancy ($y = 0$) transmits third-moment intensity skewness directly into first- and second-order queue moments. For the CIR-preserving diffusion, the third-moment error at $\mu = \theta$ follows the predicted slope $-(1+\alpha)/2 = -0.80$ from Theorem~\ref{thm:cir_smooth_bound}, whereas at $\mu \neq \theta$ the queue-moment error lies strictly below the $O(\lambda^{-(2+\alpha)/4})$ ($-0.65$) upper bound of Corollary~\ref{cor:cir_moments} and decays near slope $-0.80$, suggesting that smooth queue-only observables may benefit from additional cancellation across full occupancy (Section~\ref{sec:conclusion}). Panel~(b) confirms the opposing $\alpha$-shifts across $\alpha \in \{0.25, 0.50, 0.75\}$: larger $\alpha$ flattens the OU decay trajectories toward slopes $-0.375$, $-0.25$, and $-0.125$ while steepening the CIR trajectories, yielding nearly four orders of magnitude higher accuracy at large $\lambda$ and strong overdispersion.

\section{Conclusion}
\label{sec:conclusion}

For many-server queues driven by overdispersed, temporally dependent DSPP arrivals, accurate pre-limit diffusion approximations require coupling two distinct fluctuation scales across $\alpha \in [0, 1)$. While calibrating an $O(\lambda^{-\alpha/2})$ queue-length volatility eliminates the leading counting-variance error, replacing the CIR arrival intensity with its heavy-traffic OU limit loses positive skewness and degrades steady-state accuracy as overdispersion increases, a penalty that retaining the exact CIR intensity reverses. Methodologically, our directional and boundary-weighted Poisson-equation bounds demonstrate how drift-transmitted fluctuations and path-adapted Brownian perturbations extend Stein's method to multi-dimensional diffusions with vanishing volatilities and nondifferentiable drifts.

Two open problems merit further investigation. First, extending our smooth-test guarantees to steady-state delay probabilities requires handling discontinuous congestion indicators: specifically, the time-congestion probability $\pr(X_\lambda(\infty) \ge n_\lambda)$ and, because Poisson arrivals see time averages (PASTA) fails under a stochastic arrival intensity, the customer-perceived delay probability
$
\lambda^{-1} \E\bigl[I_\lambda(t)\,\ind\{X_\lambda(t-) \ge n_\lambda\}\bigr]
$
in stationarity, which combines the discontinuity at $y = 0$ with the multiplicative weight $I_\lambda(t)/\lambda$. Second, whereas Proposition~\ref{prop:ou-sharpness} establishes sharpness of the OU rate across $m \ge 2$ using arrival-intensity observables, it remains open whether a smooth observable $h(y)$ of the queue length alone attains the $O(\lambda^{-(2+\alpha)/4})$ variance-matched CIR bound in Theorem~\ref{thm:cir_smooth_bound} when $\mu \neq \theta$, or whether the queue-side rate can be sharpened further.

\bibliographystyle{informs2014}
\bibliography{references}

\ECSwitch

\ECHead{Supplementary Material}
\numberwithin{equation}{section}

\section{Generator Comparison and Stationary Moment Bounds}
\label{sec:intermediate}

\subsection{Exact and Diffusion Generators}
\label{subsec:generator_comparison}

For any scaled state $(y,u)\in\mathsf{D}_\lambda$ of the exact process, write the corresponding unscaled queue length and arrival intensity as
\begin{equation}
x=n_\lambda+\delta_\lambda^{-1}y,
\qquad
\iota=\lambda+\delta_\lambda^{-1}u,
\label{eq:state_unscaled}
\end{equation}
with $x$ restricted to the nonnegative integer lattice. For a twice continuously differentiable function $f$ on $\mathsf{D}_\lambda$, the generator of the exact scaled process $(\wt X_\lambda,\wt I_\lambda)$ is
\begin{align*}
\calG_X f(y,u)
={}&\iota\bigl[f(y+\delta_\lambda,u)-f(y,u)\bigr]
+d_\lambda(x)\bigl[f(y-\delta_\lambda,u)-f(y,u)\bigr]\\
&-\kappa u f_u(y,u)
+\frac{\sigma^2}{2}\frac{\iota}{\lambda}f_{uu}(y,u),
\end{align*}
where $d_\lambda(x)$ is the total departure rate in \eqref{eq:death_rate}. By comparison, the generators of the OU-based and CIR-preserving diffusions in \eqref{eq:ou_diffusion_general_c} and \eqref{eq:cir_diffusion_general_c} are
\begin{align*}
\calG_{\mathrm{OU}}f(y,u)
={}&[-\beta_\lambda+u-g(y)]f_y(y,u)-\kappa u f_u(y,u)
+\frac{c_\lambda^2}{2}f_{yy}(y,u)
+\frac{\sigma^2}{2}f_{uu}(y,u),
\\
\calG_{\mathrm{CIR}}f(y,u)
={}&[-\beta_\lambda+u-g(y)]f_y(y,u)-\kappa u f_u(y,u)
+\frac{c_\lambda^2}{2}f_{yy}(y,u)
+\frac{\sigma^2}{2}v_\lambda(u)f_{uu}(y,u).
\end{align*}

To expand the unit-customer jumps of $\calG_X$ around $(y,u)$, define the second-order Taylor averages
\begin{equation}
\mathcal{T}_{\lambda,\pm}f_{yy}(y,u)
:=2\int_0^1(1-s)f_{yy}(y\pm s\delta_\lambda,u)\,\dd s,
\label{eq:integral-taylor-average}
\end{equation}
so that Taylor's formula with integral remainder gives
\[
f(y\pm\delta_\lambda,u)-f(y,u)
=\pm\delta_\lambda f_y(y,u)
+\frac{\delta_\lambda^2}{2}\mathcal{T}_{\lambda,\pm}f_{yy}(y,u).
\]
Because $\delta_\lambda\{\iota-d_\lambda(x)\}=-\beta_\lambda+u-g(y)$ by \eqref{eq:death_rate} and \eqref{eq:staffing}, the first-order jump terms in $\calG_X f$ match the queue-length drift of both diffusions. Adding and subtracting $f_{yy}(y,u)$ inside each Taylor average and noting that $\iota/\lambda = 1 + \eps_\lambda u = v_\lambda(u)$ yields the explicit form of \eqref{eq:master-comparison}:
\begin{align*}
(\calG_X-\calG_{\mathrm{OU}})f(y,u)
={}&\frac{\delta_\lambda^2\iota}{2}
\bigl\{\mathcal{T}_{\lambda,+}f_{yy}(y,u)-f_{yy}(y,u)\bigr\}\\
&+\frac{\delta_\lambda^2 d_\lambda(x)}{2}
\bigl\{\mathcal{T}_{\lambda,-}f_{yy}(y,u)-f_{yy}(y,u)\bigr\}\\
&+\frac{\delta_\lambda}{2}[u+\beta_\lambda+g(y)]f_{yy}(y,u)\\
&+\left(\lambda^{-\alpha}-\frac{c_\lambda^2}{2}\right)f_{yy}(y,u)
+\frac{\sigma^2}{2}\eps_\lambda u f_{uu}(y,u).
\end{align*}
For the CIR-preserving diffusion, the arrival-intensity generators coincide, so the final term vanishes and we obtain
\begin{align}
(\calG_X-\calG_{\mathrm{CIR}})f(y,u)
={}&\frac{\delta_\lambda^2\iota}{2}
\bigl\{\mathcal{T}_{\lambda,+}f_{yy}(y,u)-f_{yy}(y,u)\bigr\}\nonumber\\
&+\frac{\delta_\lambda^2 d_\lambda(x)}{2}
\bigl\{\mathcal{T}_{\lambda,-}f_{yy}(y,u)-f_{yy}(y,u)\bigr\}\nonumber\\
&+\frac{\delta_\lambda}{2}[u+\beta_\lambda+g(y)]f_{yy}(y,u)\nonumber\\
&+\left(\lambda^{-\alpha}-\frac{c_\lambda^2}{2}\right)f_{yy}(y,u).
\label{eq:gen_diff_cir_general_c}
\end{align}

Because the queue-jump and counting-variance terms take the same form in both generator differences, we bound them simultaneously. Abbreviate $\delta=\delta_\lambda$, $r=\lambda^{-\alpha}$, and $c=c_\lambda$, and define the volatility-scaling factors
\[
s(c)=1+|\theta-\mu|c^{-1/2},\qquad
M(c)=\begin{cases}
c^{-5/2},&\mu\ne\theta,\\
1+c^{-1/2},&\mu=\theta.
\end{cases}
\]
At any exact scaled state $(y,u)$, substituting $d_\lambda(x)=\lambda+\delta^{-1}\{\beta_\lambda+g(y)\}$ and $\lambda\delta^2=r$ gives
\[
\delta^2\iota \le r+\delta|u|,
\qquad
\delta^2 d_\lambda(x)\le r+\delta(B_\beta+a_1|y|).
\]
Since $\delta/r=\varepsilon_\lambda\le1$ for $\lambda\ge1$, summing these two bounds yields
\begin{equation}
\delta^2\{\iota+d_\lambda(x)\}\le C r w(y,u).
\label{eq:cir-proof-rate-bound}
\end{equation}

Now suppose that $f$ satisfies the second-derivative and local increment bounds
\[
\begin{aligned}
|f_{yy}(y,u)|&\le CH s(c)w(y,u)^p\omega(u),\\
|f_{yy}(z,u)-f_{yy}(y,u)|&\le CH M(c)|z-y|w(y,u)^p\omega(u)
\quad (|z-y|\le\delta),
\end{aligned}
\]
for fixed exponents $p,q\ge0$ (enlarging $p$ if necessary to cover both estimates) and a constant $C$ independent of $\lambda$ and $c$, where $\omega(u)=1$ for the OU Poisson solution and $\omega(u)=1+v_\lambda(u)^{-q}$ for the CIR Poisson solution. Assume further that $\pi^X(w^{p+1}\omega)$ is uniformly bounded over $0<c\le\bar c$ and all sufficiently large $\lambda$. Because every intermediate point $z=y\pm s\delta$ ($s\in[0,1]$) in \eqref{eq:integral-taylor-average} shares the same intensity coordinate $u$ and satisfies $w(z,u)\le 2w(y,u)$ when $\delta\le1$, the integrand $|\mathcal{T}_{\lambda,\pm}f_{yy}(y,u)-f_{yy}(y,u)|$ is bounded both by $CH\delta M(c)w(y,u)^p\omega(u)$ via the increment estimate and by $CHs(c)w(y,u)^p\omega(u)$ via the triangle inequality. Letting
\[
T_+=\tfrac12\delta^2\pi^X\!\left(\iota\,|\mathcal{T}_{\lambda,+}f_{yy}-f_{yy}|\right),
\qquad
T_-=\tfrac12\delta^2\pi^X\!\left(d_\lambda(x)\,|\mathcal{T}_{\lambda,-}f_{yy}-f_{yy}|\right)
\]
denote the stationary expectations of the upward and downward Taylor remainders, combining the smaller of the two bounds with \eqref{eq:cir-proof-rate-bound} and the uniform bound on $\pi^X(w^{p+1}\omega)$ yields
\begin{equation}
T_++T_-\le CHr\min\{s(c),\delta M(c)\}.
\label{eq:taylor-minimum-corrected}
\end{equation}
Similarly, since $|u+\beta_\lambda+g(y)|\le (1+B_\beta+a_1)w(y,u)$, applying the size bound on $|f_{yy}|$ and the same moment estimate controls the two counting-variance terms:
\begin{align}
\tfrac12\delta\,\pi^X\!\left(|u+\beta_\lambda+g(y)|\,|f_{yy}|\right)
&\le CH\delta s(c), \label{eq:local-second-order-corrected}\\
|r-c^2/2|\,\pi^X(|f_{yy}|)
&\le CH|r-c^2/2|s(c). \label{eq:variance-mismatch-corrected}
\end{align}

\subsection{Stationary Distributions and Polynomial Moment Bounds}
\label{subsec:stationarity}

We next establish uniqueness of the stationary distributions and uniform polynomial moment bounds for the two-scale diffusions (Lemma~\ref{lem:diffusion_drift}) and the exact scaled process (Lemma~\ref{lem:exact_scaled_moments}) using the Foster--Lyapunov criteria of \citet{meyn1993stability}.

\begin{lemma}
\label{lem:diffusion_drift}
Suppose Assumptions~\ref{assump:model_primitives} and~\ref{assump:staffing_c} hold. There exist constants $a_V,b_V,m_0>0$, independent of $\lambda$ and $c_\lambda$, such that $V(y,v)=y^2+m_0v^2$ satisfies
\begin{equation}
\calG_D V(y,v)\le -a_V V(y,v)+b_V
\label{eq:uniform_lyapunov_diffusion}
\end{equation}
for both diffusions $D\in\{\mathrm{OU},\mathrm{CIR}\}$. Consequently, both diffusions admit unique stationary distributions whose polynomial moments of every fixed order are uniformly bounded over sufficiently large $\lambda$ and $0<c_\lambda\le\bar c$. Moreover, for every fixed $r\ge1$ and initial state $(y,v)$, the expected value of $w^r$ at time $t$ is at most $C_r w(y,v)^r$, uniformly over $t\ge0$, admissible $\lambda$, and $0<c_\lambda\le\bar c$.
\end{lemma}

\proof{Proof of Lemma \ref{lem:diffusion_drift}}
Throughout this proof, $C$ depends only on the fixed model parameters, $B_\beta$, and $\bar c$, and may change from line to line; a subscript $p$ allows additional dependence on $p$. In both diffusions, the queue-length drift is $-\beta_\lambda+v-g(y)$, where $yg(y)\ge a_0y^2$. Applying $\calG_{\mathrm{OU}}$ to $V(y,v)=y^2+m_0v^2$ gives
\[
\calG_{\mathrm{OU}} V
=2y[-\beta_\lambda+v-g(y)]-2m_0\kappa v^2
+c_\lambda^2+m_0\sigma^2.
\]
By Young's inequality,
\[
2|yv|\le \frac{a_0}{2}y^2+\frac{2}{a_0}v^2,
\qquad
2B_\beta|y|\le \frac{a_0}{2}y^2+\frac{2B_\beta^2}{a_0},
\]
and hence
\[
\calG_{\mathrm{OU}} V
\le -a_0y^2-\left(2m_0\kappa-\frac{2}{a_0}\right)v^2+C.
\]
Choosing $m_0>1/(\kappa a_0)$ establishes \eqref{eq:uniform_lyapunov_diffusion} for the OU diffusion. For the CIR-preserving diffusion, replacing $\sigma^2$ with $\sigma^2(1+\eps_\lambda v)$ changes the intensity diffusion term from $m_0\sigma^2$ to $m_0\sigma^2(1+\eps_\lambda v)$. Since $\eps_\lambda\le1$, Young's inequality bounds the additional linear term $m_0\sigma^2\eps_\lambda v$ by a constant plus an arbitrarily small multiple of $v^2$; absorbing that multiple into the negative quadratic coefficient yields \eqref{eq:uniform_lyapunov_diffusion} for the CIR diffusion as well.

To deduce positive Harris recurrence at a fixed admissible $\lambda$ and $c_\lambda>0$, observe first that both diffusions are nonexplosive and Feller because the autonomous arrival-intensity coordinate has these properties and the queue-length drift is globally Lipschitz with linear growth. On the open state space, positive continuous transition densities imply irreducibility and render interior compact sets petite. For the CIR-preserving diffusion, we extend the state space to $\overline{\mathsf{D}}_\lambda=\R\times[-\eps_\lambda^{-1},\infty)$ by including the lower boundary $u=-\eps_\lambda^{-1}$, from which the Feller extension enters the interior $\mathsf{D}_\lambda$ at every positive time under Assumption~\ref{assump:model_primitives}. At any fixed positive time, the Feller property and a finite open cover of any compact subset of $\overline{\mathsf{D}}_\lambda$ provide a uniform positive probability of reaching a bounded interior open set with compact closure; combining this step with the interior transition-density bound shows that compact sets meeting the boundary are also petite. Since the sublevel sets of $V$ are compact and petite on $\overline{\mathsf{D}}_\lambda$ and $V$ is bounded on each sublevel set, \eqref{eq:uniform_lyapunov_diffusion} and \citet[Theorem~4.2]{meyn1993stability} yield positive Harris recurrence and a unique stationary distribution $\pi^{\mathrm{CIR}}$, which assigns zero probability to $\{u=-\eps_\lambda^{-1}\}$ because its intensity marginal is Gamma-distributed on $(-\eps_\lambda^{-1},\infty)$. The OU diffusion requires no boundary extension on $\R^2$, and neither case requires a lower bound on $c_\lambda$ that is uniform in $\lambda$.

For higher polynomial moments, fix $p\ge1$ and set $W=1+V$. Applying the diffusion chain rule together with $0<\eps_\lambda\le1$ gives, for either $D\in\{\mathrm{OU},\mathrm{CIR}\}$,
\begin{align*}
\calG_D W^p
&\le pW^{p-1}\calG_D V
 +C_p W^{p-2}(y^2+v^2+|v|^3)\\
&\le -a_V pW^p+p(a_V+b_V)W^{p-1}
 +C_p W^{p-1/2}\\
&\le -\frac{a_V p}{2}W^p+b_p,
\end{align*}
for a constant $b_p<\infty$ independent of $\lambda$ and $0<c_\lambda\le\bar c$. Stopped Dynkin's formula and Gronwall's inequality bound the expected value of $W^p$ at time $t$ by $e^{-a_V pt/2}W(y,v)^p + 2b_p/(a_V p)$, while \citet[Theorem~4.3(ii)]{meyn1993stability} yields the stationary bound $\pi^D(W^p)\le 2b_p/(a_V p)$. Because $W^{r/2}$ is comparable to $w^r$, these estimates establish the stated transient and stationary bounds for $r\ge2$, and orders $1\le r<2$ follow by Jensen's inequality.
\Halmos
\endproof

\begin{lemma}
\label{lem:exact_scaled_moments}
Under Assumptions~\ref{assump:model_primitives} and~\ref{assump:staffing_c}, the exact scaled process $(\wt X_\lambda,\wt I_\lambda)$ has a unique stationary distribution $\pi^X$ at every admissible $\lambda$. Moreover, for each $m\ge1$, there exists $\lambda_0<\infty$ such that
\begin{equation}
\sup_{\lambda\ge\lambda_0}\E\left[|\wt X_\lambda(\infty)|^m+|\wt I_\lambda(\infty)|^m\right]<\infty.
\label{eq:exact_moment_bound}
\end{equation}
\end{lemma}

\proof{Proof of Lemma \ref{lem:exact_scaled_moments}}
Retain $V(y,u)=y^2+m_0u^2$ with $m_0$ chosen as in Lemma~\ref{lem:diffusion_drift}, and recall $x,\iota$ from \eqref{eq:state_unscaled}. Applying $\calG_X$ to $V$ gives
\[
\calG_X V(y,u)
=2y[-\beta_\lambda+u-g(y)]-2m_0\kappa u^2
+\delta_\lambda^2\{\iota+d_\lambda(x)\}
+m_0\sigma^2(1+\eps_\lambda u),
\]
where the second-order jump contribution expands as
\[
\delta_\lambda^2\{\iota+d_\lambda(x)\}
=2\lambda^{-\alpha}
+\delta_\lambda\{u+\beta_\lambda+g(y)\}.
\]
Since $\lambda^{-\alpha},\delta_\lambda,\eps_\lambda\le1$, Young's inequality absorbs the linear terms into the negative quadratic drift to yield $\calG_X V\le -a_X V+b_X$ for constants $a_X,b_X>0$ independent of sufficiently large $\lambda$. For $p\ge1$ and $W=1+V$, applying $\calG_X$ to $W^p$ via Taylor's theorem for the queue jumps and the diffusion chain rule for the intensity coordinate bounds the remaining lower-order terms by $C_p W^{p-1/2}$, where $C_p$ depends only on $p$ and the fixed model parameters. Enlarging the additive constant therefore gives
\[
\calG_X W^p\le -\frac{a_X p}{2}W^p+b_p
\]
uniformly for sufficiently large $\lambda$; the same inequality holds with $\lambda$-dependent finite constants at every admissible $\lambda$.

To verify ergodicity at a fixed admissible $\lambda$, note first that pathwise continuity of $I_\lambda$ ensures almost surely finitely many arrivals on bounded time intervals; because departures cannot exceed arrivals plus the initial queue length, $(\wt X_\lambda,\wt I_\lambda)$ is nonexplosive, and the CIR Feller property together with locally bounded, continuous jump rates makes the joint process Feller. Next, any compact subset of $\mathsf{D}_\lambda$ contains only finitely many unscaled queue states $x$. On the event that the arrival intensity stays within a compact interval bounded away from zero while no arrivals occur and all initial customers depart, the positive CIR transition density restricted to such intensity paths and the positive departure rates yield a uniform minorization into $x=0$ paired with an interior intensity interval; combining this minorization with prescribed finite arrival and departure sequences establishes irreducibility and shows that interior compact sets are petite. Extending the state space to include the boundary $u=-\eps_\lambda^{-1}$ as in the proof of Lemma~\ref{lem:diffusion_drift} makes all compact sublevel sets of $V$ petite, so \citet[Theorems~4.2 and~4.3(ii)]{meyn1993stability} establish positive Harris recurrence, uniqueness of $\pi^X$ (with zero probability on $\{u=-\eps_\lambda^{-1}\}$), and the uniform bound \eqref{eq:exact_moment_bound} for $m\ge2$ and sufficiently large $\lambda$. Orders $1\le m<2$ follow by Jensen's inequality.
\Halmos
\endproof

\subsection{Stationary Inverse Moments of the Arrival Intensity}
\label{subsec:inverse-moment-proofs}

\begin{lemma}
\label{lem:cir_inverse_moments}
Under Assumptions~\ref{assump:model_primitives} and~\ref{assump:staffing_c}, write
\[
(\widehat X_\lambda,\widehat I_\lambda)
:=\bigl(\wt X_\lambda(\infty),\wt I_\lambda(\infty)\bigr),
\qquad
V_\lambda:=v_\lambda(\widehat I_\lambda).
\]
Then $V_\lambda\sim\operatorname{Gamma}(k_\lambda,1/k_\lambda)$, where $k_\lambda=2\kappa/(\sigma^2\eps_\lambda^2)\to\infty$ as $\lambda\to\infty$. For every fixed $q>0$ and all $\lambda$ satisfying $k_\lambda>q$,
\begin{equation*}
\E V_\lambda^{-q}
=k_\lambda^q\frac{\Gamma(k_\lambda-q)}{\Gamma(k_\lambda)}
\longrightarrow 1
\quad\text{as }\lambda\to\infty.
\end{equation*}
Consequently, for every fixed $m\ge0$ and $q>0$, there exists $\lambda_{m,q}<\infty$ such that
\begin{equation*}
\sup_{\lambda\ge\lambda_{m,q}}
\E\left[w(\widehat X_\lambda,\widehat I_\lambda)^m V_\lambda^{-q}\right]<\infty.
\end{equation*}
The same statements hold under $\pi^{\mathrm{CIR}}$.
\end{lemma}

\proof{Proof of Lemma~\ref{lem:cir_inverse_moments}}
Since $\eps_\lambda\delta_\lambda=\lambda^{-1}$, we have $V_\lambda=1+\eps_\lambda\widehat I_\lambda=I_\lambda(\infty)/\lambda$, so the stationary Gamma distribution of $I_\lambda(\infty)$ implies that $V_\lambda$ has density
\[
p_\lambda(v)=\frac{k_\lambda^{k_\lambda}}{\Gamma(k_\lambda)}
v^{k_\lambda-1}e^{-k_\lambda v},\qquad v>0.
\]
Integrating $v^{-q}$ against $p_\lambda(v)$ converges whenever $k_\lambda>q$ (and diverges for $k_\lambda\le q$), giving
\[
\E V_\lambda^{-q}
=\frac{k_\lambda^{k_\lambda}}{\Gamma(k_\lambda)}
\int_0^\infty v^{k_\lambda-q-1}e^{-k_\lambda v}\,\dd v
=k_\lambda^q\frac{\Gamma(k_\lambda-q)}{\Gamma(k_\lambda)}.
\]
Because $\Gamma(k-q)/\Gamma(k)=k^{-q}(1+O(k^{-1}))$ as $k\to\infty$ for fixed $q$, it follows that $\E V_\lambda^{-q}\to 1$ as $\lambda\to\infty$; in particular, every fixed inverse moment is uniformly bounded once $\lambda$ is sufficiently large.

To bound the joint expectation without assuming independence between $\widehat X_\lambda$ and $\widehat I_\lambda$, fix $m\ge0$ and $q>0$ and apply the Cauchy--Schwarz inequality:
\[
\E\left[w(\widehat X_\lambda,\widehat I_\lambda)^m V_\lambda^{-q}\right]
\le
\left(\E w(\widehat X_\lambda,\widehat I_\lambda)^{2m}\right)^{1/2}
\left(\E V_\lambda^{-2q}\right)^{1/2}.
\]
The first factor is uniformly bounded for large $\lambda$ by Lemma~\ref{lem:exact_scaled_moments} (bounding noninteger $2m$ by a larger integer power), and the second factor is uniformly bounded once $k_\lambda>2q$. Under $\pi^{\mathrm{CIR}}$, the normalized intensity $v_\lambda(Z_\lambda(\infty))$ shares the exact same $\operatorname{Gamma}(k_\lambda,1/k_\lambda)$ marginal while Lemma~\ref{lem:diffusion_drift} supplies the uniform polynomial moment bounds, so the identical argument applies.
\Halmos
\endproof

\subsection{Stationary Generator Identity}
\label{subsec:stationary-generator-identity}

The following lemma justifies taking stationary expectations in the Poisson equation once second-derivative growth bounds have been established.

\begin{lemma}
\label{lem:petite_generator_identity}
Under Assumptions~\ref{assump:model_primitives} and~\ref{assump:staffing_c}, let $\pi^X$ be the stationary distribution of $\widehat Z_\lambda(t):=\bigl(\wt X_\lambda(t),\wt I_\lambda(t)\bigr)$, and recall $\mathsf{D}_\lambda=\R\times(-\eps_\lambda^{-1},\infty)$ and $v_\lambda(u)=1+\eps_\lambda u$ from Section~\ref{sec:model}. Let $f=f_\lambda\in C^2(\mathsf{D}_\lambda)$, and suppose that there exist constants $p,q\ge0$ independent of $\lambda$ and a finite constant $K_\lambda$ (which may depend on $\lambda$ and $f$) such that
\begin{equation}
\max_{r,s\in\{y,u\}}|f_{rs}(y,u)|
\le K_\lambda w(y,u)^p\bigl(1+v_\lambda(u)^{-q}\bigr),
\qquad (y,u)\in\mathsf{D}_\lambda.
\label{eq:exact_generator_hessian_growth}
\end{equation}
Then, for all sufficiently large $\lambda$, both $f$ and $\calG_X f$ belong to $L^1(\pi^X)$, and
\begin{equation}
\int\calG_X f\,\dd\pi^X
=\E\calG_X f\bigl(\wt X_\lambda(\infty),\wt I_\lambda(\infty)\bigr)=0.
\label{eq:smooth_stein_generator_identity}
\end{equation}
In particular, \eqref{eq:smooth_stein_generator_identity} holds for either Poisson solution $f_h^D$, $D\in\{\mathrm{OU},\mathrm{CIR}\}$, once its $C^2$ regularity and Hessian bound \eqref{eq:exact_generator_hessian_growth} have been verified.
\end{lemma}

\proof{Proof of Lemma \ref{lem:petite_generator_identity}}
Fix a sufficiently large $\lambda$ and abbreviate $\pi=\pi^X$, $\delta=\delta_\lambda$, and $\omega(u)=1+v_\lambda(u)^{-q}$; throughout this proof, constants $C_\lambda$ may depend on $\lambda,f,K_\lambda,p,q$, and the fixed model parameters. Along the line segment from $(0,0)$ to $(y,u)\in\mathsf{D}_\lambda$, convexity gives
\[
w(sy,su)\le w(y,u),\qquad
v_\lambda(su)=(1-s)+s v_\lambda(u)
\ge\min\{1,v_\lambda(u)\},\quad 0\le s\le1.
\]
Integrating the Hessian bound \eqref{eq:exact_generator_hessian_growth} twice along this segment and absorbing $|f(0,0)|$ and $|\nabla f(0,0)|$ into $C_\lambda$ yields, with $\ell=p+2$,
\[
|f(y,u)|+|\nabla f(y,u)|+|f_{uu}(y,u)|
\le C_\lambda w(y,u)^\ell\omega(u).
\]
Next, write $\Delta_\pm f(y,u)=f(y\pm\delta,u)-f(y,u)$ for the unit-jump increments. Because queue-length jumps leave $u$ unchanged and satisfy $w(y+a,u)\le 2w(y,u)$ for $|a|\le\delta\le1$, the derivative bound implies
\[
|\Delta_\pm f(y,u)|
\le C_\lambda\delta w(y,u)^\ell\omega(u).
\]
Meanwhile, Lemmas~\ref{lem:exact_scaled_moments} and~\ref{lem:cir_inverse_moments} ensure that, for every fixed $r\ge0$ and all sufficiently large $\lambda$,
\[
\int w(y,u)^r\omega(u)^2\,\pi(\dd y,\dd u)<\infty.
\]
At any exact state $(y,u)$ with unscaled coordinates $(x,\iota)$ from \eqref{eq:state_unscaled} and departure rate $d=d_\lambda(x)$ (which vanishes at $x=0$), the transition rates satisfy $0\le\iota+d\le C_\lambda w(y,u)$ while $0<v_\lambda(u)\le w(y,u)$. Combining the formula for $\calG_X f$ with the preceding growth and integrability bounds gives
\begin{equation}
\int\bigl(|f|+|\calG_X f|\bigr)\,\dd\pi<\infty,
\label{eq:exact_generator_L1}
\end{equation}
and shows that the predictable quadratic-variation density
\begin{align}
Q_f(y,u)
&:=\iota\bigl(\Delta_+f(y,u)\bigr)^2+d\bigl(\Delta_-f(y,u)\bigr)^2
   +\sigma^2 v_\lambda(u)f_u(y,u)^2\nonumber\\
&\le C_\lambda w(y,u)^{2\ell+1}\omega(u)^2
\label{eq:exact_generator_quadratic_integrability}
\end{align}
also satisfies $\int Q_f\,\dd\pi<\infty$.

Now initialize $\widehat Z_\lambda$ with its stationary distribution $\pi$. Let $M_f(t)$ be the sum of the compensated arrival and departure stochastic integrals of $\Delta_+ f$ and $\Delta_- f$ and the Brownian stochastic integral of $\sigma\sqrt{v_\lambda}\,f_u$, each evaluated at $\widehat Z_\lambda(s-)$. Because stationarity gives $\E\int_0^t Q_f(\widehat Z_\lambda(s))\,\dd s = t\int Q_f\,\dd\pi<\infty$, the stochastic-integral isometries imply that $M_f$ is a zero-mean square-integrable martingale. Applying It\^o's formula for jump--diffusions up to the localization stopping times
\[
\tau_k=\inf\bigl\{t\ge0 : w(\widehat Z_\lambda(t))\ge k
\ \text{or}\ v_\lambda(\wt I_\lambda(t))\le k^{-1}\bigr\}
\]
gives the stopped version of
\begin{equation}
f(\widehat Z_\lambda(t))-f(\widehat Z_\lambda(0))
=\int_0^t\calG_X f(\widehat Z_\lambda(s))\,\dd s+M_f(t).
\label{eq:exact_generator_martingale_decomposition}
\end{equation}
Nonexplosion and nonattainment of the lower boundary $u=-\eps_\lambda^{-1}$ imply that $\tau_k\uparrow\infty$ almost surely as $k\to\infty$. As $k\to\infty$, the stopped drift integrals converge in $L^1$ by \eqref{eq:exact_generator_L1} and stationarity, and the stopped martingales converge in $L^2$ by \eqref{eq:exact_generator_quadratic_integrability}. Consequently, the stopped left-hand side converges in $L^1$ to its almost-sure limit $f(\widehat Z_\lambda(t))-f(\widehat Z_\lambda(0))$, establishing \eqref{eq:exact_generator_martingale_decomposition} without localization. Taking expectations in \eqref{eq:exact_generator_martingale_decomposition} and applying stationarity and Fubini's theorem yields $0=t\int\calG_X f\,\dd\pi$ for every $t>0$, which proves \eqref{eq:smooth_stein_generator_identity}.
\Halmos
\endproof

\section{Proof of Theorem~\ref{thm:smooth_bound_general_c} and Related OU Results}
\label{sec:main-result-proofs}

\subsection{Steady-State Error Bound for the OU-Based Diffusion}
\label{subsec:ou-main-proof}

\proof{Proof of Theorem \ref{thm:smooth_bound_general_c}.}
Following Section~\ref{subsec:generator_comparison}, write
\[
(\widehat X,\widehat I)
=
\bigl(\wt X_\lambda(\infty),\wt I_\lambda(\infty)\bigr),
\qquad
\delta=\delta_\lambda,
\qquad
r=\lambda^{-\alpha},
\qquad
\eps=\eps_\lambda,
\qquad
c=c_\lambda,
\]
and set $\widehat w=w(\widehat X,\widehat I)$, retaining $s(c)$ and $M(c)$ from Section~\ref{subsec:generator_comparison}. Let $f$ be the centered Poisson solution for $h$, satisfying $\calG_{\mathrm{OU}}f=h-\pi^{\mathrm{OU}}(h)$. Combining this Poisson equation with the stationary generator identity in Lemma~\ref{lem:petite_generator_identity} yields
\begin{align}
\E h(\widehat X,\widehat I)
-\E h\bigl(Y_\lambda^{\mathrm{OU}}(\infty),U(\infty)\bigr)
=
\E\left[
(\calG_{\mathrm{OU}}-\calG_X)
f(\widehat X,\widehat I)
\right].
\label{eq:ou-stein-identity-corrected}
\end{align}

By Lemma~\ref{lemma:stein-factor}, the second derivatives of $f$ satisfy
\[
|f_{yy}(y,u)| \le C H s(c) w(y,u)^{m-1},
\qquad
|f_{uu}(y,u)| \le C H w(y,u)^{m-1},
\]
while Lemma~\ref{lem:ou-fyy-lipschitz} gives the weighted increment bound
\begin{equation*}
|f_{yy}(v,u)-f_{yy}(y,u)|
\le C H M(c) |v-y|
(1+|v|+|y|+|u|)^m.
\end{equation*}
When $|v-y|\le\delta\le1$, the spatial factor $(1+|v|+|y|+|u|)^m$ is bounded by $2^m w(y,u)^m$; this common bound also covers the equal-rate case $\mu=\theta$, for which Lemma~\ref{lem:ou-fyy-lipschitz} provides the smaller exponent $m-2$. Moreover, Lemma~\ref{lem:exact_scaled_moments} implies that $\E\widehat w^\ell\le C_\ell$ for every fixed $\ell\ge0$ and all sufficiently large $\lambda$ (with noninteger orders controlled by larger integer moments). In particular, taking $\ell = m+1$ controls all queue-jump remainder and counting-variance expectations, so the hypotheses of \eqref{eq:taylor-minimum-corrected}--\eqref{eq:variance-mismatch-corrected} hold with $p=m$ and $\omega=1$. Finally, the $c$-uniform bound on $f_{uu}$ yields
\begin{align}
T_2
&:=
\frac{\sigma^2}{2}\eps
\E\left[
|\widehat I f_{uu}(\widehat X,\widehat I)|
\right]
\le
C H\eps\E\widehat w^m\le
C H\eps.
\label{eq:environmental-term-corrected}
\end{align}
Combining \eqref{eq:ou-stein-identity-corrected}, \eqref{eq:taylor-minimum-corrected}, \eqref{eq:local-second-order-corrected}, \eqref{eq:variance-mismatch-corrected}, and \eqref{eq:environmental-term-corrected}, we obtain
\begin{align}
\left|
\E h(\widehat X,\widehat I)
-\E h\bigl(Y_\lambda^{\mathrm{OU}}(\infty),U(\infty)\bigr)
\right|\le
C H
\left\{
r\min\{s(c),\delta\,M(c)\}
+
\left(
\delta+
\left|r-\frac{c^2}{2}\right|
\right)s(c)
+\eps
\right\}.
\label{eq:ou-general-bound-corrected-proof}
\end{align}

We now express the Taylor and local counting-variance terms in \eqref{eq:ou-general-bound-corrected-proof} in terms of the quantities appearing in the theorem. If $r\le c^2\le4r$, then $r\le1$, $\delta r^{-1/4}=\lambda^{-(2+\alpha)/4}$, and $\delta r^{-1/4}\le\eps$ for $\alpha\ge0$, which gives
\[
r\min\{s(c),\delta M(c)\}+\delta s(c)
\le C
\begin{cases}
\delta r^{-1/4},&\mu\ne\theta,\\
\delta,&\mu=\theta,
\end{cases}
\le C\eps.
\]
Conversely, if $c^2<r$ or $c^2>4r$, then $|r-c^2/2|\ge r/2$; using $\delta\le r$ together with the size alternative $s(c)$ in $\min\{s(c),\delta M(c)\}$ yields
\[
r\min\{s(c),\delta M(c)\}+\delta s(c)
\le C\left|r-\frac{c^2}{2}\right|s(c).
\]
Since $s(c)=1$ when $\mu=\theta$ and $s(c)\le Cc^{-1/2}$ when $\mu\ne\theta$ and $c\le\bar c$, substituting these bounds into \eqref{eq:ou-general-bound-corrected-proof} and absorbing the fixed factor $H$ into $C$ proves \eqref{eq:ou_bound_second_order_only}.

Under variance matching, $c=c_\lambda^*=\sqrt{2}\lambda^{-\alpha/2}$, we have $r=c^2/2$. Substituting the $r \le c^2 \le 4r$ estimate directly into \eqref{eq:ou-general-bound-corrected-proof} gives
\begin{align*}
\left|\E h\bigl(\wt X_\lambda(\infty),\wt I_\lambda(\infty)\bigr)
-\E h\bigl(Y_\lambda^{\mathrm{OU}}(\infty),U(\infty)\bigr)\right|
\le C H
\begin{cases}
\lambda^{-(2+\alpha)/4}+\eps_\lambda,&\mu\ne\theta,\\
\delta_\lambda+\eps_\lambda,&\mu=\theta,
\end{cases}
\end{align*}
where both branches are bounded by $C H\eps_\lambda$. Absorbing $H$ into the constant in the theorem proves \eqref{eq:ou_bound_variance_matched}.
\Halmos
\endproof

\subsection{Queue-Length Moments and Sharpness}
\label{subsec:ou-moment-proofs}

\proof{Proof of Corollary~\ref{cor:mean_variance}.}
Let $X=\widetilde X_\lambda(\infty)$ and $Y=Y_\lambda^{\mathrm{OU}}(\infty)$. Since the test functions $y$ and $y^2$ satisfy \eqref{eq:test_class} with $m=2$ and $H=2$, Theorem~\ref{thm:smooth_bound_general_c} implies
\[
 |\E X-\E Y|\le C\{\Psi_\lambda(c_\lambda)+\eps_\lambda\},
 \qquad
 |\E X^2-\E Y^2|\le C\{\Psi_\lambda(c_\lambda)+\eps_\lambda\}.
\]
Writing the variance difference as
\begin{equation}
 |\Var(X)-\Var(Y)|\le|\E X^2-\E Y^2|+
 |\E X-\E Y|\bigl(|\E X|+|\E Y|\bigr),
 \label{eq:variance-decomposition-corrected}
\end{equation}
and noting that $\sup_{\lambda\ge\lambda_0}(|\E X|+|\E Y|)<\infty$ by Lemmas~\ref{lem:diffusion_drift} and~\ref{lem:exact_scaled_moments}, we obtain \eqref{eq:ou_queue_moment_bound}. Under variance matching, $\Psi_\lambda(c_\lambda^*)=0$.
\Halmos
\endproof

\proof{Proof of Proposition~\ref{prop:ou-sharpness}.}
Write $V = I_\lambda(\infty)/\lambda$ and $k = 2\kappa/(\sigma^2\eps_\lambda^2)$. Under Assumption~\ref{assump:model_primitives}, $V$ follows a Gamma stationary distribution with shape $k$ and scale $1/k$, so $\E[(V - 1)^3] = 2/k^2$. Since $\wt I_\lambda(\infty) = \eps_\lambda^{-1}(V - 1)$, we obtain
\[
\E\bigl[\wt I_\lambda(\infty)^3\bigr]
= \eps_\lambda^{-3}\E\bigl[(V - 1)^3\bigr]
= \frac{\sigma^4}{2\kappa^2}\eps_\lambda,
\qquad
\E\bigl[U(\infty)^3\bigr] = 0,
\]
where the second identity follows from the symmetry of the Gaussian distribution. Because $h(y, u) = u^3$ satisfies \eqref{eq:test_class} for every $m \ge 3$ with $H = 6$, rescaling $h$ by $\min\{1, H/6\}$ yields an admissible test function attaining the exact rate $\eps_\lambda$ for any $H > 0$.

For $m=2$, take $h(y,u)=H\sin u$, which satisfies all three inequalities in \eqref{eq:test_class}. Evaluating the characteristic function of the Gamma distribution of $V$ gives
\[
\E\sin\bigl(\wt I_\lambda(\infty)\bigr)
=\left(1+\frac{1}{k^2\eps_\lambda^2}\right)^{-k/2}
 \sin\left(k\arctan\frac{1}{k\eps_\lambda}
              -\eps_\lambda^{-1}\right).
\]
Substituting $k=2\kappa/(\sigma^2\eps_\lambda^2)$ shows that the prefactor outside the sine equals $e^{-\sigma^2/(4\kappa)}+O(\eps_\lambda^2)$ while the sine argument expands as $-\sigma^4\eps_\lambda/(12\kappa^2)+O(\eps_\lambda^3)$. Since $\E\sin(U(\infty))=0$ by Gaussian symmetry, it follows that
\[
|\pi^X(h)-\pi^{\mathrm{OU}}(h)|
=\frac{H\sigma^4}{12\kappa^2}e^{-\sigma^2/(4\kappa)}
  \eps_\lambda+O(\eps_\lambda^3),
\]
whose leading coefficient is strictly positive and independent of $c_\lambda$. This establishes sharpness for $m=2$ as well.
\Halmos
\endproof

\subsection{Time Spent near Full Occupancy and Second-Derivative Bounds}
\label{subsec:ou-hessian-proofs}

\begin{lemma}
\label{lemma:stein-factor}
Suppose that $\mu,\theta,\kappa,\sigma>0$, $B_\beta<\infty$, and $0<c_\lambda\le \bar c$, and let $P_t$ and $\pi$ denote the semigroup and stationary distribution of the OU diffusion in \eqref{eq:ou_diffusion_general_c}, respectively. Fix an integer $m\ge2$, and assume that $h\in C^2(\mathbb R^2)$ satisfies, for some $H<\infty$,
\[
|h(y,u)|\le H w(y,u)^m,
\qquad
|\nabla h(y,u)|\le H w(y,u)^{m-1},
\qquad
\|\nabla^2 h(y,u)\|\le H w(y,u)^{m-2}.
\]
Define the centered Poisson solution by
\[
f(y,u)
=
-\int_0^\infty
\left(P_t h(y,u)-\pi(h)\right)\,dt.
\]
Then $f\in C^2(\mathbb R^2)$. Moreover, there exists a constant $C<\infty$, allowed to depend on $m$ but independent of $\lambda,c_\lambda,h,y$, and $u$, such that
\begin{align}
\left|f_{yy}(y,u)\right|
&\le
C H
\left(1+|\theta-\mu|\,c_\lambda^{-1/2}\right)
w(y,u)^{m-1},
\label{eq:stein-factor-yy}\\
\left|f_{yu}(y,u)\right|
+
\left|f_{uu}(y,u)\right|
&\le
C H w(y,u)^{m-1}.
\label{eq:stein-factor-yu-uu}
\end{align}
\end{lemma}

\proof{Proof of Lemma \ref{lemma:stein-factor}}
Fix $\lambda$ and write $c=c_\lambda$, $\beta=\beta_\lambda$, $z=(y,u)$, and $\gamma=\tfrac12\min\{a_0,\kappa\}$. Throughout the proof, constants may depend on the fixed model parameters, $B_\beta,\bar c$, and the indicated fixed moment orders, but not on $\lambda,c,h$, or the smoothing parameter.

Choose an even, nonnegative mollifier $\rho\in C_c^\infty((-1,1))$ with $\int\rho=1$, set $\rho_\varepsilon(r)=\varepsilon^{-1}\rho(r/\varepsilon)$, and smooth the queue-length drift via $g_\varepsilon=g*\rho_\varepsilon$ for $0<\varepsilon\le1$. Then
\begin{equation}
 a_0\le g_\varepsilon'\le a_1,\qquad
 g_\varepsilon''=(\theta-\mu)\rho_\varepsilon,\qquad
 \|g_\varepsilon-g\|_\infty\le a_1\varepsilon.
 \label{eq:stein-completion-smoothing}
\end{equation}
Let $h_\varepsilon$ be the smooth convolution of $h$ with a nonnegative unit-mass kernel supported in the ball of radius $\varepsilon$, so that $h_\varepsilon$, $\nabla h_\varepsilon$, and $\nabla^2 h_\varepsilon$ satisfy the polynomial bounds with weights $w^m,w^{m-1},w^{m-2}$ and constant $CH$, and $|h_\varepsilon-h|\le CH\varepsilon w^{m-1}$. Write $Z_t^\varepsilon=(Y_t^\varepsilon,U_t)$ for the diffusion with smoothed drift $g_\varepsilon$, and let $P_t^\varepsilon$ and $\pi_\varepsilon$ denote its semigroup and stationary distribution. Because $|g_\varepsilon(0)|\le a_1$ and $y[g_\varepsilon(y)-g_\varepsilon(0)]\ge a_0y^2$, the Lyapunov argument in Lemma~\ref{lem:diffusion_drift} applies uniformly in $\varepsilon\in(0,1]$, giving, for every fixed $r\ge0$,
\begin{equation*}
 \mathbb E_z w(Z_t^\varepsilon)^r\le C_r w(z)^r,\qquad
 \mathbb E[w(Z_t^\varepsilon)^r\mid\mathcal F_s]
 \le C_r w(Z_s^\varepsilon)^r\quad(s\le t),\qquad
 \pi_\varepsilon(w^r)\le C_r,
\end{equation*}
with the same bounds holding for the unsmoothed diffusion.

Coupling two trajectories from initial states $z=(y,u)$ and $z'=(y',u')$ with the same Brownian motions gives $U_t^u-U_t^{u'}=e^{-\kappa t}(u-u')$, while the secant-slope bounds in \eqref{eq:stein-completion-smoothing} imply
\[
 |Y_t^{\varepsilon,z}-Y_t^{\varepsilon,z'}|
 \le e^{-a_0t}|y-y'|
 +|u-u'|\int_0^t e^{-a_0(t-s)}e^{-\kappa s}\,ds.
\]
Hence $|Z_t^{\varepsilon,z}-Z_t^{\varepsilon,z'}|\le Ce^{-\gamma t}|z-z'|$ uniformly in $\varepsilon$. Drawing $z'$ from $\pi_\varepsilon$ and applying the derivative bound on $h_\varepsilon$ together with the order-$m$ moment estimates yields
\[
 |P_t^\varepsilon h_\varepsilon(z)-\pi_\varepsilon(h_\varepsilon)|
 \le CH e^{-\gamma t}w(z)^m,
\]
which ensures absolute convergence of the centered Poisson integrals
\[
 f_\varepsilon(z)=-\int_0^\infty
 [P_t^\varepsilon h_\varepsilon(z)-\pi_\varepsilon(h_\varepsilon)]\,dt
 \quad\hbox{and}\quad
 f(z)=-\int_0^\infty[P_th(z)-\pi(h)]\,dt,
\]
along with the pointwise bound $|f_\varepsilon(z)|, |f(z)| \le CHw(z)^m$. For fixed $\varepsilon,c>0$, the semigroup identity
\[
 P_s^\varepsilon f_\varepsilon-f_\varepsilon
 =\int_0^s[P_t^\varepsilon h_\varepsilon
          -\pi_\varepsilon(h_\varepsilon)]\,dt
\]
implies that $f_\varepsilon$ solves the smoothed Poisson equation in distributions, and local elliptic regularity yields a smooth version of $f_\varepsilon$. Moreover, synchronous coupling of $Z^\varepsilon$ and $Z$ gives
\begin{equation*}
 \sup_{t\ge0}|Y_t^\varepsilon-Y_t|
 \le \frac{a_1}{a_0}\varepsilon,\qquad U_t^\varepsilon=U_t,
\end{equation*}
so the bounds on $h_\varepsilon$ and $h$ imply $|P_t^\varepsilon h_\varepsilon(z)-P_th(z)|\le CH\varepsilon w(z)^{m-1}$ uniformly in $t\ge0$. Evaluating at $z=0$ and letting $t\to\infty$ gives $|\pi_\varepsilon(h_\varepsilon)-\pi(h)|\le CH\varepsilon$; combining this finite-time convergence with the uniform $O(e^{-\gamma t})$ tail bound yields $f_\varepsilon\to f$ locally uniformly as $\varepsilon\downarrow0$.

We next establish the short-time derivative bounds and occupation estimates needed for both the Hessian bounds and the queue-direction increments. Define the short-time kernel
\begin{equation}
 k_c(t)=
 \begin{cases}
 (c^2t+t^3)^{-1/2},&0<t\le1,\\
 1,&t>1.
 \end{cases}
 \label{eq:ou-smoothing-kernel}
\end{equation}
For $0<t\le1$, set $A_s=\exp\{-\int_0^s g_\varepsilon'(Y_r^\varepsilon)\,dr\}$. Perturbing the driving Brownian motions via $dB_i(s)+\eta v_i(s)\,ds$ with an adapted control $v=(v_1,v_2)$ induces the first variation $\mathcal V=(\mathcal V^Y,\mathcal V^U)$ solving
\[
 \dot{\mathcal V}_s^Y
 =-g_\varepsilon'(Y_s^\varepsilon)\mathcal V_s^Y
   +\mathcal V_s^U+cv_1(s),\qquad
 \dot{\mathcal V}_s^U=-\kappa\mathcal V_s^U+\sigma v_2(s),
 \qquad \mathcal V_0=0.
\]
To match the initial-queue tangent $\partial_yZ_t^{\varepsilon,z}=(A_t,0)$ at time $t$, we construct two adapted controls. The first perturbs only $B_1$:
\[
 v_1(s)=\frac{A_s}{ct},\qquad v_2(s)=0,\qquad
 \mathcal V_s^Y=\frac{s}{t}A_s,\quad\mathcal V_s^U=0,
 \qquad \int_0^t|v(s)|^2\,ds\le(c^2t)^{-1}.
\]
The second perturbs only $B_2$ using $\ell(s)=3(s/t)^2-2(s/t)^3$:
\[
 v_1(s)=0,\qquad
 v_2(s)=\frac{A_s}{\sigma}
 [\ell''(s)+(\kappa-g_\varepsilon'(Y_s^\varepsilon))\ell'(s)],
\]
for which $\mathcal V_s^Y=\ell(s)A_s$ and $\mathcal V_s^U=\ell'(s)A_s$ solve the variation system, the endpoint conditions on $\ell$ give $\mathcal V_t=(A_t,0)$, and the bounds on $g_\varepsilon'$ yield $\int_0^t|v(s)|^2\,ds\le Ct^{-3}$. Because both controls are adapted and square integrable, Malliavin integration by parts \citep{nualart2006malliavin} gives, initially for $\psi\in C_c^\infty(\mathbb R^2)$,
\begin{equation*}
 (P_t^\varepsilon\psi)_y(z)
 =\mathbb E_z\left[
 \psi(Z_t^\varepsilon)\int_0^t v(s)\cdot dB(s)\right].
\end{equation*}
Applying Cauchy--Schwarz, the order-$2r$ moment bound, and the sharper of the two control-energy bounds (using $\min\{c^{-1}t^{-1/2},t^{-3/2}\}\le\sqrt2\,k_c(t)$ on $(0,1]$) shows that whenever $|\psi|\le K w^r$,
\begin{equation}
 |(P_t^\varepsilon\psi)_y(z)|
 \le C_r K w(z)^r k_c(t),\qquad t>0.
 \label{eq:ou-lip-semigroup-gradient}
\end{equation}
Here smooth truncation extends the bound to all continuous $\psi$ with $|\psi|\le K w^r$ via $L^2$ convergence under the transition distribution, while the case $t>1$ follows from $P_t^\varepsilon=P_1^\varepsilon P_{t-1}^\varepsilon$ and the transient moment bound.

The same pair of controls yields a short-time occupation bound without any additional construction. For nonnegative $\varphi\in C_c^\infty(\mathbb R)$, set $\Psi(x)=\int_{-\infty}^x\varphi(a)\,da$. Because the terminal arrival-intensity variation $\mathcal V_t^U = D_v U_t$ vanishes under either control,
\[
 D_v\left[\Psi(Y_t^\varepsilon)(1+U_t^2)^{r/2}\right]
 =\varphi(Y_t^\varepsilon)A_t(1+U_t^2)^{r/2}.
\]
Using $A_t\ge e^{-a_1t}$ and $0\le\Psi\le\|\varphi\|_1$, Malliavin duality and Cauchy--Schwarz imply
\[
 \begin{aligned}
 \mathbb E_z[(1+U_t^2)^{r/2}\varphi(Y_t^\varepsilon)]
 &\le e^{a_1t}\|\varphi\|_1
 \left(\mathbb E_z(1+U_t^2)^r\right)^{1/2}
 \left(\mathbb E_z\int_0^t|v(s)|^2\,ds\right)^{1/2}\\
 &\le C_r(1+|u|)^r\|\varphi\|_1
 \min\{c^{-1}t^{-1/2},t^{-3/2}\},
 \qquad 0<t\le1,
 \end{aligned}
\]
where the explicit OU representation bounds $\mathbb E_z(1+U_t^2)^r$ by $C_r(1+|u|)^{2r}$ and justifies truncating the terminal weight (whose Malliavin derivative along $v$ remains zero since $D_vU_t=0$). Density extends this inequality to every nonnegative $\varphi\in L^1(\mathbb R)$, and conditioning on $\mathcal F_{t-1}$ over the final unit interval extends it to $t>1$, giving, for every fixed $r\ge0$,
\begin{equation}
 \mathbb E_{y,u}[(1+|U_t|)^r\varphi(Y_t^\varepsilon)]
 \le C_r(1+|u|)^r k_c(t)\|\varphi\|_1,\qquad t>0,
 \label{eq:ou-sharp-density}
\end{equation}
with constants independent of $\varepsilon$.

To bound the Hessian of $f_\varepsilon$, write $J_a(t)=\partial_a Z_t^{\varepsilon,z}$, $j_a(t)=\partial_aY_t^{\varepsilon,z}$, and $A(s,t)=\exp\{-\int_s^t g_\varepsilon'(Y_r^\varepsilon)\,dr\}$ for $a,b\in\{y,u\}$. The first-order flow variations
\[
 j_y(t)=A(0,t),\qquad
 j_u(t)=\int_0^t A(s,t)e^{-\kappa s}\,ds,\qquad
 \partial_yU_t=0,\quad\partial_uU_t=e^{-\kappa t},
\]
satisfy $|J_a(t)|\le Ce^{-\gamma t}$ and $|j_a(t)|\le Ct^{\mathbf1_{\{a=u\}}}e^{-\gamma t}$, whereas the second variations of $U_t$ vanish and
\begin{equation*}
 \partial_{ab}Y_t^\varepsilon
 =-(\theta-\mu)\int_0^t
 A(s,t)\rho_\varepsilon(Y_s^\varepsilon)j_a(s)j_b(s)\,ds.
\end{equation*}
Setting $m_{ab}=\mathbf1_{\{a=u\}}+\mathbf1_{\{b=u\}}$ and differentiating under the expectation (justified at fixed $\varepsilon>0$ by boundedness of $g_\varepsilon',g_\varepsilon''$ and the polynomial moment bounds) gives
\begin{equation}
 (P_t^\varepsilon h_\varepsilon)_{ab}(z)
 =\mathbb E_z\left[
 J_a(t)^\top\nabla^2h_\varepsilon(Z_t^\varepsilon)J_b(t)
 +(h_\varepsilon)_y(Z_t^\varepsilon)\partial_{ab}Y_t^\varepsilon
 \right].
 \label{eq:stein-completion-semigroup-hessian}
\end{equation}
The first term in \eqref{eq:stein-completion-semigroup-hessian} is bounded in magnitude by $CH e^{-2\gamma t}w(z)^{m-2}$. For the second term, conditioning on $\mathcal F_s$ and noting that $|Y_s^\varepsilon|\le1$ on the support of $\rho_\varepsilon(Y_s^\varepsilon)$ gives, via \eqref{eq:ou-sharp-density} and $\|\rho_\varepsilon\|_1=1$,
\begin{equation*}
 \begin{aligned}
 \mathbb E_z[w(Z_t^\varepsilon)^{m-1}
                    \rho_\varepsilon(Y_s^\varepsilon)]
 &\le C\mathbb E_z[w(Z_s^\varepsilon)^{m-1}
                    \rho_\varepsilon(Y_s^\varepsilon)]\\
 &\le C(1+|u|)^{m-1}k_c(s),\qquad 0<s\le t.
 \end{aligned}
\end{equation*}
Combining the flow bounds and integrating first over $t\ge s$ (using $\int_s^\infty e^{-a_0(t-s)}\,dt=a_0^{-1}$) yields
\begin{equation}
 \begin{aligned}
 \int_0^\infty |(P_t^\varepsilon h_\varepsilon)_{ab}(z)|\,dt
 \le CH\bigg[
 w(z)^{m-2}
 +|\theta-\mu|(1+|u|)^{m-1}
 \int_0^\infty s^{m_{ab}}e^{-2\gamma s}k_c(s)\,ds
 \bigg].
 \end{aligned}
 \label{eq:stein-completion-integrated-semigroup}
\end{equation}
On $(0,1]$, substituting $s=cr$ and splitting the integral at $r=1$ gives
\begin{equation}
 \int_0^1 k_c(s)\,ds
 =c^{-1/2}\int_0^{1/c}(r+r^3)^{-1/2}\,dr
 \le4c^{-1/2},
 \label{eq:stein-completion-Ic0}
\end{equation}
whereas for $\ell=1,2$ we have $s^\ell k_c(s)\le s^{\ell-3/2}$ on $(0,1]$. Since the exponentially damped integrals over $(1,\infty)$ are $O(1)$ and $c\le\bar c$, we obtain
\begin{equation}
 \int_0^\infty e^{-2\gamma s}k_c(s)\,ds\le Cc^{-1/2},
 \qquad
 \int_0^\infty s^\ell e^{-2\gamma s}k_c(s)\,ds\le C
 \quad(\ell=1,2).
 \label{eq:stein-completion-scaling}
\end{equation}
Substituting \eqref{eq:stein-completion-scaling} into \eqref{eq:stein-completion-integrated-semigroup} justifies differentiating the finite- and infinite-horizon Poisson integrals and yields
\begin{equation*}
 |(f_\varepsilon)_{yy}(z)|
 \le CH(1+|\theta-\mu|c^{-1/2})w(z)^{m-1},\qquad
 |(f_\varepsilon)_{yu}(z)|+|(f_\varepsilon)_{uu}(z)|
 \le CHw(z)^{m-1},
\end{equation*}
while the first-variation formula similarly gives $|\nabla f_\varepsilon(z)|\le CHw(z)^{m-1}$ uniformly in $\varepsilon$.

Finally, since $f_\varepsilon\to f$ locally uniformly, integrating by parts against smooth compactly supported test functions passes these uniform derivative bounds to the weak derivatives of $f$. Moreover, in the identity
\[
 \mathcal G_{\mathrm{OU}}f_\varepsilon
 =h_\varepsilon-\pi_\varepsilon(h_\varepsilon)
 +(g_\varepsilon-g)(f_\varepsilon)_y,
\]
the remainder $(g_\varepsilon-g)(f_\varepsilon)_y$ converges locally uniformly to zero as $\varepsilon\downarrow0$, so $\mathcal G_{\mathrm{OU}}f=h-\pi(h)$ holds in distributions. Local boundedness of the weak Hessian, Sobolev embedding, and local elliptic regularity for fixed $c>0$ imply $f\in C^2(\mathbb R^2)$, so its continuous second derivatives satisfy \eqref{eq:stein-factor-yy}--\eqref{eq:stein-factor-yu-uu} pointwise.
\Halmos
\endproof

\subsection{Queue-Direction Increments of $f_{yy}$}
\label{subsec:ou-increment-proofs}

\begin{lemma}
\label{lem:ou-fyy-lipschitz}
Under the assumptions of Lemma~\ref{lemma:stein-factor}, there exists $C<\infty$, allowed to depend on the fixed degree $m$ but independent of $\lambda,c_\lambda$, and $h$, such that, for every $x,y,u\in\mathbb R$,
\begin{align}
\left|f_{yy}(y,u)
      -f_{yy}(x,u)\right|\le
C H c_\lambda^{-5/2}|y-x|
\bigl(1+|x|+|y|+|u|\bigr)^m.
\label{eq:stein-factor-yy-modulus}
\end{align}
In particular, for every fixed $u$, the map $y\mapsto f_{yy}(y,u)$ is locally Lipschitz. Moreover, $f_{yy}$ admits a weak derivative in the queue coordinate, \(q:=D_y^{\mathrm w}\bigl(f_{yy}\bigr)\in L^\infty_{\mathrm{loc}}(\mathbb R^2)\) (where $L^\infty_{\mathrm{loc}}$ denotes functions essentially bounded on compact sets and $D_y^{\mathrm w}$ denotes weak differentiation in $y$), satisfying
\[
  \int_{\mathbb R^2} q(z)\varphi(z)\,dz
  =
  -\int_{\mathbb R^2}
    f_{yy}(z)\,\varphi_y(z)\,dz,
  \qquad
  \varphi\in C_c^\infty(\mathbb R^2),
\]
and
\begin{equation}
  |q(y,u)|
  \le C H c_\lambda^{-5/2}w(y,u)^m
  \quad\text{for almost every }(y,u)\in\mathbb R^2.
\label{eq:ou-fyyy-weak}
\end{equation}
If $\mu=\theta$, the sharper estimate holds:
\begin{equation}
  \left|f_{yy}(y,u)
        -f_{yy}(x,u)\right|
  \le C H(1+c_\lambda^{-1/2})|y-x|
  \bigl(1+|x|+|y|+|u|\bigr)^{m-2},
  \qquad x,y,u\in\mathbb R.
\label{eq:stein-factor-yy-modulus-linear}
\end{equation}
\end{lemma}

\proof{Proof of Lemma \ref{lem:ou-fyy-lipschitz}}
Retain the smoothed drift $g_\varepsilon$, test function $h_\varepsilon$, and centered Poisson solution $f_\varepsilon$ from the proof of Lemma~\ref{lemma:stein-factor}, and write $c=c_\lambda$, $\beta=\beta_\lambda$, and $b_\varepsilon(y,u)=-\beta+u-g_\varepsilon(y)$. As shown in that proof, uniformly in $\varepsilon\in(0,1]$,
\begin{equation*}
 \begin{aligned}
 |(f_\varepsilon)_y(z)|&\le CHw(z)^{m-1},\\
 |(f_\varepsilon)_{yy}(z)|
 &\le CH(1+|\theta-\mu|c^{-1/2})w(z)^{m-1},\\
 |(f_\varepsilon)_{yu}(z)|+|(f_\varepsilon)_{uu}(z)|
 &\le CHw(z)^{m-1}.
 \end{aligned}
\end{equation*}

For fixed $\varepsilon,c>0$, differentiating the smoothed Poisson equation $\mathcal L_\varepsilon f_\varepsilon = h_\varepsilon - \pi_\varepsilon(h_\varepsilon)$ twice in $u$ yields
\begin{equation*}
 (\mathcal L_\varepsilon-2\kappa)(f_\varepsilon)_{uu}
 =F_\varepsilon,\qquad
 F_\varepsilon=(h_\varepsilon)_{uu}-2(f_\varepsilon)_{yu},
 \quad |F_\varepsilon|\le CHw^{m-1},
\end{equation*}
where $\mathcal L_\varepsilon$ is the smoothed OU generator. At each fixed $\varepsilon>0$, smooth-flow differentiation provides polynomial growth bounds on the higher derivatives of $f_\varepsilon$, while the globally Lipschitz drift and constant diffusion matrix yield polynomial moments of the path maximum over finite horizons. These bounds justify removing localization in discounted Dynkin's formula and letting the terminal horizon tend to infinity via the exponential discount $e^{-2\kappa t}$, giving
\[
 (f_\varepsilon)_{uu}(z)
 =-\int_0^\infty e^{-2\kappa t}
 P_t^\varepsilon F_\varepsilon(z)\,dt.
\]
Differentiating under the integral via the queue-direction derivative bound \eqref{eq:ou-lip-semigroup-gradient} and the kernel integrals \eqref{eq:stein-completion-Ic0}--\eqref{eq:stein-completion-scaling} then gives
\begin{equation*}
 |(f_\varepsilon)_{yuu}(z)|
 \le CH(1+c^{-1/2})w(z)^{m-1}.
\end{equation*}
Next, differentiating the smoothed Poisson equation once in $y$ and rearranging yields
\begin{equation}
 \frac{c^2}{2}(f_\varepsilon)_{yyy}
 =(h_\varepsilon)_y+g_\varepsilon'(f_\varepsilon)_y
 -b_\varepsilon(f_\varepsilon)_{yy}
 +\kappa u(f_\varepsilon)_{yu}
 -\frac{\sigma^2}{2}(f_\varepsilon)_{yuu}.
 \label{eq:ou-lip-fyyy-equation}
\end{equation}
Since $|b_\varepsilon(z)|\le Cw(z)$, $|g_\varepsilon'|\le a_1$, and $0<c\le\bar c$, combining the preceding derivative bounds gives
\begin{equation}
 |(f_\varepsilon)_{yyy}(z)|
 \le CHc^{-2}(1+c^{-1/2})w(z)^m
 \le CHc^{-5/2}w(z)^m.
 \label{eq:ou-lip-smooth-fyyy}
\end{equation}

Because $f_\varepsilon\to f$ locally uniformly and $f\in C^2(\mathbb R^2)$, integrating by parts against any $\varphi\in C_c^\infty(\mathbb R^2)$ gives
\[
 \left|-\int_{\mathbb R^2}f_{yy}(z)\varphi_y(z)\,dz\right|
 =\left|\lim_{\varepsilon\downarrow0}
 \int_{\mathbb R^2}(f_\varepsilon)_{yyy}(z)\varphi(z)\,dz\right|
 \le CHc^{-5/2}\int_{\mathbb R^2}w(z)^m|\varphi(z)|\,dz.
\]
By $L^1$--$L^\infty$ duality on compact sets, $D_y^{\mathrm w}f_{yy}$ is represented by a function $q\in L^\infty_{\mathrm{loc}}(\mathbb R^2)$ satisfying $|q|\le CHc^{-5/2}w^m$ almost everywhere, proving \eqref{eq:ou-fyyy-weak}. By the slicing characterization of weak derivatives, for almost every $u\in\mathbb R$, $f_{yy}(\,\cdot\,,u)$ admits an absolutely continuous representative with derivative $q(\,\cdot\,,u)$ almost everywhere; continuity of $f_{yy}$ identifies $f_{yy}(\,\cdot\,,u)$ with that representative everywhere along each such slice, so the fundamental theorem of calculus gives
\[
 |f_{yy}(y,u)-f_{yy}(x,u)|
 \le CHc^{-5/2}|y-x|(1+|x|+|y|+|u|)^m.
\]
Approximating any $u$ in the measure-zero exceptional set by nonexceptional values and invoking joint continuity of $f_{yy}$ extends this bound to all $x,y,u\in\mathbb R$, proving \eqref{eq:stein-factor-yy-modulus}.

Finally, when $\mu=\theta$, the queue tangent reduces to $e^{-\mu t}$ and its second variation vanishes identically, so
\[
 (f_\varepsilon)_{yy}(z)
 =-\int_0^\infty e^{-2\mu t}
 P_t^\varepsilon(h_\varepsilon)_{yy}(z)\,dt.
\]
Applying \eqref{eq:ou-lip-semigroup-gradient} with $r=m-2$ together with $|(h_\varepsilon)_{yy}|\le CHw^{m-2}$ yields $|(f_\varepsilon)_{yyy}(z)|\le CH(1+c^{-1/2})w(z)^{m-2}$, and the same weak-limit, slicing, and continuity argument establishes \eqref{eq:stein-factor-yy-modulus-linear} (including the spatially uniform bound at $m=2$).
\Halmos
\endproof

\section{Proof of Theorem~\ref{thm:cir_smooth_bound} and Related CIR Results}
\label{sec:cir-result-proofs}

\subsection{Steady-State Error Bound for the CIR-Preserving Diffusion}
\label{subsec:cir-proofs}

\proof{Proof of Theorem~\ref{thm:cir_smooth_bound}.}
Throughout this proof, $C$ may depend on the fixed growth order $m \ge 2$ and the fixed model parameters, but is independent of $\lambda$, $h$, and $H$; we keep the factor $H$ explicit until the final theorem bounds. Let $P_t^{\mathrm{CIR}}$ denote the CIR-preserving semigroup and let $f$ be its centered Poisson solution,
\begin{equation*}
\begin{aligned}
f(y, u) &= -\int_0^\infty \bigl(P_t^{\mathrm{CIR}}h(y, u) - \pi^{\mathrm{CIR}}(h)\bigr)\,\dd t, \\
\calG_{\mathrm{CIR}}f &= h - \pi^{\mathrm{CIR}}(h).
\end{aligned}
\end{equation*}

Retain $\delta, r, c, s(c), M(c)$ from Section~\ref{subsec:generator_comparison}, put $\Gamma = |r - c^2/2|$, and let $(\widehat X, \widehat I)$ have distribution $\pi^X$. Lemma~\ref{lem:cir_queue_stein} establishes the second-derivative and queue-increment bounds for $f$ with fixed exponents $p, q$ and boundary weight $\omega(u) = 1 + v_\lambda(u)^{-q}$. Setting $\widehat w = w(\widehat X, \widehat I)$ and $\widehat\omega = \omega(\widehat I)$, Lemma~\ref{lem:cir_inverse_moments} implies that, for every fixed $\ell \ge 0$,
\begin{equation}
\sup_{\lambda \ge \lambda_\ell} \E[\widehat w^\ell \widehat\omega] < \infty.
\label{eq:cir-proof-weight-integrability}
\end{equation}

Because customer arrivals and departures leave $\widehat I$ unchanged, the same boundary weight $\widehat\omega$ appears across both Taylor increments. Taking $\ell = p + 1$ in \eqref{eq:cir-proof-weight-integrability} absorbs the linear rate factor in \eqref{eq:cir-proof-rate-bound}, so \eqref{eq:taylor-minimum-corrected} yields
\begin{equation*}
T_+ + T_- \le C H r \min\{s(c), \delta M(c)\}.
\end{equation*}
Meanwhile, \eqref{eq:local-second-order-corrected} and \eqref{eq:variance-mismatch-corrected} bound the two counting-variance terms by $C H(\delta + \Gamma)s(c)$, whereas the arrival-intensity generators match identically in \eqref{eq:gen_diff_cir_general_c}.

Using $r \le \Gamma + c^2/2$ together with $c^2 M(c) \le C s(c)$ gives
\[
r \min\{s(c), \delta M(c)\}
\le \Gamma s(c) + \frac{c^2}{2}\delta M(c)
\le C(\delta + \Gamma)s(c).
\]
Combining the Poisson equation with Lemma~\ref{lem:petite_generator_identity} and absorbing the fixed factor $H$ into the theorem constant proves \eqref{eq:cir-general-c-bound}. Under variance matching, $\Gamma = 0$, and substituting $c = \sqrt{2}\,\lambda^{-\alpha/2}$ and $\delta = \lambda^{-(1+\alpha)/2}$ yields \eqref{eq:cir-smooth-bound}.
\Halmos
\endproof

\subsection{Queue-Length Moments and Exact Moments at Equal Rates}
\label{subsec:cir-moment-proofs}

\proof{Proof of Corollary~\ref{cor:cir_moments}.}
Applying Theorem~\ref{thm:cir_smooth_bound} with $c_\lambda = c_\lambda^*$ to the test functions $y$ and $y^2$ (both of which satisfy \eqref{eq:test_class} with $m = H = 2$) and substituting the resulting bounds into \eqref{eq:variance-decomposition-corrected} with $X = \widetilde X_\lambda(\infty)$ and $Y = Y_\lambda^{\mathrm{CIR}}(\infty)$, together with the uniform first-moment bounds in Lemmas~\ref{lem:diffusion_drift} and~\ref{lem:exact_scaled_moments}, yields \eqref{eq:cir-moment-bound}.
\Halmos
\endproof

\proof{Proof of Proposition~\ref{prop:cir-exact-moments}.}
Fix an admissible $\lambda$ and put $\mu = \theta = a$. At this fixed $\lambda$, the polynomial drift estimates in the proofs of Lemmas~\ref{lem:diffusion_drift} and~\ref{lem:exact_scaled_moments} hold with constants allowed to depend on $\lambda$: the queue jump rates and arrival-intensity variance rate grow at most linearly, while the linear queue and intensity drifts give $\mathcal G V \le -a_V V + b_V$ for the quadratic Lyapunov function $V$. Applying the same estimates to powers of $1 + V$ provides finite polynomial moments of all orders needed below. Moreover, for test polynomials of total degree at most two, both the generators and the quadratic-variation rates of the martingale terms have polynomial growth, so localization and stationarity yield the exact stationary generator identities at this $\lambda$.

Let $(Y_\infty, Z_\infty)$ denote any of $(\widehat X_\lambda, \widehat I_\lambda)$, $(Y_\lambda^{\mathrm{OU}}(\infty), U(\infty))$, or $(Y_\lambda^{\mathrm{CIR}}(\infty), Z_\lambda(\infty))$. Because the stationary distributions of $U(\infty)$ and the scaled CIR intensity share the same first two moments, all three stationary pairs satisfy
\[
\E Z_\infty = 0, \qquad \E Z_\infty^2 = \frac{\sigma^2}{2\kappa}.
\]
Applying the stationary generator identity to $y$ and $yu$ then gives
\[
\E Y_\infty = -\frac{\beta_\lambda}{a}, \qquad
\E(Y_\infty Z_\infty) = \frac{\sigma^2}{2\kappa(a + \kappa)}.
\]
For the test function $y^2$, the second-order jump contribution to $\calG_X y^2$ in the exact scaled queue is $\delta_\lambda^2\{\iota + d_\lambda(x)\}$, which has stationary mean $2\lambda^{-\alpha}$ because $\E\iota = \E d_\lambda(X_\lambda(\infty)) = \lambda$, whereas the corresponding second-order contribution in either diffusion is $c_\lambda^2$. Applying the stationary generator identity to $y^2$ therefore yields
\begin{align*}
\Var(\widehat X_\lambda)
&= \frac{\lambda^{-\alpha}}{a} + \frac{\sigma^2}{2\kappa a(a + \kappa)}, \\
\Var\bigl(Y_\lambda^{\mathrm{OU}}(\infty)\bigr)
= \Var\bigl(Y_\lambda^{\mathrm{CIR}}(\infty)\bigr)
&= \frac{c_\lambda^2}{2a} + \frac{\sigma^2}{2\kappa a(a + \kappa)}.
\end{align*}
Thus the stationary means and mixed moments agree for every $c_\lambda$, and the stationary variances agree identically when $c_\lambda^2 = 2\lambda^{-\alpha}$. This establishes exact agreement on the basis $\{1, y, u, y^2, yu, u^2\}$ and hence on every polynomial of total degree at most two.
\Halmos
\endproof

\subsection{Intensity Path Estimates near the Lower Boundary}
\label{subsec:cir-stein-proofs}

The proof of Lemma~\ref{lem:cir_queue_stein} spans this subsection and the next two: here we establish the weighted pathwise moment bounds and construct the centered Poisson solution, Section~\ref{subsec:cir-control-proofs} proves the derivative bounds and weighted occupation-time estimates, and Section~\ref{subsec:cir-derivative-proofs} derives the second-derivative and queue-increment bounds.

\begin{lemma}
\label{lem:cir_queue_stein}
Suppose Assumptions~\ref{assump:model_primitives} and~\ref{assump:staffing_c} hold and $0 < c_\lambda \le \bar c$. Fix $m \ge 2$ and let $h \in C^2(\mathbb R^2)$ satisfy \eqref{eq:test_class} with constant $H$. There exist fixed exponents $p \ge m$ and $q > 0$, together with constants $C, \lambda_0 < \infty$ that may depend on $m$ and the fixed model parameters but are independent of $\lambda, h, H, c_\lambda, x, y, u$, such that for all $\lambda \ge \lambda_0$, the Poisson integral
\[
f(y, u) := -\int_0^\infty \bigl(P_t^{\mathrm{CIR}}h(y, u) - \pi^{\mathrm{CIR}}(h)\bigr)\,dt
\]
converges absolutely. Moreover, $f \in C^2(\mathsf D_\lambda)$ solves $\calG_{\mathrm{CIR}}f = h - \pi^{\mathrm{CIR}}(h)$ and satisfies
\begin{align}
|f_{yy}(y, u)|
&\le C H\left(1 + |\theta - \mu|c_\lambda^{-1/2}\right) w(y, u)^p\left[1 + v_\lambda(u)^{-q}\right],
\label{eq:cir-f11-state-dependent} \\
|f_{yu}(y, u)| + |f_{uu}(y, u)|
&\le C H w(y, u)^p\left[1 + v_\lambda(u)^{-q}\right],
\label{eq:cir-mixed-hessian} \\
|f_{yy}(y, u) - f_{yy}(x, u)|
&\le C H |y - x|(1 + |x| + |y| + |u|)^p\left[1 + v_\lambda(u)^{-q}\right]
\begin{cases}
c_\lambda^{-5/2}, & \mu \ne \theta, \\
1 + c_\lambda^{-1/2}, & \mu = \theta
\end{cases},
\label{eq:cir-f11-modulus-state-dependent}
\end{align}
where $x, y \in \R$ and $u > -\eps_\lambda^{-1}$.
\end{lemma}

\proof{Proof of Lemma~\ref{lem:cir_queue_stein}.}
Fix $\lambda$ and abbreviate $c = c_\lambda$, $\beta = \beta_\lambda$, $v(u) = v_\lambda(u) = 1 + \eps_\lambda u$, $v_0 = v(u)$, and $\nu = \sigma\eps_\lambda$. Let $g_\eta$ and $h_\eta$ ($0 < \eta \le 1$) denote the smooth approximations constructed in the proof of Lemma~\ref{lemma:stein-factor}, so $g_\eta$ satisfies \eqref{eq:stein-completion-smoothing}, $h_\eta$ satisfies \eqref{eq:test_class} with constant $C H$ uniformly in $\eta$, $|h_\eta - h| \le C H \eta w^{m-1}$, and $h_\eta \to h$ in $C^2_{\mathrm{loc}}$. To accommodate the resolvent argument used later for $(f_\eta)_{uu}$, we introduce a shifted generator indexed by $j \in \{0, 2\}$ and write
\[
(Y_t, Z_t) := \bigl(Y_\lambda^{\mathrm{CIR}, \eta, [j]}(t), Z_\lambda^{[j]}(t)\bigr),
\qquad (Y_0, Z_0) = (y, u),
\]
where $j = 0$ corresponds to the scaled arrival intensity $Z_\lambda$ in \eqref{eq:cir_diffusion_general_c}. The corresponding smoothed, shifted generator is
\[
\mathcal L_\eta^{[j]}
= [-\beta + u - g_\eta(y)]\partial_y
+ \left(-\kappa u + \frac{j}{2}\sigma^2\eps_\lambda\right)\partial_u
+ \frac{c^2}{2}\partial_{yy}
+ \frac{\sigma^2 v(u)}{2}\partial_{uu},
\]
with transition semigroup $P_t^{\eta, [j]}$; thus $j = 0$ is the smoothed CIR operator and $j = 2$ adds the drift shift $\sigma^2\eps_\lambda \partial_u$. Under either index $j$, the normalized intensity $V_t = 1 + \eps_\lambda Z_t$ satisfies
\[
dV_t = (\kappa_j - \kappa V_t)\,dt + \nu\sqrt{V_t}\,dB_2(t),
\qquad \kappa_j = \kappa + \frac{j}{2}\nu^2,
\]
whose one-dimensional generator in the $V$ coordinate is $\mathcal A^{[j]} = (\kappa_j - \kappa v)\partial_v + \frac{\nu^2 v}{2}\partial_{vv}$.

Throughout the proof, constants denoted by $C$ may increase from line to line and depend only on $m$, the fixed model parameters, and the fixed moment orders in use, remaining independent of $\eta, c, h, H$, the state, and all sufficiently large $\lambda$. We allow the spatial and boundary weight exponents to be enlarged finitely many times as functions of $m$, choosing the final threshold $\lambda_0$ large enough to satisfy all resulting inverse-drift and moment conditions.

We begin with the weighted transient and pathwise moment bounds. Define
\[
\Phi_Q(v) = 1 + v^{-Q}, \qquad W_{P, Q}(y, u) = w(y, u)^P \Phi_Q(v(u)).
\]
For each fixed $Q > 0$, a direct calculation gives
\begin{equation}
\mathcal A^{[j]} v^{-Q}
= Q\kappa v^{-Q} - Q\left(\kappa_j - \frac{Q + 1}{2}\nu^2\right)v^{-Q-1}.
\label{eq:cir-proof-inverse-drift}
\end{equation}
Because $\nu = \sigma\eps_\lambda \to 0$, the parenthesized coefficient is at least $\kappa/2$ for all sufficiently large $\lambda$, which yields $\mathcal A^{[j]} v^{-Q} \le -b_Q v^{-Q} + C_Q$ for fixed constants $b_Q, C_Q > 0$. Applying localized It\^o's formula and Fatou's lemma then gives
\[
\E_{v_0} V_t^{-Q}
\le e^{-b_Q t}v_0^{-Q} + \frac{C_Q}{b_Q}(1 - e^{-b_Q t})
\le C\Phi_Q(v_0).
\]
Meanwhile, since $|g_\eta(0)| \le a_1$ and the additional intensity drift in $\mathcal L_\eta^{[j]}$ is uniformly bounded, the higher-power Lyapunov calculation in Lemma~\ref{lem:diffusion_drift} applies uniformly across $\eta$ and $j \in \{0, 2\}$. Combining its transient polynomial moment bounds with the preceding inverse-moment estimate via Cauchy--Schwarz at orders $2P$ and $2Q$ yields
\begin{equation}
\sup_{t \ge 0}\E_{y, u}\bigl[w(Y_t, Z_t)^P \Phi_Q(V_t)\bigr]
\le C W_{P, Q}(y, u), \qquad P, Q \ge 0.
\label{eq:cir-proof-transient-weights}
\end{equation}

To bound running path suprema, fix $R > Q > 0$; for sufficiently large $\lambda$, \eqref{eq:cir-proof-inverse-drift} implies that $e^{-R\kappa s}V_s^{-R}$ is a nonnegative supermartingale, whose weak maximal inequality at exponent $Q/R < 1$ gives
\[
\E_{v_0}\sup_{0 \le s \le 1} V_s^{-Q}
\le \frac{R}{R - Q}e^{Q\kappa}v_0^{-Q}.
\]
Writing $S_t = (Y_t, Z_t)$ and letting $(\mathcal F_t)$ denote its natural filtration, the uniform linear growth of the SDE coefficients together with stopped Burkholder--Davis--Gundy and Gronwall estimates gives, for every fixed integer $p \ge 2$,
\begin{equation*}
\E_{y, u}\sup_{0 \le s \le 1} w(S_s)^p \le C w(y, u)^p.
\end{equation*}
Given any fixed moment order $b > 0$, choose an integer $P \ge \max\{2, b\}$ and $Q \ge b$, and set
\[
T_k = 1 + \sup_{k \le s \le k+1}(|Y_s| + |Z_s|) + \sup_{k \le s \le k+1} V_s^{-1},
\qquad k = 0, 1, \ldots.
\]
The two unit-interval supremum bounds imply
\begin{equation*}
\E_{y, u} T_0^b \le C\{w(y, u)^P + \Phi_Q(v_0)\} \le C W_{P, Q}(y, u),
\end{equation*}
and conditioning on $\mathcal F_k$ via the Markov property gives $\E[T_k^b \mid \mathcal F_k] \le C W_{P, Q}(S_k)$, so \eqref{eq:cir-proof-transient-weights} yields
\begin{equation*}
\E_{y, u} T_k^b \le C \E_{y, u} W_{P, Q}(S_k) \le C W_{P, Q}(y, u), \qquad k \ge 0.
\end{equation*}
Covering $[0, t]$ by $N = \lfloor t \rfloor + 1 \le 1 + t$ unit intervals gives
\[
1 + \sup_{s \le t}(|Y_s| + |Z_s|) + \sup_{s \le t} V_s^{-1}
\le 3 \max_{0 \le k < N} T_k,
\]
and bounding $\max_{k < N} T_k^b$ by $\sum_{k=0}^{N-1} T_k^b$ yields
\begin{equation}
\E_{y, u}\left[\left(1 + \sup_{s \le t}(|Y_s| + |Z_s|) + \sup_{s \le t} V_s^{-1}\right)^b\right]
\le C(1 + t)W_{P, Q}(y, u), \qquad t \ge 0.
\label{eq:cir-proof-path-moments}
\end{equation}

For the unshifted arrival intensity ($j = 0$), set $R_t = 2\sqrt{V_t}/(\sigma\eps_\lambda)$ and recall $k_\lambda = 2\kappa/(\sigma^2\eps_\lambda^2)$. By It\^o's formula, $R_t$ satisfies the Lamperti equation
\begin{equation}
\dd R_t = \left[\frac{k_\lambda - \tfrac{1}{2}}{R_t} - \frac{\kappa}{2}R_t\right]\dd t + \dd B_2(t).
\label{eq:cir-lamperti}
\end{equation}

We now establish convergence of the Poisson integral at $j = 0$ via an $L^1$ synchronous coupling. Write $z = (y, u)$, $z' = (y', u')$, and $S_t^z = (Y_t^{y, u}, Z_t^u)$. By \eqref{eq:cir-lamperti}, two synchronously coupled Lamperti paths satisfy
\[
\frac{d}{dt}(R_t^u - R_t^{u'})
= -\left(\frac{\kappa}{2} + \frac{k_\lambda - \tfrac{1}{2}}{R_t^u R_t^{u'}}\right)(R_t^u - R_t^{u'}),
\]
so $R_t^u - R_t^{u'}$ (and hence $Z_t^u - Z_t^{u'}$) preserves its initial sign for all $t \ge 0$. Combining sign preservation with the linear drift of $Z_t$ gives
\begin{equation*}
\E|Z_t^u - Z_t^{u'}| = \left|\E(Z_t^u - Z_t^{u'})\right| = e^{-\kappa t}|u - u'|.
\end{equation*}
Applying the secant-slope queue comparison from Lemma~\ref{lemma:stein-factor} then yields
\[
\E|Y_t^{y, u} - Y_t^{y', u'}|
\le e^{-a_0 t}|y - y'| + |u - u'|\int_0^t e^{-a_0(t-s)}e^{-\kappa s}\,ds,
\]
and hence, with $\gamma_0 = \tfrac{1}{4}\min\{a_0, \kappa\} > 0$,
\[
\E|S_t^z - S_t^{z'}| \le C e^{-2\gamma_0 t}|z - z'|.
\]

Let $\pi_\eta$ denote the stationary distribution at $j = 0$, couple $S_t^z$ synchronously with a stationary trajectory started at $\Xi \sim \pi_\eta$, and put
\[
D_t = |S_t^z - S_t^\Xi|, \qquad W_t = w(S_t^z) + w(S_t^\Xi).
\]
The preceding $L^1$ contraction, \eqref{eq:cir-proof-transient-weights}, and the stationary polynomial moment bounds in Lemma~\ref{lem:diffusion_drift} imply
\[
\E D_t \le C e^{-2\gamma_0 t}w(z), \qquad \E W_t^{2m-1} \le C w(z)^{2m-1}.
\]
Because $D_t \le W_t$, the derivative bound on $h_\eta$ and Cauchy--Schwarz yield
\begin{equation}
\begin{aligned}
|P_t^{\eta, [0]}h_\eta(z) - \pi_\eta(h_\eta)|
&\le C H \E[D_t W_t^{m-1}] \\
&\le C H(\E D_t)^{1/2}(\E[D_t W_t^{2m-2}])^{1/2} \\
&\le C H(\E D_t)^{1/2}(\E W_t^{2m-1})^{1/2} \\
&\le C H e^{-\gamma_0 t}w(z)^m.
\end{aligned}
\label{eq:cir-proof-poisson-convergence}
\end{equation}
Consequently, the centered Poisson integral
\[
f_\eta(y, u) = -\int_0^\infty [P_t^{\eta, [0]}h_\eta(y, u) - \pi_\eta(h_\eta)]\,dt
\]
converges absolutely, and the finite-time semigroup argument in Lemma~\ref{lemma:stein-factor} together with \eqref{eq:cir-proof-poisson-convergence} shows that $\mathcal L_\eta^{[0]}f_\eta = h_\eta - \pi_\eta(h_\eta)$ in distributions, so interior elliptic regularity on $\mathsf D_\lambda$ makes $f_\eta$ smooth for each fixed $\lambda, c, \eta > 0$. Finally, since the queue tangent satisfies
\[
A_t := \partial_y Y_t = \exp\left\{-\int_0^t g_\eta'(Y_s)\,ds\right\},
\qquad e^{-a_1 t} \le A_t \le e^{-a_0 t},
\]
the derivative bound on $h_\eta$ and the transient moment estimate justify differentiating under the Poisson integral to obtain
\[
|(f_\eta)_y(y, u)|
\le C H \int_0^\infty e^{-a_0 t}\E_{y, u}w(S_t)^{m-1}\,dt
\le C H w(y, u)^{m-1}.
\]

\subsection{Brownian Perturbations and Weighted Time near Full Occupancy}
\label{subsec:cir-control-proofs}

We first verify Malliavin differentiability with respect to the driving Brownian motions. For either $j \in \{0, 2\}$ and $0 \le r \le s$, the Lamperti coordinate in \eqref{eq:cir-lamperti} yields
\[
D_r^{(2)}Z_s
= \sigma\sqrt{V_s}\exp\left\{
-\frac{\kappa}{2}(s - r)
- \left(\frac{\kappa_j}{2} - \frac{\nu^2}{8}\right)\int_r^s V_a^{-1}\,da
\right\},
\qquad D_r^{(1)}Z_s = 0.
\]
To justify this formula rigorously, approximate the singular reciprocal term in the Lamperti drift from below by smooth, nonnegative, decreasing functions that agree with it outside a shrinking neighborhood of zero. Each approximating drift has derivative at most $-\kappa/2$, and pathwise comparison traps the approximating process between the linear-drift process with drift $-\kappa R/2$ and the exact positive process $R$, keeping its moments uniformly bounded while nonattainment of zero gives pathwise convergence on bounded intervals. Because the approximating Malliavin derivatives are bounded by $e^{-\kappa(s-r)/2}$, uniform integrability and closedness of the Malliavin derivative yield the displayed formula for $D_r^{(2)}Z_s$. Differentiating the smoothed queue equation in the same manner gives
\[
D_r^{(1)}Y_s = c\exp\left\{-\int_r^s g_\eta'(Y_a)\,da\right\},
\qquad
D_r^{(2)}Y_s = \int_r^s \exp\left\{-\int_a^s g_\eta'(Y_b)\,db\right\} D_r^{(2)}Z_a\,da,
\]
which, combined with \eqref{eq:cir-proof-path-moments}, ensures that $(Y_s, Z_s)$ possesses Malliavin derivatives with finite moments of every fixed order required below for sufficiently large $\lambda$.

Fix $0 < t \le 1$, set $r_s = t - s$, and write $a_s = g_\eta'(Y_s)$ and $\xi_s = \nu/(2\sqrt{V_s})$. For any predictable Brownian control $b = (b_1, b_2)$, define the directional Malliavin variations
\[
D_b F := \sum_{i=1}^2 \int_0^t D_s^{(i)}F\,b_i(s)\,ds,
\qquad \mathcal V_s^y := D_b Y_s, \quad \mathcal V_s^u := D_b Z_s.
\]
Whenever $F \in \mathbb D^{1, 2}$ and $\E\int_0^t |b(s)|^2\,ds < \infty$, the Skorohod divergence reduces to the It\^o integral $\mathcal I_b = \int_0^t b(s)\cdot\dd B(s)$ and satisfies $\E D_b F = \E[F\mathcal I_b]$ \citep{nualart2006malliavin}. Differentiating the coupled SDEs shows that $(\mathcal V_s^y, \mathcal V_s^u)$ evolves according to
\begin{align*}
d\mathcal V_s^y &= (-a_s\mathcal V_s^y + \mathcal V_s^u + c b_1(s))\,ds, \\
d\mathcal V_s^u &= -\kappa\mathcal V_s^u\,ds + \xi_s\mathcal V_s^u\,dB_2(s) + \sigma\sqrt{V_s}b_2(s)\,ds,
\qquad \mathcal V_0 = 0.
\end{align*}
Because of the multiplicative martingale term $\xi_s\mathcal V_s^u\,dB_2(s)$ in the second equation, the deterministic open-loop intensity perturbation used in the OU proof cannot hit the prescribed terminal state variation $(A_t, 0)$.

To steer $(\mathcal V_s^y, \mathcal V_s^u)$ via adapted feedback, consider the fixed matrix
\[
B_0 = \begin{pmatrix} 0 & -1 \\ 72 & -17 \end{pmatrix},
\]
which has eigenvalues $-8$ and $-9$, and choose the unique symmetric positive-definite matrix $P_0$ solving $(B_0 + 4I)^\top P_0 + P_0(B_0 + 4I) = -I$, so $B_0^\top P_0 + P_0 B_0 \le -8 P_0$. On $s \in [0, t)$, let $(\chi_s, \omega_s)$ solve the adapted linear SDE
\begin{equation*}
d\chi_s = -\frac{\omega_s}{r_s}\,ds, \qquad
d\omega_s = \frac{72\chi_s - 17\omega_s}{r_s}\,ds + \xi_s\omega_s\,dB_2(s),
\qquad (\chi_0, \omega_0) = (1, 0),
\end{equation*}
and define the feedback control
\begin{equation}
b_1(s) = 0, \qquad
b_2(s) = \frac{A_s}{\sigma\sqrt{V_s}}
\left\{\frac{72\chi_s - 16\omega_s}{r_s^2} + \frac{(\kappa - a_s)\omega_s}{r_s}\right\}.
\label{eq:cir-proof-feedback-control}
\end{equation}
Direct It\^o differentiation verifies that
\begin{equation*}
\mathcal V_s^y = A_s(1 - \chi_s), \qquad \mathcal V_s^u = A_s\omega_s / r_s,
\end{equation*}
for both $j = 0$ and $j = 2$, since the constant drift shift in $\mathcal L_\eta^{[2]}$ has zero pathwise derivative.

Set $F_s = (\chi_s, \omega_s)P_0(\chi_s, \omega_s)^\top$. For every fixed integer $\ell \ge 1$, It\^o's formula and $B_0^\top P_0 + P_0 B_0 \le -8 P_0$ give
\[
d(F_s^\ell)
\le \left(-\frac{8\ell}{r_s} + C_\ell\frac{\nu^2}{V_s}\right)F_s^\ell\,ds
+ F_s^\ell\zeta_s\frac{\nu}{\sqrt{V_s}}\,dB_2(s),
\qquad |\zeta_s| \le C_\ell.
\]
Fix $L \ge \max\{\ell, 1\}$ and apply the It\^o product rule to $X_s := F_s^\ell\Phi_L(V_s)$. All positive singular drift and covariation terms are bounded by $C_{\ell, L}\nu^2 F_s^\ell\{V_s^{-1} + V_s^{-L-1}\}$, whereas by \eqref{eq:cir-proof-inverse-drift} the inverse-power drift of $\Phi_L(V_s)$ contributes $-L(\kappa_j - (L+1)\nu^2/2)F_s^\ell V_s^{-L-1}$. Since $V_s^{-1} \le 1 + V_s^{-L-1}$ and $\kappa_j \ge \kappa$, this negative drift absorbs every singular $V_s^{-1}$ and $V_s^{-L-1}$ term once $\lambda$ is sufficiently large, leaving
\[
dX_s \le \left(C_{\ell, L} - \frac{8\ell}{r_s}\right)X_s\,ds + dM_s,
\]
where $M$ is a local martingale. Consequently,
\[
e^{-C_{\ell, L}s}\left(\frac{t}{r_s}\right)^{8\ell} X_s, \qquad 0 \le s < t,
\]
is a nonnegative local supermartingale; taking expectations after localization and applying Fatou's lemma (using $t \le 1$ and the deterministic initial value $F_0 = (1, 0)P_0(1, 0)^\top$) yields
\begin{equation*}
\E[F_s^\ell\Phi_L(V_s)]
\le C\left(\frac{r_s}{t}\right)^{8\ell}\Phi_L(v_0),
\qquad 0 \le s < t \le 1.
\end{equation*}
Because $r_s \le 1$, $A_s \le 1$, and $L \ge \ell$, substituting this decay estimate into \eqref{eq:cir-proof-feedback-control} gives
\[
\E|b_2(s)|^{2\ell} \le C r_s^{-4\ell}(r_s / t)^{8\ell}\Phi_L(v_0),
\]
and Minkowski's integral inequality yields the control-cost bound
\begin{equation}
\E\left(\int_0^t |b_2(s)|^2\,ds\right)^\ell \le C t^{-3\ell}\Phi_L(v_0).
\label{eq:cir-proof-control-energy}
\end{equation}

To verify the terminal condition at $s = t$, applying It\^o's formula to $\log F_s$ gives
\[
d\log F_s \le -8 r_s^{-1}\,ds + C\xi_s^2\,ds + H_s\,dB_2(s),
\qquad |H_s| \le C|\xi_s|.
\]
Because $V$ is almost surely bounded away from zero on the compact interval $[0, t]$, the drift and martingale integrals involving $\xi_s$ converge to finite limits as $s \uparrow t$, forcing $F_s = O(r_s^8)$ almost surely. Positive definiteness of $P_0$ implies $F_s \ge c_0(\chi_s^2 + \omega_s^2)$ for some $c_0 > 0$, whence $\chi_s \to 0$ and $\omega_s / r_s \to 0$ almost surely as $s \uparrow t$, confirming that
\[
\mathcal V_t = (A_t, 0) = \partial_y(Y_t, Z_t).
\]

We also record a pathwise maximal bound for $(\mathcal V_s^y, \mathcal V_s^u)$ needed in the occupation estimate. Since
\[
M_s^{\ell, L} = e^{-C_{\ell, L}s}(t / r_s)^{8\ell}F_s^\ell\Phi_L(V_s), \qquad s < t,
\]
is a nonnegative supermartingale, choosing $2\ell > b$ for any fixed $b > 0$ and applying its weak maximal inequality at exponent $b/(2\ell) < 1$ (together with $\frac{|\omega_s|}{r_s} \le C t^{-1}\bigl(M_s^{\ell, L}\bigr)^{1/(2\ell)}(r_s/t)^3$ up to a bounded exponential factor) proves
\begin{equation}
\left\|\sup_{s < t}|\mathcal V_s^u|\right\|_{L^b} \le C t^{-1}\Phi_L(v_0)^{1/(2\ell)},
\qquad
\left\|\sup_{s < t}|\mathcal V_s^y|\right\|_{L^b} \le C\Phi_L(v_0)^{1/(2\ell)}.
\label{eq:cir-proof-control-maximal}
\end{equation}

Malliavin integration by parts now yields, initially for smooth compactly supported $\psi$,
\begin{equation*}
\partial_y P_t^{\eta, [j]}\psi(y, u) = \E_{y, u}\left[\psi(Y_t, Z_t)\mathcal I_b\right],
\qquad \mathcal I_b = \int_0^t b(s)\cdot dB(s).
\end{equation*}
Pairing the feedback control \eqref{eq:cir-proof-feedback-control} with the direct queue-Brownian control
\[
b_1(s) = A_s / (c t), \qquad b_2(s) = 0,
\]
which satisfies $\mathcal V_s^y = (s/t)A_s$, $\mathcal V_s^u = 0$, and $\int_0^t |b(s)|^2\,ds \le (c^2 t)^{-1}$, and applying \eqref{eq:cir-proof-control-energy} and \eqref{eq:cir-proof-transient-weights} separately to each control gives
\begin{equation}
|\partial_y P_t^{\eta, [j]}\psi(y, u)|
\le C K_\psi W_{P', Q'}(y, u) k_c(t)
\quad\hbox{if } |\psi(y, u)| \le K_\psi W_{P, Q}(y, u),
\label{eq:cir-proof-semigroup-gradient}
\end{equation}
for fixed enlarged exponents $P', Q'$ and $k_c$ from \eqref{eq:ou-smoothing-kernel}: for $0 < t \le 1$ we take the smaller of the two control bounds, whereas for $t > 1$ we factor $P_t^{\eta, [j]} = P_1^{\eta, [j]} P_{t-1}^{\eta, [j]}$ and apply \eqref{eq:cir-proof-transient-weights}. Smooth truncation together with the transient weighted moment bounds at enlarged orders extends \eqref{eq:cir-proof-semigroup-gradient} to all continuous $\psi$ satisfying the stated growth condition, with $C, P', Q'$ independent of $c, \eta$, and sufficiently large $\lambda$.

At $j = 0$, the same pair of controls bounds the weighted occupation time near full occupancy. Put $\Theta_t = 1 + \sup_{s \le t}|Z_s| + \sup_{s \le t}V_s^{-1}$. For every fixed integer $\ell_0 \ge 0$, we show that there exist fixed exponents $P, Q$ such that, for every nonnegative $\varphi \in L^1(\mathbb R)$,
\begin{equation}
\E_{y, u}[\Theta_t^{\ell_0}\varphi(Y_t)]
\le C W_{P, Q}(y, u)\|\varphi\|_1
\begin{cases}
k_c(t), & 0 < t \le 1, \\
1 + t, & t > 1.
\end{cases}
\label{eq:cir-proof-weighted-occupation}
\end{equation}
First take $0 < t \le 1$ and a nonnegative $\varphi \in C_c^\infty(\mathbb R)$, and set
\[
\Psi(x) = \int_{-\infty}^x \varphi(a)\,da,
\qquad 0 \le \Psi \le \|\varphi\|_1, \qquad \Psi' = \varphi.
\]
To approximate $\Theta_t$ smoothly in $\mathbb D^{1, 2}$, sample $[0, t]$ along finite grids of vanishing mesh, replace $|x|$ by $\sqrt{x^2 + n^{-2}}$ and $V_s^{-1}$ by $(V_s + n^{-1})^{-1}$, and convolve each finite-grid maximum with a nonnegative smooth kernel of radius $n^{-1}$, denoting the resulting functional by $\Theta_{t, n}$. Because the partial derivatives of a smoothed finite maximum are nonnegative and sum to one, the chain rule together with $D_b V_s = \eps_\lambda \mathcal V_s^u$ and $|(V_s + n^{-1})^{-2}D_b V_s| \le V_s^{-2}|\mathcal V_s^u|$ gives
\[
\Theta_{t, n} \le C\Theta_t, \qquad
|D_b\Theta_{t, n}^{\ell_0}|
\le C_{\ell_0}\Theta_t^{\ell_0 + 1}\sup_{s \le t}|\mathcal V_s^u|.
\]
Path continuity and \eqref{eq:cir-proof-path-moments} ensure $\Theta_{t, n}^{\ell_0} \to \Theta_t^{\ell_0}$ in every required $L^b$ space while placing $\Psi(Y_t)\Theta_{t, n}^{\ell_0}$ in $\mathbb D^{1, 2}$. Applying Malliavin duality to $\Psi(Y_t)\Theta_{t, n}^{\ell_0}$ yields
\[
\E\!\left[\varphi(Y_t)A_t\Theta_{t, n}^{\ell_0}\right]
= \E\!\left[\Psi(Y_t)\Theta_{t, n}^{\ell_0}\mathcal I_b\right]
- \E\!\left[\Psi(Y_t)D_b\Theta_{t, n}^{\ell_0}\right].
\]
Under the arrival-intensity feedback control \eqref{eq:cir-proof-feedback-control}, \eqref{eq:cir-proof-control-energy} gives $\|\mathcal I_b\|_2 \le C t^{-3/2}\Phi_L(v_0)^{1/2}$; combining this bound with \eqref{eq:cir-proof-control-maximal} and \eqref{eq:cir-proof-path-moments} bounds the two terms on the right-hand side by
\[
\E[\varphi(Y_t)A_t\Theta_{t, n}^{\ell_0}]
\le C\|\varphi\|_1 W_{P, Q}(y, u)(t^{-3/2} + t^{-1})
\le C\|\varphi\|_1 W_{P, Q}(y, u) t^{-3/2}.
\]
Under the direct queue control $b_1(s) = A_s/(c t)$, $b_2(s) = 0$, the arrival-intensity path is unperturbed, so $D_b\Theta_{t, n}^{\ell_0} = 0$ and its energy bound gives
\[
\E[\varphi(Y_t)A_t\Theta_{t, n}^{\ell_0}]
\le C\|\varphi\|_1 W_{P, Q}(y, u) c^{-1}t^{-1/2}.
\]
Since $A_t \ge e^{-a_1 t}$ for $t \le 1$, letting $n \to \infty$ and taking the smaller of the two bounds yields
\[
\E_{y, u}[\Theta_t^{\ell_0}\varphi(Y_t)]
\le C\|\varphi\|_1 W_{P, Q}(y, u) k_c(t), \qquad 0 < t \le 1,
\]
which extends to every nonnegative $\varphi \in L^1(\mathbb R)$ by density. For $t > 1$, conditioning on $\mathcal F_{t-1}$, splitting the path supremum across $[0, t-1]$ and $[t-1, t]$, and applying the unit-interval estimate on $[t-1, t]$ together with \eqref{eq:cir-proof-path-moments} at an enlarged fixed order gives
\[
\E_{y, u}\!\left[\Theta_t^{\ell_0}\varphi(Y_t)\right]
\le C\|\varphi\|_1 \E_{y, u}\!\left[\Theta_{t-1}^{\ell_0} W_{P, Q}(Y_{t-1}, Z_{t-1})\right]
\le C(1 + t)W_{P', Q'}(y, u)\|\varphi\|_1,
\]
completing the proof of \eqref{eq:cir-proof-weighted-occupation}.

\subsection{Second Derivatives and Queue-Direction Increments}
\label{subsec:cir-derivative-proofs}

We now combine the pathwise flow variations with \eqref{eq:cir-proof-semigroup-gradient} and \eqref{eq:cir-proof-weighted-occupation} to bound the Hessian of $f_\eta$ and its queue-direction increments. Working at $j = 0$, write $b_0(r) = (k_\lambda - \tfrac{1}{2})/r - \kappa r/2$ for the Lamperti drift in \eqref{eq:cir-lamperti}, and put
\[
d = \frac{\kappa}{2} - \frac{\nu^2}{8} \ge 0, \qquad
L_t = \partial_{R_0}R_t, \qquad N_t = \partial_{R_0}^2 R_t.
\]
Differentiating the Lamperti flow with respect to $R_0$ gives
\begin{align*}
L_t &= \exp\left\{-\frac{\kappa}{2}t - d\int_0^t V_s^{-1}\,ds\right\} \le e^{-\kappa t/2}, \\
N_t &= \int_0^t \exp\left\{\int_s^t b_0'(R_r)\,dr\right\} b_0''(R_s) L_s^2\,ds,
\qquad b_0''(R_s) = \sigma\eps_\lambda d\,V_s^{-3/2},
\end{align*}
so in particular
\[
0 \le N_t \le \sigma\eps_\lambda d\,e^{-\kappa t/2}\int_0^t e^{-\kappa s/2}V_s^{-3/2}\,ds.
\]
Transforming back to $J_t = \partial_u Z_t$ and $K_t = \partial_u^2 Z_t$ via $Z_t = \eps_\lambda^{-1}(\tfrac{1}{4}\sigma^2\eps_\lambda^2 R_t^2 - 1)$ yields
\begin{equation}
J_t = \sqrt{V_t / v_0}\,L_t, \qquad
K_t = \frac{\eps_\lambda}{2v_0}L_t^2 + \frac{\sqrt{V_t}}{\sigma v_0}N_t - \frac{\eps_\lambda\sqrt{V_t}}{2v_0^{3/2}}L_t,
\label{eq:cir-proof-environment-flows}
\end{equation}
and hence
\[
|J_t| \le e^{-\kappa t/2}\Theta_t, \qquad |K_t| \le C e^{-\kappa t/2}\Theta_t^3,
\]
where $K_t$ captures the nonlinear second-order sensitivity of the CIR arrival intensity to its initial state.

Setting $A(s, t) = \exp\{-\int_s^t g_\eta'(Y_r)\,dr\}$ and retaining $\gamma_0 = \tfrac{1}{4}\min\{a_0, \kappa\}$, the first-order queue variations satisfy
\[
j_y(t) = A(0, t), \qquad j_u(t) = \int_0^t A(s, t)J_s\,ds,
\]
so that
\[
|j_y(t)| \le e^{-2\gamma_0 t}, \qquad |j_u(t)| \le t e^{-2\gamma_0 t}\Theta_t.
\]
For $a, b \in \{y, u\}$, put $m_{ab} = \mathbf{1}_{\{a=u\}} + \mathbf{1}_{\{b=u\}}$. Differentiating the smoothed queue equation once more gives
\begin{equation}
\begin{aligned}
\partial_{ab} Y_t
&= -(\theta - \mu)\int_0^t A(s, t)\rho_\eta(Y_s) j_a(s)j_b(s)\,ds
+ \mathbf{1}_{\{a=b=u\}}\int_0^t A(s, t)K_s\,ds, \\
\partial_{ab} Z_t &= \mathbf{1}_{\{a=b=u\}} K_t.
\end{aligned}
\label{eq:cir-proof-queue-second-flows}
\end{equation}

Applying the chain rule in \eqref{eq:stein-completion-semigroup-hessian} to $S_t = (Y_t, Z_t)$ and $h_\eta$ (now retaining the additional intensity-curvature term $(h_\eta)_u(S_t)\partial_{ab}Z_t$) is justified for each fixed $\eta$ by localization and the pathwise moment estimates in \eqref{eq:cir-proof-path-moments}. Using the first-variation bounds, the Hessian growth bound on $h_\eta$, and \eqref{eq:cir-proof-path-moments} gives
\[
\begin{aligned}
&\left|\E\left[(\partial_a S_t)^\top \nabla^2 h_\eta(S_t)\partial_b S_t\right]\right| \\
&\quad \le C H(1 + t)^2 e^{-4\gamma_0 t}\E[w(S_t)^{m-2}\Theta_t^2]
\le C H(1 + t)^3 e^{-4\gamma_0 t}W_{P, Q}(y, u),
\end{aligned}
\]
whose integral over $t \ge 0$ is bounded by $C H W_{P, Q}(y, u)$ uniformly in $c$. Likewise, conditioning on $\mathcal F_s$ (which replaces $|(h_\eta)_y(S_t)|$ with $C H w(S_s)^{m-1}$ for $t \ge s$ by the transient moment bound), integrating $A(s, t) \le e^{-a_0(t-s)}$ over $t \ge s$, and applying \eqref{eq:cir-proof-environment-flows} and \eqref{eq:cir-proof-path-moments} bounds the time integrals of both $K$-terms by
\[
C H \int_0^\infty e^{-\kappa s/2}\E[w(S_s)^{m-1}\Theta_s^3]\,ds \le C H W_{P, Q}(y, u).
\]
For the full-occupancy mollifier integral in \eqref{eq:cir-proof-queue-second-flows}, the same conditioning on $\mathcal F_s$, together with $|j_a(s)j_b(s)| \le s^{m_{ab}}e^{-4\gamma_0 s}\Theta_s^{m_{ab}}$ and $w(S_s) \le 2\Theta_s$ on the support of $\rho_\eta(Y_s)$, bounds the remaining integrand at time $s$ by
\[
C H s^{m_{ab}}e^{-4\gamma_0 s}\E[\Theta_s^{m - 1 + m_{ab}}\rho_\eta(Y_s)].
\]
Applying \eqref{eq:cir-proof-weighted-occupation} with any integer $\ell_0 \ge m - 1 + m_{ab}$ (using $\Theta_s \ge 1$ and $\|\rho_\eta\|_1 = 1$) and integrating $e^{-a_0(t-s)}$ over $t \ge s$ yields
\[
\int_0^\infty |(P_t^{\eta, [0]}h_\eta)_{ab}(y, u)|\,dt
\le C H W_{P, Q}(y, u)\left[1 + |\theta - \mu|\int_0^1 s^{m_{ab}}k_c(s)\,ds\right],
\]
since the long-time integral over $s > 1$ is bounded uniformly in $c$ by the integrability of $s^{m_{ab}}(1 + s)e^{-4\gamma_0 s}$. Invoking the short-time kernel integrals $\int_0^1 k_c(s)\,ds \le C c^{-1/2}$ and $\int_0^1 s^r k_c(s)\,ds \le C$ ($r = 1, 2$) from \eqref{eq:stein-completion-Ic0}--\eqref{eq:stein-completion-scaling} establishes \eqref{eq:cir-f11-state-dependent}--\eqref{eq:cir-mixed-hessian} for $f_\eta$ uniformly in $\eta$.

To bound the queue-direction increments of $(f_\eta)_{yy}$, we adapt the resolvent argument of Lemma~\ref{lem:ou-fyy-lipschitz} using the shifted operator $\mathcal L_\eta^{[2]}$. Differentiating $\mathcal L_\eta^{[0]}f_\eta = h_\eta - \pi_\eta(h_\eta)$ twice in $u$ gives
\[
(\mathcal L_\eta^{[2]} - 2\kappa)(f_\eta)_{uu} = F_\eta := (h_\eta)_{uu} - 2(f_\eta)_{yu},
\]
where $|F_\eta| \le C H W_{P, Q}$ after enlarging $P$. Applying discounted Dynkin's formula up to the first exit from a compact subset of $\mathsf D_\lambda$ and a fixed horizon $T$, the bounds on $(f_\eta)_{uu}$ and $F_\eta$ together with \eqref{eq:cir-proof-path-moments} at enlarged orders bound the stopped terminal and drift terms in $L^2$ on $[0, T]$, allowing removal of the stopping times by uniform integrability; sending $T \to \infty$ via \eqref{eq:cir-proof-transient-weights} and the discount factor $e^{-2\kappa T}$ then gives
\begin{equation*}
(f_\eta)_{uu} = -\int_0^\infty e^{-2\kappa t}P_t^{\eta, [2]}F_\eta\,dt.
\end{equation*}
Differentiating in $y$ and applying \eqref{eq:cir-proof-semigroup-gradient} and \eqref{eq:stein-completion-Ic0} yields $|(f_\eta)_{yuu}| \le C H(1 + c^{-1/2})W_{P, Q}$. Writing $b_\eta = -\beta + u - g_\eta(y)$, the CIR counterpart of \eqref{eq:ou-lip-fyyy-equation} is
\begin{equation*}
\frac{c^2}{2}(f_\eta)_{yyy}
= (h_\eta)_y + g_\eta'(f_\eta)_y - b_\eta(f_\eta)_{yy} + \kappa u(f_\eta)_{yu} - \frac{\sigma^2 v(u)}{2}(f_\eta)_{yuu}.
\end{equation*}
As in \eqref{eq:ou-lip-smooth-fyyy}, using $|(h_\eta)_y| + |(f_\eta)_y| \le C H w^{m-1}$, $|b_\eta| \le C w$, $v(u) \le w$, and $0 < c \le \bar c$ enlarges only the fixed polynomial exponent $P$ and gives $|(f_\eta)_{yyy}| \le C H c^{-5/2}W_{P, Q}$. When $\mu = \theta = a$, the deterministic queue tangent instead gives
\[
(f_\eta)_{yy} = -\int_0^\infty e^{-2at}P_t^{\eta, [0]}[(h_\eta)_{yy}]\,dt,
\]
and applying \eqref{eq:cir-proof-semigroup-gradient} to $|(h_\eta)_{yy}| \le C H w^{m-2}$ yields $|(f_\eta)_{yyy}| \le C H(1 + c^{-1/2})W_{P, Q}$, matching the two branches of \eqref{eq:cir-f11-modulus-state-dependent}.

Finally, we pass to $\eta \downarrow 0$ by the limiting arguments in Lemmas~\ref{lemma:stein-factor} and~\ref{lem:ou-fyy-lipschitz}. Coupling $Y_\lambda^{\mathrm{CIR}, \eta, [0]}$ and $Y_\lambda^{\mathrm{CIR}}$ with the same arrival-intensity path and queue-length Brownian motion gives $\sup_{t \ge 0}|Y_\lambda^{\mathrm{CIR}, \eta, [0]}(t) - Y_\lambda^{\mathrm{CIR}}(t)| \le a_1\eta / a_0$. Combined with $|h_\eta - h| \le C H \eta w^{m-1}$, the derivative bound on $h$, and the uniform polynomial moments, this pathwise bound implies convergence of both the stationary means $\pi_\eta(h_\eta) \to \pi^{\mathrm{CIR}}(h)$ and the finite-horizon Poisson integrals locally uniformly in $(y, u)$, while \eqref{eq:cir-proof-poisson-convergence} controls the infinite-horizon tails uniformly in $\eta$ so that $f_\eta \to f$ locally uniformly on $\mathsf D_\lambda$. Because the weight $W_{P, Q}$ is bounded on every compact subset of $\mathsf D_\lambda$, integration by parts passes the derivative bounds to the weak derivatives of $f$, and passing to the limit in
\[
\calG_{\mathrm{CIR}}f_\eta = h_\eta - \pi_\eta(h_\eta) + (g_\eta - g)(f_\eta)_y
\]
shows that $\calG_{\mathrm{CIR}}f = h - \pi^{\mathrm{CIR}}(h)$ in distributions. Local Sobolev embedding and interior Schauder regularity for fixed $c > 0$ then yield $f \in C^2(\mathsf D_\lambda)$, while the weak-derivative slicing argument of Lemma~\ref{lem:ou-fyy-lipschitz} establishes \eqref{eq:cir-f11-modulus-state-dependent} for every $u > -\eps_\lambda^{-1}$. Choosing fixed exponents $p \ge m$ and $q > 0$ to dominate the finitely many spatial and boundary weights encountered above (and enlarging $\lambda_0$ to satisfy all required inverse-drift, moment, and control conditions without altering any power of $c$) completes the proof.
\Halmos
\endproof

\section{Laguerre--QBD Computation of the Exact Stationary Distribution}
\label{sec:ec-laguerre-qbd}

This section details the spectral quasi-birth-and-death (QBD) method used in Section~\ref{sec:numerical-experiments} to compute stationary expectations under $\pi^X$ without Monte Carlo simulation error.

\subsection{Generalized Laguerre Eigenbasis and Exact Moment Balance}
\label{subsec:ec-laguerre-basis}

Recall from Section~\ref{sec:model} and Lemma~\ref{lem:cir_inverse_moments} that the normalized arrival intensity $V_\lambda(t) := I_\lambda(t)/\lambda = 1 + \eps_\lambda \wt I_\lambda(t)$ takes values in $(0, \infty)$ under Assumption~\ref{assump:model_primitives} and has stationary $\operatorname{Gamma}(k_\lambda, 1/k_\lambda)$ density
\[
p_\lambda(v) = \frac{k_\lambda^{k_\lambda}}{\Gamma(k_\lambda)} v^{k_\lambda - 1} e^{-k_\lambda v},
\qquad v > 0,
\]
where $k_\lambda = 2\kappa/(\sigma^2\eps_\lambda^2) = (2\kappa/\sigma^2)\lambda^{1-\alpha} \ge 1$. On smooth functions $\phi : (0, \infty) \to \R$, the infinitesimal generator of $V_\lambda$ acts as
\begin{equation}
\mathcal{A}_\lambda \phi(v)
= \kappa(1 - v)\phi'(v) + \frac{\kappa}{k_\lambda} v \phi''(v)
= \frac{\kappa}{k_\lambda p_\lambda(v)} \frac{\dd}{\dd v}\bigl[v\,p_\lambda(v)\phi'(v)\bigr],
\label{eq:ec-cir-sturm-liouville}
\end{equation}
which is symmetric in $L^2\bigl((0, \infty), p_\lambda(v)\,\dd v\bigr)$ and admits a discrete polynomial eigenbasis \citep{FormanSorensen2008}. Specifically, let $L_j^{(k_\lambda - 1)}(z)$ ($j \ge 0$) denote the generalized Laguerre polynomials of index $k_\lambda - 1 \ge 0$, and define the orthonormal polynomials
\begin{equation}
\psi_j(v)
:= (-1)^j \sqrt{\frac{j!\,\Gamma(k_\lambda)}{\Gamma(j + k_\lambda)}}\, L_j^{(k_\lambda - 1)}(k_\lambda v),
\qquad j = 0, 1, 2, \ldots,
\label{eq:ec-orthonormal-laguerre}
\end{equation}
so that $\psi_0(v) \equiv 1$, $\psi_1(v) = \sqrt{k_\lambda}\,(v - 1)$, and $\int_0^\infty \psi_i(v)\psi_j(v)p_\lambda(v)\,\dd v = \delta_{ij}$. By Laguerre's differential equation and \eqref{eq:ec-cir-sturm-liouville}, each $\psi_j$ satisfies the eigenvalue relation $\mathcal{A}_\lambda \psi_j = -j\kappa\,\psi_j$. Moreover, the standard three-term recurrence for $L_j^{(k_\lambda - 1)}(k_\lambda v)$ expresses multiplication by $v$ in the basis $\{\psi_j\}_{j \ge 0}$ via a symmetric tridiagonal Jacobi operator $J_\lambda$:
\begin{equation}
v\,\psi_j(v) = b_{j-1}\psi_{j-1}(v) + a_j\psi_j(v) + b_j\psi_{j+1}(v),
\qquad j \ge 0 \quad (b_{-1} := 0),
\label{eq:ec-laguerre-jacobi}
\end{equation}
with strictly positive entries
\[
a_j = 1 + \frac{2j}{k_\lambda},
\qquad
b_j = \frac{\sqrt{(j + 1)(j + k_\lambda)}}{k_\lambda},
\qquad j \ge 0.
\]

Now consider the unscaled queue length paired with the normalized intensity, $\bigl(X_\lambda(t), V_\lambda(t)\bigr)$ on $\mathbb{Z}_+ \times (0, \infty)$, whose generator is
\[
\mathcal{L}_\lambda F(x, v)
= \lambda v\bigl[F(x + 1, v) - F(x, v)\bigr]
+ d_\lambda(x)\bigl[F(x - 1, v) - F(x, v)\bigr]
+ \mathcal{A}_\lambda F(x, v),
\]
with $d_\lambda(x)$ given by \eqref{eq:death_rate} and $d_\lambda(0) = 0$. Define the stationary Laguerre--Fourier coefficients
\[
\xi_{x, j}
:= \E\bigl[\mathbb{I}_{\{X_\lambda(\infty) = x\}}\,\psi_j\bigl(V_\lambda(\infty)\bigr)\bigr],
\qquad x \in \mathbb{Z}_+,\quad j \ge 0,
\]
and collect them into column vectors $\boldsymbol{\xi}_x = (\xi_{x, 0}, \xi_{x, 1}, \ldots)^\top$ (with $\boldsymbol{\xi}_{-1} := \mathbf{0}$). For each fixed $x \in \mathbb{Z}_+$ and $j \ge 0$, the test function $F(x', v) = \mathbb{I}_{\{x' = x\}}\psi_j(v)$ has finite support in $x'$ and polynomial growth of degree $j$ in $v$, so the stationary generator identity $\E[\mathcal{L}_\lambda F(X_\lambda(\infty), V_\lambda(\infty))] = 0$ holds by the localization argument of Lemma~\ref{lem:petite_generator_identity}. Substituting $\mathcal{A}_\lambda \psi_j = -j\kappa\psi_j$ and \eqref{eq:ec-laguerre-jacobi} yields the exact infinite-dimensional level balance equations
\begin{equation}
\lambda J_\lambda \boldsymbol{\xi}_{x-1}
- \bigl(\lambda J_\lambda + d_\lambda(x)I + \kappa D\bigr)\boldsymbol{\xi}_x
+ d_\lambda(x + 1)\boldsymbol{\xi}_{x+1}
= \mathbf{0},
\qquad x \in \mathbb{Z}_+,
\label{eq:ec-qbd-balance}
\end{equation}
where $D := \operatorname{diag}(0, 1, 2, \ldots)$. Summing \eqref{eq:ec-qbd-balance} over $x \in \mathbb{Z}_+$ telescopes the arrival and departure terms to give $\kappa D \sum_{x=0}^\infty \boldsymbol{\xi}_x = \mathbf{0}$, which enforces the marginal orthogonality conditions $\sum_{x=0}^\infty \xi_{x, j} = \E[\psi_j(V_\lambda(\infty))] = 0$ for all $j \ge 1$, together with the normalization $\sum_{x=0}^\infty \xi_{x, 0} = 1$.

\subsection{Spectral Truncation, Quadrature Interpretation, and Active Lattice Scaling}
\label{subsec:ec-qbd-truncation}

To solve \eqref{eq:ec-qbd-balance} numerically, we retain the first $M$ Laguerre modes ($j = 0, \ldots, M - 1$), replacing $J_\lambda$ and $D$ by their leading $M \times M$ principal submatrices $J_\lambda^{(M)}$ and $D^{(M)} = \operatorname{diag}(0, 1, \ldots, M - 1)$. This spectral truncation preserves Markovian positivity. By the Golub--Welsch theorem \citep{GolubWelsch1969}, the symmetric tridiagonal matrix $J_\lambda^{(M)}$ admits an orthogonal eigendecomposition $J_\lambda^{(M)} = Q \Lambda_M Q^\top$, where $\Lambda_M = \operatorname{diag}(v_0, \ldots, v_{M-1})$ contains the $M$ Gauss--Laguerre quadrature nodes in $(0, \infty)$ and $w_m = Q_{0, m}^2 > 0$ are the corresponding quadrature weights. Transforming $\boldsymbol{\xi}_x^{(M)}$ into quadrature-node coordinates diagonalizes the upward transition block $\lambda J_\lambda^{(M)}$ to $\lambda\Lambda_M$ while mapping $-\kappa Q^\top D^{(M)} Q$ to the generator of a reversible $M$-state continuous-time Markov environment whose stationary distribution $(w_0, \ldots, w_{M-1})$ matches the first $2M - 1$ moments of $V_\lambda(\infty)$. Consequently, the mode-truncated system is equivalent to the stationary equations of an irreducible $M$-phase QBD process \citep{LatoucheRamaswami1999} with level-dependent departure rates $d_\lambda(x)$, guaranteeing a unique, strictly positive stationary solution normalized by $\sum_{x=0}^\infty \xi_{x, 0}^{(M)} = 1$.

Once $\bigl\{\boldsymbol{\xi}_x^{(M)}\bigr\}$ is computed, marginal queue-length probabilities follow directly from $\mathbb{P}(X_\lambda(\infty) = x) \approx \xi_{x, 0}^{(M)}$ since $\psi_0 \equiv 1$, and joint queue--intensity moments follow from low-degree modes. For instance, $\wt I_\lambda(\infty) = (\sigma/\sqrt{2\kappa})\,\psi_1(V_\lambda(\infty))$ gives $\E[\wt X_\lambda(\infty)\wt I_\lambda(\infty)] \approx (\sigma/\sqrt{2\kappa})\sum_x \delta_\lambda(x - n_\lambda)\xi_{x, 1}^{(M)}$. In our experiments, we set $M \in [36, 52]$ (increasing $M$ when $k_\lambda$ is small as $\alpha \uparrow 1$) and restrict $x$ to a reflecting active window $[x_{\min}, x_{\max}] \subset \mathbb{Z}_+$ extending $14$ standard deviations beyond the quasi-stationary queue equilibria associated with the extreme quadrature nodes of $\Lambda_M$, so that the boundary probability mass $\xi_{x_{\min}, 0}^{(M)} + \xi_{x_{\max}, 0}^{(M)}$ remains below $10^{-14}$. Ordering the unknowns level-by-level produces a banded linear system of half-bandwidth $M$ whose active level count $L_\lambda = x_{\max} - x_{\min} + 1$ scales as $\Theta(\lambda^{(1+\alpha)/2})$, requiring $O(L_\lambda M^3) = O(M^3 \lambda^{(1+\alpha)/2})$ operations via banded LU factorization.

For the diffusion approximations $(Y_\lambda^{\mathrm{OU}}, U)$ and $(Y_\lambda^{\mathrm{CIR}}, Z_\lambda)$, an analogous orthogonal-polynomial representation in the intensity coordinate (using normalized probabilists' Hermite polynomials for $U$ and the same Laguerre basis $\{\psi_j\}$ in \eqref{eq:ec-orthonormal-laguerre} for $Z_\lambda$, mapped to Gauss--Hermite or Gauss--Laguerre quadrature nodes) reduces the two-dimensional stationary Fokker--Planck equation to $M$ coupled one-dimensional advection--diffusion equations in the scaled queue coordinate $y$. Discretizing $y$ on a bounded domain $[y_{\min}, y_{\max}]$ (with $y = 0$ placed at a cell interface) via the exponentially fitted finite-volume scheme of \citet{ChangCooper1970} yields a banded linear system of fixed dimension $N_y M$ independent of $\lambda$, which remains second-order accurate across $c_\lambda > 0$ and degenerates smoothly to exact upwind differencing at $c_\lambda = 0$.

\end{document}